\documentclass[draftcls,11pt,onecolumn]{IEEEtran}

\usepackage{amsmath}
\usepackage{amssymb}
\usepackage{amsthm}
\usepackage{cite}
\usepackage{graphicx}
\usepackage{verbatim}
\usepackage{subfigure}
\usepackage{color}

\usepackage{epsfig}

\usepackage{citesort}
\ifx\pdfoutput\undefined
\usepackage{graphicx}
\else
\fi

\usepackage{psfrag}
\usepackage{dsfont}
\usepackage{subfigure}

\usepackage{url}

\usepackage{stfloats}
\usepackage{amsmath,amsfonts,amssymb,amsxtra}
\usepackage{array}

\usepackage{pstricks,pst-node,pst-text,pst-3d,pst-plot,pst-all}

\newcommand{\tr}[1]{{\rm tr} \! \left(#1 \right)}
\newcommand{\be}{\begin{equation}}
\newcommand{\ee}{\end{equation}}

\newcommand{\ist}{\hspace*{.2mm}}
\newcommand{\rmv}{\hspace*{-.2mm}}
\newcommand{\remark}[1]{}
\renewcommand{\a}{{\bf a}}

\newcommand{\e}{{\bf e}}

\newcommand{\m}{{\bf m}}
\renewcommand{\u}{{\bf u}}
\renewcommand{\v}{{\bf v}}
\newcommand{\x}{{\bf x}}
\newcommand{\y}{{\bf y}}
\newcommand{\A}{{\bf A}}
\newcommand{\B}{{\bf B}}

\newcommand{\bE}{{\bf E}}
\newcommand{\I}{{\bf I}}
\newcommand{\J}{{\bf J}}
\newcommand{\bM}{{\bf M}}
\renewcommand{\P}{{\bf P}}
\newcommand{\Q}{{\bf Q}}
\newcommand{\V}{{\bf V}}
\newcommand{\hx}{{\hat{\x}}}

\DeclareMathOperator*{\argmin}{arg\,min}
\DeclareMathOperator*{\argmax}{arg\,max}

\newcommand{\normconstgauss}{\frac{1}{(2 \pi \sigma^{2})^{N/2}} }
\def\CR{Cram\'er--Rao }

\def\MSE{\varepsilon(\xtrue; \hat{\mathbf{x}})}
\def\LS_est{\hat{\mathbf{x}}_{\text{LS}}}
\def\ML_est{\hat{\mathbf{x}}_{\text{ML}}}
\newtheorem{theorem}{Theorem}

\newtheorem{lemma}[theorem]{Lemma}

\newtheorem{proposition}[theorem]{Proposition}

\newcommand{\eq}{\,=\,}
\renewcommand{\jmath}{j}

\newcommand{\pinv}{\dagger}
\newcommand{\XS}{{\mathcal{X}_S}}
\newcommand{\RR}{{\mathbb R}}
\newcommand{\one}{{\boldsymbol 1}}
\newcommand{\zero}{{\boldsymbol 0}}
\newcommand{\bdel}{{\boldsymbol\delta}}
\newcommand{\expp}[1]{\exp\!\left(#1\right)}
\newcommand{\eps}{\varepsilon}

\newcommand{\beq}{\begin{equation}}
\newcommand{\eeq}{\end{equation}}

\DeclareMathOperator{\supp}{supp}

\newcommand{\xtrue}{\mathbf{x}}
\newcommand{\xitrue}{\xi}

\allowdisplaybreaks

\begin{document}
%%%%%%%%%%%%%%%%%%%%%%%%%%%%%%%%%%%%%%%%%%%%%%%%%%%%%%%%%%

\title{Unbiased Estimation of a Sparse Vector in\\[-3mm]White Gaussian Noise\vspace*{4mm}}
\author{\emph{ Alexander Jung\textsuperscript{a}, Zvika Ben-Haim\textsuperscript{b}, Franz Hlawatsch\textsuperscript{a}, and Yonina C.\ Eldar\textsuperscript{b} 
\thanks{This work was supported by the FWF under Grant S10603-N13 (Statistical Inference) within the National Research Network SISE, by the WWTF under Grant MA 07-004 (SPORTS), by the Israel Science Foundation under Grant 1081/07, and by the European Commission under the FP7 Network of Excellence in Wireless Communications NEWCOM++ (contract no. 216715).
Parts of this work were previously presented at ICASSP 2010, Dallas, TX, March 2010.}
} \\[2mm]
{\normalsize\emph{\textsuperscript{a}}Institute of Communications and Radio-Frequency Engineering, Vienna University of Technology}\\[-2.3mm]
{\normalsize Gusshausstrasse\ 25/389, 1040 Vienna, Austria}\\[-2.3mm]
{\normalsize Phone: +43 1 58801 38963, Fax: +43 1 58801 38999, E-mail: \{ajung,\,fhlawats\}@nt.tuwien.ac.at}\\[.5mm]
{\normalsize\emph{\textsuperscript{b}}Technion---Israel Institute of Technology}\\[-2.3mm]
{\normalsize Haifa 32000, Israel; e-mail: \{zvikabh@tx,\,yonina@ee\}.technion.ac.il} %% \\[-2mm]
}

%\address{\normalsize \textsuperscript{a}Institute of Communications and Radio-Frequency Engineering, Vienna University of Technology\\[-1.2mm]
%\normalsize Gusshausstrasse 25/389, A-1040 Vienna, Austria; e-mail: ajung@nt.tuwien.ac.atÊ\\[1mm]
%\normalsize \textsuperscript{b}Technion---Israel Institute of Technology, %\\[-1.2mm]
%\normalsize Haifa 32000, Israel; e-mail: zvikabh@tx.technion.ac.il}

 %\date{Submitted to IEEE Transactions on Information Theory, \today}

%%%%%%%%%%%%%%%%%%%%%%%%%%%%%%%%%%%%%%%%%%%%%%%%%%%%%%%%%%
\maketitle
%%%%%%%%%%%%%%%%%%%%%%%%%%%%%%%%%%%%%%%%%%%%%%%%%%%%%%%%%%

\vspace*{-10mm}

\begin{abstract}
\vspace*{-1.5mm}
We consider unbiased estimation of a sparse nonrandom vector corrupted by additive white Gaussian noise.
We show that while there are infinitely many unbiased estimators for this problem, 
none of them has uniformly minimum variance. Therefore, we focus on locally minimum variance unbiased %%% AJ
(LMVU) estimators. We derive simple closed-form lower and upper bounds on the variance of LMVU estimators or, equivalently,
on the Barankin bound (BB). Our bounds allow an estimation
%% identification 
of the threshold region separating the low-SNR and high-SNR regimes, and they indicate the asymptotic behavior of the BB at high SNR.
We also develop numerical lower and upper bounds which are tighter than the closed-form bounds and thus characterize the BB more %%% AJ
accurately. 
Numerical studies compare our characterization of the BB with established 
%% though 
biased estimation schemes, and demonstrate that while unbiased estimators perform poorly at low SNR, they may perform better than biased estimators at high SNR. 
An interesting conclusion of our analysis is that the high-SNR behavior of the BB depends solely on the value 
of the smallest nonzero component of the sparse vector, and that this type of dependence is also exhibited by
%% also appears to be true for 
the performance of certain practical estimators.
\vspace*{-1mm}
\end{abstract}

\begin{keywords}
\vspace*{-1.5mm}
Sparsity, unbiased estimation, denoising, Cram\'{e}r--Rao bound, Barankin bound, Hammersley--Chapman--Robbins 
bound, locally minimum variance unbiased estimator.
%% Parameter estimation
\vspace*{2mm}
\end{keywords}

%ALEX: Suggestion: I think it would be a good idea to order the appendices in the order in which they are referenced in the text. Or is there some logic in the current order?

%%%%%%%%%%%%%%%%%%%%%%%%%%%%%%%%%%%%%%%%%%%%%%%%%%%%%%%%%%
\section{{Introduction}}\label{sec.intro}
%%%%%%%%%%%%%%%%%%%%%%%%%%%%%%%%%%%%%%%%%%%%%%%%%%%%%%%%%%

\vspace{1mm}

Research in the past few years has led to a recognition that the performance of signal processing algorithms can be boosted by exploiting the tendency of 
many
signals to have sparse representations.
%Perhaps the most celebrated application of this principle 
%% of this fact 
%is 
%% in the field of 
%, where sparsity is used for the purpose of signal compression. %%% AJ 
Applications of this principle include signal reconstruction (e.g. in the context of compressed sensing \cite{Don06,NearOptimalSigRecCanTao06}) and signal enhancement (e.g. in the context of image denoising and deblurring \cite{Donoho94idealspatial,dabov08,protter09}). %%% AJ
%In addition, there are other studies which investigated how one can leverage the sparsity assumption from an estimation point of view, see e.g. \cite{Larsson07,DonohoJohnstone94,NearOracleBPDN,DantzigCandes}.\footnote{ALEX: I am not really sure what this sentence means. I think maybe it could be removed.}

In this work, we consider the estimation of an $S$-sparse, finite-dimensional vector $\xtrue \!\in\! \mathbb{R}^{N}\!$.
By ``$S$-sparse'' we mean 
%% in our context 
that the vector $\xtrue$ has 
at most
%% only 
$S$ nonzero entries, which is denoted by ${\| \xtrue \|}_{0} \triangleq | \supp( \xtrue ) | \le S$, 
where $\supp( \xtrue )$ denotes the set of indices of the nonzero entries of $\xtrue$.     %%% AJ
The ``sparsity'' $S$ is assumed to be known, and typically $S \!\ll\! N$. However, the positions of the nonzero entries (i.e., $\supp( \xtrue )$) as well as the values of the nonzero entries are unknown. %%% AJ
%% and the values of the non-zeros can be arbitrary.
%can be represented by much fewer than $N$ columns of a given fixed dictionary $\mathbb{C}^{N \times L}$ (where $L \geq N$), i.e., we have the following \emph{signal model} for the vector $\mathbf{x}$ which we would like to estimate:
%\begin{equation}
%\label{equ_sparse_signal_model}
%\mathbf{x} = \mathbf{D}  \mathbf{\alpha} \quad \quad \| \mathbf{\alpha} \|_{0}=S \ll N
%\end{equation}
%where $\| \mathbf{x} \|_{0}$ denotes the number of nonzeros of the vector $\mathbf{x}$.
%% In this paper, we 
%% want to 
We investigate how much we can gain in estimation accuracy by knowing \emph{a priori} that the vector $\xtrue$ is $S$-sparse. We will use the frequentist setting \cite{LC} 
of
%% for 
estimation theory, i.e., we will model $\xtrue$ as unknown but
deterministic. This is in contrast to 
%% the 
Bayesian estimation theory, where one treats $\xtrue$ as a random vector whose probability density function (pdf) 
or certain moments thereof are assumed to be known. In the Bayesian setting, the sparsity can be
%% is 
modeled by using a pdf that favors sparse 
%% parameter 
vectors, see e.g.\ \cite{Larsson07,Tipping01sparsebayesian,FBMP}.

A fundamental concept in the frequentist setting is that of unbiasedness \cite{LC,scharf91,kay}. An unbiased estimator is one whose expectation always equals the true underlying vector $\x$.
%% parameter. 
The restriction to unbiased estimators is important as
it excludes trivial and practically useless estimators,
%% . Such an approach 
and it allows us to study the difficulty of the estimation problem using established techniques such as the \CR bound (CRB) \cite{cramer45,scharf91,kay}. Another justification of unbiasedness
is that for typical estimation problems, when the variance of the noise is low, it is necessary for an estimator to be unbiased in order to achieve a small mean-squared estimation error (MSE) \cite{LC}. 

%%% AJ
These reasons notwithstanding, there is no guarantee that unbiased estimators are necessarily optimal. In fact, in many settings, including the scenario described in this paper, there exist biased estimators which are strictly better than any unbiased technique in terms of MSE \cite{1956_Stein,BlindMinimaxZvika,RethinkingBiasedEldar}. Nevertheless, for simplicity and because of the reasons stated above, we focus on bounds for unbiased estimation in this work. As we will see, bounds on unbiased techniques give some indication of the general difficulty of the setting, and as such some of our conclusions will be shown empirically to characterize biased techniques as well.

%The estimation of the sparse parameter vector $\mathbf{x}$ will be based on a linear observation model:
%\begin{equation}
%\label{equ_lin_observation}
%\mathbf{y} = \mathbf{A} \mathbf{x} + \mathbf{n}
%\end{equation}
%with the fixed known matrix $\mathbf{A} \in \mathbb{C}^{M \times N}$ and the noise vector $\mathbf{n} \sim \mathcal{N}(\mathbf{0}, \sigma^{2} \mathbf{I}_{NÊ\times N})$ which models additive white Gaussian noise (AWGN).
%The estimation problem given by the composition of the signal model \eqref{equ_sparse_signal_model} and the linear observation model \eqref{equ_lin_observation}:
%\begin{equation}
%\label{equ_slm}
%\mathbf{y} = \mathbf{A} \mathbf{x} + \mathbf{n}  \quad \quad  \mathbf{x} \in \{ \mathbf{D} \mathbf{\alpha} \, | \,  \| \mathbf{\alpha} \|_{0} = S \}
%\end{equation}
%is known as the ``sparse linear model''.
%For a channel estimation application we have the following interpretation of \eqref{equ_slm}: the vector $\mathbf{y}$ represents the observed received signal, the vector $\mathbf{x}$ consists of the unknown channel taps and
%the matrix $\mathbf{A}$ represents the known training signal.

%\paragraph{Contribution and Related Work}

Our main contribution is a characterization of the optimal performance of unbiased estimators $\hat{\mathbf{x}}(\mathbf{y})$ that are based on 
\vspace{-2mm}
observing
\be
\label{equ_slm}
\mathbf{y} \ist=\ist \mathbf{A} \xtrue + \mathbf{n}
\ee
where $\mathbf{A} \!\in\! \mathbb{R}^{M \times N}$ ($M \!\geq\! N$) is a known 
%% orthonormal matrix, 
matrix with orthonormal columns, 
i.e., $\mathbf{A}^{T} \rmv\mathbf{A}= \mathbf{I}_N$, and $\mathbf{n} \sim \mathcal{N} (\mathbf{0},\sigma^{2} \mathbf{I}_M)$ denotes zero-mean
%% additive 
white Gaussian noise with known variance $\sigma^{2}$
(here, $\mathbf{I}_N$ denotes the identity matrix of size $N \times N$). Note that without loss of generality we can then assume that $\mathbf{A} =Ê\mathbf{I}_N$ and $M=N$, 
i.e., $\mathbf{y} = \xtrue + \mathbf{n}$,
since premultiplication of the model \eqref{equ_slm} by $\mathbf{A}^{T}$ will reduce the estimation problem
to an equivalent problem $\mathbf{y}' = \mathbf{A}' \xtrue + \mathbf{n}'$ in which $\mathbf{A}' = \mathbf{A}^{T} \rmv \mathbf{A} = \mathbf{I}_N$ and the noise $\mathbf{n}' = \mathbf{A}^{T}\mathbf{n}$ is again zero-mean white Gaussian with variance $\sigma^{2}$.
Such a sparse signal model can be used, e.g., for channel estimation \cite{Carbonelli06} when %%% AJ
%% whenever 
%% it is assumed that 
the 
%% time-invariant 
channel consists only of 
%% very 
few significant taps and an orthogonal training signal is used \cite{OptimalTrainingDong}. Another application that fits our scope is image denoising using an orthonormal wavelet basis \cite{Donoho94idealspatial}. We note 
that parts of this work were previously presented in \cite{AlexZvikaICASSP}.

The estimation problem \eqref{equ_slm} with $\mathbf{A} = \mathbf{I}_N$ was studied by Donoho and Johnstone \cite{DonohoJohnstone94,Donoho92minimaxestimation}. %%% AJ
Their work was aimed at demonstrating asymptotic minimax optimality, i.e., they considered estimators having optimal worst-case behavior when the problem dimensions $N,S$ tend to infinity.
By contrast, we consider the finite-dimensional setting, and attempt to characterize the performance at each 
%% parameter 
value 
of $\x$, rather than analyzing worst-case behavior.
Such a ``pointwise'' approach was also advocated by the authors of \cite{ZvikaCRB,ZvikaSSP}, who studied the CRB for the sparse linear model \eqref{equ_slm} with arbitrary $\mathbf{A}$.
However, the CRB is a local bound, 
in the sense that
%% i.e., 
the performance characterization it provides
%% its performance characterization 
is 
only 
based on the statistical properties in the neighborhood of the 
%% parameter 
specific value of $\x$ being examined. In particular, the CRB 
for a given $\x$
is 
only 
based on a 
\emph{local} 
unbiasedness assumption, meaning that the estimator is only required to be unbiased at $\x$ and 
in
its infinitesimal neighborhood. Our goal  %%% AJ
in this paper
is to obtain performance bounds for the more restrictive case of globally unbiased estimators, i.e., estimators whose expectation equals the true $\x$
%%  parameter value 
for each $S$-sparse vector $\x$.
Since any globally unbiased estimator is also locally unbiased, our lower bounds will be  %%% AJ
tighter 
%% higher 
than those of \cite{ZvikaCRB,ZvikaSSP}.

Our contributions and the organization of this paper can be summarized as follows.
%% The outline of our contribution is as follows. 
In Section \ref{sec_est_theory_slm}, we show that whereas only one unbiased estimator exists for
%% in 
the ordinary (nonsparse) 
%% linear 
signal in noise model, there are infinitely many unbiased estimators for the sparse signal in noise model;
%%  We also briefly review
%%  %%  discuss 
%%  some
%%  established estimator designs which are not 
%%  %% directly 
%%  based on the unbiasedness rationale. 
on the other hand, 
none of them has uniformly minimum variance.
%% there exists no uniformly minimum variance unbiased estimator.
In Sections \ref{sec_lowerbound} and \ref{sec_upperbound}, 
%% we characterize the achievable performance 
%% of unbiased estimators in terms of the mean-squared error (MSE). 
%% in terms of mean-squared error (MSE) of unbiased estimators.
we characterize the performance of \emph{locally} minimum variance unbiased estimators by providing, respectively, lower and upper bounds on their mean-squared error (MSE). 
These bounds can equivalently be viewed as lower and upper bounds on the \emph{Barankin bound} \cite{GormanHero,barankin49}.
Finally, numerical studies exploring and extending our performance bounds and comparing them with established estimator designs are presented in Section \ref{sec_sim}.

\emph{Notation}: Throughout the paper, boldface lowercase letters (e.g., $\mathbf{x}$) denote column vectors while boldface uppercase letters (e.g., $\mathbf{M}$) denote matrices. We denote 
by $\mbox{tr}(\mathbf{M})$, $\mathbf{M}^{T}\!$, and $\mathbf{M}^{\dagger}$ the trace, transpose, and Moore-Penrose pseudoinverse of $\mathbf{M}$, respectively. The identity matrix of size $N\!\times\! N$ is denoted by $\mathbf{I}_N$. The notation $\mathbf{M} \succeq \mathbf{N}$ indicates that $\mathbf{M} \!-\! \mathbf{N}$ is a positive 
semidefinite matrix. %The number and the set of indices of the nonzero entries of a vector $\mathbf{x}$ are denoted by ${\| \mathbf{x} \|}_{0}$ and $\supp(\mathbf{x})$, respectively. 
The set of indices of the nonzero entries of a vector $\x$ is denoted by $\supp(\x)$, and ${\| \x \|}_0$ is defined as the size of this set.
The $k\ist$th entry of $\mathbf{x}$ is written $x_k$. We also use the signum function of a real number $y$, $\mbox{sgn}(y) \triangleq y/|y|$. The sets of nonnegative, nonpositive, and positive real numbers will be denoted by $\mathbb{R}_{+}$, $\mathbb{R}_{-}$, and $\mathbb{R}_{++}$, respectively.

%% \pagebreak %%%%%%%%%%%

%%%%%%%%%%%%%%%%%%%%%%%%%%%%%%%%%%%%%%%%%%%%%%%%%%%%%%%%%%
\section{The Sparse Signal in Noise Model}
\label{sec_est_theory_slm}
%%%%%%%%%%%%%%%%%%%%%%%%%%%%%%%%%%%%%%%%%%%%%%%%%%%%%%%%%%

\subsection{Problem Setting}
%%%%%%%%%%%%%%%%%%%%%%%%%%%%%%%%%%%%%%%%%%%%%%%%%%%%%%%%%%

Let $\xtrue \rmv\in\rmv \RR^N$ be an unknown deterministic 
%% parameter 
vector which is known to be $S$-sparse, i.e.,
\[
\xtrue \!\in\! \mathcal{X}_{S} \,, \qquad \text{with} \;\; \mathcal{X}_{S} \triangleq \{ \mathbf{x} \!\in\! \mathbb{R}^{N} : {\| \mathbf{x} \|}_{0} \rmv\leq\rmv S \} \,.
\]
The 
%% parameter 
vector $\xtrue$ is to be estimated based on the observation of a vector $\mathbf{y}$ which is 
the sum of $\xtrue$ and zero-mean
%% obtained from $\mathbf{x}$ with the addition of 
white Gaussian noise. Thus
\begin{equation}
\label{equ_ssnm}
\mathbf{y} \ist=\ist \xtrue + \mathbf{n} \,,  \quad\quad \text{with} \;\; \xtrue \!\in\! \mathcal{X}_{S} \, \mbox{,} \;\;
  \mathbf{n} \sim \mathcal{N} ( \mathbf{0} , \sigma^{2} \mathbf{I}_{N} ) 
\end{equation}
where the noise variance $\sigma^{2}$ is assumed to be nonzero and known. It follows that the pdf of $\y$, parameterized by the deterministic
%% fixed 
but unknown parameter $\xtrue \!\in\! \mathcal{X}_{S}$, is  %%% AJ
\begin{equation}
\label{equ_pdf}
f ( \mathbf{y};  \xtrue) \,=\, \normconstgauss \, \exp\rmv\rmv \bigg( \!\!-\frac{1}{2 \sigma^{2}} \| \mathbf{y} \!-\! \xtrue \|_2^{2} \bigg) .
\end{equation}
We refer to \eqref{equ_ssnm} as the \emph{sparse signal in noise model} (SSNM). As explained previously, settings of the form \eqref{equ_slm} with an orthonormal matrix $\A$ can be converted to the SSNM \eqref{equ_ssnm}.
%% Note that the 
The case $S \!=\! N$ corresponds to the situation in which no sparsity assumption is made. As we will see, this case
%% the case $S \!=\! N$ 
is fundamentally different from the sparse setting $S \!<\! N$, which is our %% main 
focus in this paper. 

%For our SSNM we have that
%\be
%\label{equ_par_pdf_ssnm}
%f(\mathbf{y}; \mathbf{x}) = \normconstgauss e^{- \frac{1}{2 \sigma^{2}} \| \mathbf{y} - \mathbf{x} \|^{2}_{2}}.
%\ee
An \emph{estimator} $\hat{\mathbf{x}}(\mathbf{y})$ of the parameter $\mathbf{x}$ is a function that maps (a realization of) the observation $\mathbf{y}$ to (a realization of) the estimated 
vector
%% signal 
$\hat{\mathbf{x}}$, i.e.,
\[
\hat{\mathbf{x}}(\cdot) : \mathbb{R}^{N} \!\rmv\rightarrow \mathbb{R}^{N} \!:  \mathbf{y} \mapsto \hat{\mathbf{x}}.
\]
With an abuse of notation, we will use the symbol $\hat{\mathbf{x}}$ for both the estimator 
%% itself 
(which is a function) and 
the
%% its 
estimate (a specific 
%% realization of the 
function value). The meaning should be clear from the context. 
%% From the context the exact meaning should be clear.
The question
now
is how we can exploit the information that $\xtrue$ is $S$-sparse in order to construct ``good'' estimators. 
%% Evidently, we first have to specify what we mean by
%% %% To answer that question we must quantify the meaning of 
%% a ``good'' estimator. We will measure the quality 
Our measure of the quality of an estimator 
$\hat{\mathbf{x}}(\cdot)$ 
for
%% at 
a given parameter value
%% vector 
$\xtrue \!\in\! \mathcal{X}_{S}$ will be
%% by 
the estimator's MSE,
%% mean-squared error (MSE), 
which is defined 
\vspace{-2mm}
as
%$\varepsilon(\xtrue; \hat{\mathbf{x}}$):
\[
%% \label{eq_def_mse}
\varepsilon( \xtrue; \hat{\mathbf{x}} )Ê\,\triangleq\, \mathsf{E}_{\xtrue} \big\{ \| \hat{\mathbf{x}}( \mathbf{y}) - \xtrue \|^{2}_{2} \big\} \,.
%= \frac{1}{(2 \pi \sigma^{2})^{N/2}} \int_{\mathbf{y}} \| \hat{\mathbf{x}}(\mathbf{y})- \xtrue \|^{2}_{2} e^{-\frac{1}{2 \sigma^{2}} \| y - \xtrue \|^{2}_{2}}  d \mathbf{y}
\]
Here, the notation $\mathsf{E}_{\xtrue} \{ \cdot \}$ means that the expectation is taken with respect to the pdf $f(\mathbf{y};  \xtrue)$ of the observation $\mathbf{y}$ parameterized by $\x$.
Note that even though $\xtrue$ is known to be $S$-sparse, the estimates $\hat{\mathbf{x}}$ are not constrained to be $S$-sparse.
%%% BEGIN ALEX 15042010
%(It has been shown in \cite{} that such a constraint would result in poor estimator performance.)
%%% END ALEX 15042010

The MSE can be written as the sum of a bias  %%% AJ
term
and a variance term, i.e.,
\[
%% \label{equ_bias_var_decomp}
\varepsilon( \xtrue; \hat{\mathbf{x}} ) \,=\, \| \mathbf{b}( \xtrue; \hat{\mathbf{x}} )Ê\|^{2}_{2} \,+ V(\xtrue; \hat{\mathbf{x}} )Ê
\]
where the bias $\mathbf{b}(\xtrue; \hat{\mathbf{x}} )Ê \triangleq \mathsf{E}_{\xtrue} \{ \hat{\mathbf{x}}(\mathbf{y}) \} -  \xtrue$ accounts for systematic estimation errors
and the variance 
%% term 
$V(\xtrue; \hat{\mathbf{x}} )Ê \triangleq \mathsf{E}_{\xtrue} \{  \| \hat{\mathbf{x}}(\mathbf{y})  -  \mathsf{E}_{\xtrue} \{\hat{\mathbf{x}}(\mathbf{y}) \}  \|_2^{2} \}$ accounts for 
%% the 
errors due to 
%% the 
random fluctuations of the estimate.
%% estimator. 
Thus, for unbiased estimators ($\mathbf{b}(\xtrue; \hat{\mathbf{x}} ) \rmv=\rmv \mathbf{0}$ for all $\xtrue \!\in\! \mathcal{X}_{S}$), the MSE is equal to the variance, i.e., 
$\varepsilon( \xtrue; \hat{\mathbf{x}} )$\linebreak %%%%%%%%
$ = V(\xtrue; \hat{\mathbf{x}})$.
%In general requiring a small variance will imply a large bias and vice-versa. %This behaviour is called the bias-variance tradeoff.

We will also consider the mean power (second moment) of an estimator,
\be
\label{equ_mean_power}
P(\xtrue; \hat{\mathbf{x}}) \ist\,\triangleq\,\ist \mathsf{E}_{\xtrue} \big\{ \| \hat{\mathbf{x}}(\mathbf{y}) \|^{2}_{2} \big\} 
  \,=\, V(\xtrue; \hat{\mathbf{x}}) \ist+\ist \| \mathsf{E}_{\xtrue} \{ \hat{\mathbf{x}}(\mathbf{y}) \} \|^{2}_{2} \,.
\ee
For unbiased estimators, $\| \mathsf{E}_{\xtrue} \{ \hat{\mathbf{x}}(\mathbf{y}) \} \|^{2}_{2} = \| \xtrue \|^{2}_{2}$; thus, minimizing the variance $V(\xtrue; \hat{\mathbf{x}})$ at a fixed $\xtrue \!\in\! \mathcal{X}_{S}$ among all unbiased estimators is equivalent to minimizing $P(\xtrue; \hat{\mathbf{x}})$.

\subsection{Estimator Design}
%%%%%%%%%%%%%%%%%%%%%%%%%%%%%%%%%%%%%%%%%%%%%%%%%%%%%%%%%%

Two well-established 
%% rationales for 
estimator designs are the least squares (LS) estimator 
defined by
%% $\LS_est$, which is defined by
\begin{equation}
\label{equ_LS_est}
\LS_est (\mathbf{y}) \,\triangleq\, \argmin_{\mathbf{x}' \in \mathcal{X}_{S}} \| \mathbf{y} \!-\! \mathbf{x}' \|_2^2
\vspace{.7mm}
\end{equation}
and the maximum likelihood (ML) estimator
defined by 
%% $\ML_est$, which is given by
\be
\label{equ_ML_est}
\ML_est(\mathbf{y}) \,\triangleq\, \argmax_{\mathbf{x}' \in \mathcal{X}_{S}} f ( \mathbf{y};  \mathbf{x}' ).
\vspace{.7mm}
\ee
For the SSNM, due to \eqref{equ_pdf}, the LS and ML estimators coincide; they are easily seen to be given by
\be
\label{equ_LS_ML_coincide}
\LS_est (\mathbf{y}) \ist\ist=\ist\ist \ML_est (\mathbf{y}) \ist\ist=\ist\ist {\mathsf P}_{\!S} ( \mathbf{y}  )
\ee
where $\mathsf{P}_{\! S}$ is an operator that retains the $S$ largest (in magnitude) components and zeros %%% AJ
%% zeroes 
out
all others. The LS/ML estimator is biased unless $S \!=\! N$. 
%More explicitly the vector $\hat{\mathbf{x}} = \mathbf{P}_{S} ( \mathbf{y} ) $ obtained by applying the operator $\mathbf{P}_{S}$ on the observation vector $\mathbf{y}$ is given as
%\be
%\left( \hat{\mathbf{x}} \right)_{k} = \begin{cases}
% \left( \mathbf{y} \right)_k ,  & \text{if } k \in \mathcal{L}_{\mathbf{y}} \\
% 0, & \text{else}
%\end{cases}
%\ee
%where $\mathcal{L}_{\mathbf{y}}$ denotes the indices of the $S$ largest (in magnitude) components of $\mathbf{y}$.
Note that this estimator is
%% the LS and ML estimators are 
not based
%% neither the LS nor the ML techniques are based 
on a direct minimization of the MSE. Indeed, if the sparsity constraint is
%% s are 
removed ($S \!=\! N$) and $N \!\ge\rmv 3$, 
%% then 
it has been shown \cite{1956_Stein,BlindMinimaxZvika,RethinkingBiasedEldar} that there exist estimators which yield a better MSE performance than that of the LS/ML estimator. %%% AJ
%% LS and ML estimators.
%% approaches.

The MSE $\varepsilon( \xtrue; \hat{\mathbf{x}} )$ of a specific estimator $\hat{\mathbf{x}}(\cdot)$ depends on the value of the parameter $\xtrue$. 
This makes it difficult to define optimality in terms of minimum MSE. For example, if an estimator $\hat{\mathbf{x}}(\cdot)$ performs well (i.e., has a small MSE) for a specific parameter value $\mathbf{x}_{1}$, 
%% then 
it may 
still
%% well 
%% be that it has 
exhibit poor performance (i.e., a large MSE) for a different %%% AJ
%% another 
parameter value $\mathbf{x}_{2}$.
Ideally, an optimal estimator should have minimum MSE for all parameter values \emph{simultaneously}. However, such an optimality criterion is 
%% in general 
unobtainable since %%% AJ
the minimum MSE achievable at a specific parameter value $\mathbf{x}_{1}$ is zero; it
%% . This minimum MSE 
is achieved by the trivial estimator $\hat{\mathbf{x}}(\mathbf{y}) \equiv \mathbf{x}_{1}$ which %%% AJ
is constant and 
%% thus 
completely ignores the observation $\mathbf{y}$. Therefore, if there were a \emph{uniformly minimum MSE} estimator, it would have to achieve zero MSE for all
parameter values, which is obviously impossible. Thus, requiring the estimator to minimize the MSE at all parameter values simultaneously makes no sense.  %%% AJ

One useful optimality criterion is the minimax approach, which considers the worst-case MSE %%% AJ
\[
\sup\limits_{\xtrue \in \mathcal{X}_{S}}  \MSE
\vspace{1mm}
\]
of an estimator $\hat{\mathbf{x}}(\cdot)$.
%\footnote{We use the convention that the supremum of an unbounded set is equal to $\infty$}
An optimal estimator in the minimax sense 
%% is any estimator which 
minimizes the worst-case MSE, i.e., 
%% which
is a solution of the optimization 
\vspace{-1mm}
problem
\[
%% \label{equ_min_max_opt}
\inf_{\hat{\mathbf{x}}(\cdot)} \sup_{\xtrue \in \mathcal{X}_{S}} \MSE \,.
\vspace{1mm}
\]
Considerable effort has been spent to identify minimax estimators for sparse models such as the SSNM 
%% given 
in \eqref{equ_ssnm}; see, e.g., \cite{DonohoJohnstone94,Donoho92minimaxestimation,Johnstone04needlesand}. However, these results only apply in the asymptotic regime, i.e., when $N,S \rightarrow \infty$. By contrast, our goal is to analyze estimator performance for 
%% fixed and 
finite problem dimensions. There are no known closed-form expressions of the minimax risk or of minimax-optimal estimators for the SSNM in this case.

%Note that for the SSNM \eqref{equ_ssnm} with $S=N$, i.e., without any sparsity constraint, one specific (there are others) optimum estimator solving \eqref{equ_min_max_opt} is the LS estimator $\LS_est$ and the ML estimator $\ML_est$.

In this work, rather than pursuing
%% examining 
the minimax criterion, we consider \emph{unbiased} estimators $\hat{\mathbf{x}}(\cdot)$ for the SSNM. An
unbiased estimator is one for which the bias $\mathbf{b}(\xtrue; \hat{\mathbf{x}})$ is zero for all $S$-sparse parameter vectors %%% AJ
%% values $\xtrue \in \mathcal{X}_{S}$, 
\vspace{-2mm}
i.e.,
\begin{equation}
\label{equ_unbiased_cond1}
\mathsf{E}_{\xtrue} \{ \hat{\mathbf{x}}(\mathbf{y}) \} =  \xtrue  \qquad \text{for all }\, \xtrue \!\in\! \mathcal{X}_{S} \,.
\end{equation}
%% To constrain 
Let $\mathcal{U}$ denote the set
%% class 
of all unbiased estimators $\hat{\mathbf{x}}(\cdot)$ for the SSNM.
%%  (as defined in \eqref{equ_unbiased_cond1}).  %%% AJ
%\[
%\mathcal{U} \,\triangleq\, \big\{ \hat{\mathbf{x}}(\cdot) \,\big|\, \mathsf{E}_{\xtrue} \{ \hat{\mathbf{x}}(\mathbf{y}) \} = \xtrue  \;\, \text{for all}\; \xtrue \rmv\in\rmv \mathcal{X}_{S} \big\} \,.
%\]
Constraining an estimator to be unbiased excludes such trivial estimators as $\hat{\mathbf{x}}(\mathbf{y}) \rmv\equiv\rmv \mathbf{x}_{1}$ where $\mathbf{x}_{1} \!\rmv\in\! \mathcal{X}_{S}$ is 
%% an arbitrary but 
some fixed $S$-sparse parameter vector. 
%Another justification for the restriction to unbiased estimators  %%% AJ
%%% for our SSNM 
%is that in order to achieve small $\MSE$ 
%%% for high SNR
%in the high-SNR regime, 
%i.e., for $S$-sparse parameter vectors
%%% values 
%$\xtrue \rmv\in\rmv \mathcal{X}_{S}$ 
%%% that additionally obey 
%for which 
%%% $\mbox{SNR} \triangleq 
%%% \frac{\| \xtrue \|^{2}_{2}}{N \sigma^{2}} 
%$\| \xtrue \|^{2}_{2} \gg\rmv N \sigma^{2}$, the estimator must be 
%%% it is necessary to be 
%unbiased \cite{ZvikaCRB,ZvikaSSP}.

\subsection{Unbiased Estimation for
%% in 
the SSNM}
%%%%%%%%%%%%%%%%%%%%%%%%%%%%%%%%%%%%%%%%%%%%%%%%%%%%%%%%%%

We now study the set
%% class 
$\mathcal{U}$ of unbiased estimators for the SSNM in more detail. In particular, we will show that with the exception of the case $S \!=\! N$, 
this set
%% the class 
is uncountably large, i.e., there are infinitely many unbiased estimators. 
%However, 
We will    %%% AJ 
also
%% then 
show that there exists no uniformly minimum variance unbiased estimator
%% UMVU 
%% for the SSNM 
unless $S \!=\! N$. In what follows, we will say that an estimator $\hat{\mathbf{x}}$ has a bounded MSE if $\MSE \leq C$ for all $\xtrue \!\in\! \mathbb{R}^{N}\rmv$, 
where $C$ is a constant which may depend on $N$, $S$, and 
\vspace{1mm}
$\sigma^{2}$.

\begin{theorem}
\label{thm_exit_uniqu_unbiasedness_S_N}
Consider the SSNM \eqref{equ_ssnm} with $S \!=\! N$, i.e., without a sparsity constraint, in which case $\mathcal{X}_{S} \rmv=\rmv \mathbb{R}^{N}\rmv$. Then, there exists exactly one unbiased estimator having bounded MSE (up to deviations having zero measure). This estimator 
is given by  $\hat{\mathbf{x}} (\mathbf{y}) \rmv=\rmv \mathbf{y}$, which equals
%% coincides with 
the LS/ML estimator
%% technique 
\vspace{1mm}
in \eqref{equ_LS_est}--\eqref{equ_LS_ML_coincide}.
\end{theorem} 

The proof of this result can be found in Appendix \ref{app_proof_ thm_exit_uniqu_unbiasedness_S_N}. 
By contrast with Theorem \ref{thm_exit_uniqu_unbiasedness_S_N}, when sparsity constraints are imposed there exists a large family of unbiased estimators, as we now 
\vspace{1mm}
show. 

\begin{theorem}
\label{thm_exit_uniqu_unbiasedness_S_lt_N}
%By contrast,
%% Contrariwise, 
For 
%% any 
$1 \!\leq\! S \!<\! N$, there are uncountably infinitely many unbiased estimators for the SSNM.
% Furthermore a necessary and sufficient condition for an estimator $\hat{\mathbf{x}}(\mathbf{y})$ to
%be unbiased is that its Fourier transform $\bar{\mathbf{x}}(\bar{\mathbf{y}})$ satisfies the following set of equations:
%\be
%\int_{\bar{\mathbf{y}}_{\mathcal{A}}} \bar{\mathbf{x}} (\bar{\mathbf{y}}) \cdot e^{- \frac{2 \sigma^{2}}{\pi}  \| \mathbf{\bar{\mathbf{y}}} \|^{2}_{2}} d \bar{\mathbf{y}}_{\mathcal{A}}= \int_{\bar{\mathbf{y}}_{\mathcal{A}}} \delta(\bar{\mathbf{y}}) d \bar{\mathbf{y}}_{\mathcal{A}} \quad \quad \forall %\mathcal{A} \subset \{1,...,N\} \quad |\mathcal{A}| = N-S
%\ee
%where $\bar{\mathbf{y}}_{\mathcal{A}}$ is any length $N$ vector whose non-zeros are only the elements with indices in the index set $\mathcal{A}$.
\end{theorem}

\begin{bproof}%we consider the case $S<N$ and prove the second part of Theorem \ref{thm_exit_uniqu_unbiasedness}.
Consider the class of estimators defined by 
\begin{equation}
\label{equ_example_class_unbiased_est_123} 
\hat{x}(\mathbf{y}) \eq \mathbf{y} +Ê a \ist y_1 \Bigg[ \prod_{k=2}^{S+1} h^{(c,d)}( y_k) \Bigg] \big( 1 \;\ist 0 \,\cdots\, 0 \big)^{T} , \quad\; a \rmv\in\rmv \mathbb{R}\mbox{,} \;\, c,d \rmv\in\rmv \mathbb{R}_{++}
\vspace{-1mm}
\end{equation}
where 
\vspace{.7mm}
\begin{equation}
h^{(c,d)}(y) \,\triangleq\ist \begin{cases} \mbox{sgn}(y) \ist, & |y| \in [c,c+d] \\[-1mm] 
%% [\frac{4}{10},\frac{6}{10}]
0 \ist, & \mbox{else}. \end{cases}
\end{equation}
A straightforward calculation shows that each estimator of this uncountably infinite class is an unbiased estimator for the SSNM. \hfill $\Box$
\vspace{3mm}
\end{bproof}

This (constructive) proof 
%% of Theorem \ref{thm_exit_uniqu_unbiasedness_S_lt_N} 
points at a noteworthy fact. Consider a particular parameter value $\mathbf{x}$.
By an appropriate choice of the parameters $a,c,d$ in \eqref{equ_example_class_unbiased_est_123},
one can reduce the magnitude of the estimate $\hat{\mathbf{x}}(\mathbf{y})$ for sets of realizations $\mathbf{y}$ with 
high probability, i.e., for which $f(\mathbf{y}; \mathbf{x})$ is large. This results in a reduced mean power and (since the estimator is unbiased) in a reduced variance and MSE at the specific parameter value $\mathbf{x}$. One can thus construct an unbiased estimator that performs better than the (biased) LS/ML estimator at the given $\mathbf{x}$.  

%% In the following, we will only consider the case $S<N$.
%% consider only the SSNM for $S<N$.
In view of Theorems \ref{thm_exit_uniqu_unbiasedness_S_N} and \ref{thm_exit_uniqu_unbiasedness_S_lt_N}, we will only consider the case $S \!<\! N$ in the following.  %%% AJ
Since in this case there are infinitely many unbiased estimators, we would like to find an unbiased estimator having minimum variance
(and, thus, minimum MSE) among all unbiased estimators. If there exists an unbiased estimator
%% method 
$\hat{\mathbf{x}}(\cdot) \!\in\rmv \mathcal{U}$ which
minimizes the variance \emph{simultaneously} for all %%% AJ 
$S$-sparse 
parameter vectors
%% values 
$\xtrue \!\in\! \mathcal{X}_{S}$, then this estimator
%% technique 
is called 
a
%% the 
\emph{uniformly minimum variance unbiased} (UMVU) estimator \cite{LC}.
In other words, a UMVU estimator for the SSNM 
%% is any estimator which 
solves the optimization problem
\be
\label{equ_opt_var_proof}
\argmin_{\hat{\mathbf{x}}(\cdot) \ist\in\, \mathcal{U}} \, V(\xtrue; \hat{\mathbf{x}}) 
%% \quad \quad \mbox{s.t.} \quad \quad \hat{\mathbf{x}} \in \mathcal{U}
\vspace{1mm}
\ee
simultaneously for all $\xtrue \!\in\! \mathcal{X}_{S}$.
In the nonsparse case $S \!=\! N$, it is well known that the LS estimator is the UMVU estimator \cite{scharf91}; however, in light of Theorem \ref{thm_exit_uniqu_unbiasedness_S_N}, this is not a very strong result, since $\LS_est$ is the \emph{only} unbiased estimator in that
%% this 
case. On the other hand, for the sparse case $S \!<\! N$, the following negative result is shown in Appendix \ref{app_proof_thm_no_umvu}. %%% AJ
%% we will now show that in the sparse case $S < N$, there exists no UMVU estimator.

\begin{theorem}
\label{thm_no_umvu}
For the SSNM with $S \!<\! N$, there exists no UMVU estimator, i.e., there is
%% exists 
no unbiased estimator $\hat{\mathbf{x}} \!\in\rmv \mathcal{U}$ that minimizes $V(\mathbf{x}; \hat{\mathbf{x}})$ simultaneously for all parameter vectors $\mathbf{x} \!\in\! \mathcal{X}_{S}$.
\end{theorem}

%% The proof of this theorem is provided in Appendix \ref{app_proof_thm_no_umvu}.
 
Despite the fact that a UMVU estimator does not exist for
%% in 
the SSNM, one can still attempt to solve the optimization problem \eqref{equ_opt_var_proof} separately for each value of $\x \!\in\! \mathcal{X}_{S}$. An unbiased estimator which solves \eqref{equ_opt_var_proof} for a specific value of $\x$ is said to be \emph{locally minimum variance unbiased} (LMVU) \cite{LC}. The MSE of this estimator at $\x$ cannot be improved upon by any unbiased estimator. When viewed as a function of $\x$, this minimum
%% optimum 
MSE is known as the \emph{Barankin bound} (BB) \cite{GormanHero,barankin49}.
Thus, the BB characterizes the minimum MSE achievable by any unbiased estimator for each value of $\x \!\in\! \mathcal{X}_{S}$;
%% . Indeed, the BB 
it is the highest and tightest lower bound on the MSE of unbiased estimators. As such, the BB
%% this bound also 
serves as a measure of the difficulty of estimating $\xtrue$.
%%  in the SSNM\@.

Computing the BB is equivalent to calculating $\min_{\hat{\mathbf{x}}(\cdot) \ist\in\, \mathcal{U}} V(\xtrue; \hat{\mathbf{x}})$
%% \eqref{equ_opt_var_proof} 
for each parameter vector $\xtrue \!\in\! \mathcal{X}_{S}$ separately. Unfortunately, there does not appear to be a simple closed-form 
expression
%% solution of \eqref{equ_opt_var_proof}.
%% Instead, in 
of the BB, and the numerical computation of the BB seems to be difficult as well. 
Therefore, in the remainder of this paper,
%% the next section,
%% rather than attempting to 
%% solve \eqref{equ_opt_var_proof}, 
%% calculating the BB, 
we will provide lower and upper bounds on the BB\@.
When these %lower and upper  %%% AJ 
bounds are close to one another, they provide  %%% AJ 
%% we have 
an accurate characterization of the BB.

\vspace{1mm}

%%%%%%%%%%%%%%%%%%%%%%%%%%%%%%%%%%%%%%%%%%%%%%%%%%%%%%%%%%
\section{Lower Bounds on the Minimum MSE}  %%% AJ 
\label{sec_lowerbound}
%%%%%%%%%%%%%%%%%%%%%%%%%%%%%%%%%%%%%%%%%%%%%%%%%%%%%%%%%%

\vspace{1mm}

In this section, we will develop a lower bound on the BB 
%% $\varepsilon( \xtrue; \hat{\mathbf{x}}^{(\xtrue)})$ 
(which is thus a lower bound on the MSE of any unbiased estimator) by calculating a limiting case of the Hammersley--Chapman--Robbins  %%% AJ 
bound \cite{GormanHero} for the SSNM.

%% \pagebreak %%%%%%%%%

\subsection{Review of the CRB}
\label{sec_crb}
%%%%%%%%%%%%%%%%%%%%%%%%%%%%%%%%%%%%%%%%%%%%%%%%%%%%%%%%%%

A variety of techniques exist for developing lower bounds on the MSE of unbiased estimators. The simplest 
%% of these 
is the CRB \cite{cramer45,rao45,kay}, which was previously derived for a more general sparse estimation setting in \cite{ZvikaCRB,ZvikaSSP}. In the current setting, i.e., for the SSNM \eqref{equ_slm}, the CRB
%% this bound 
is given by
\begin{equation}
\label{equ crb}
\varepsilon(\xtrue;\hat{\mathbf{x}}) \geq
\begin{cases}
S \sigma^2, & {\|\xtrue\|}_0 \rmv=\rmv S \\[-1.5mm]
N \sigma^2, & {\|\xtrue\|}_0 \rmv<\rmv S 
\end{cases}
\end{equation}
where $\hat{\mathbf{x}} \rmv\rmv\in\rmv \mathcal{U}$, i.e., $\hat{\mathbf{x}}(\cdot)$ is any unbiased estimator for the SSNM.

In the case of parameter values $\xtrue \!\in\! \mathcal{X}_{S}$ with non-maximal support, i.e., ${\|\xtrue\|}_0 \rmv<\rmv S$, the CRB is $N \sigma^2$. This is the MSE of the trivial unbiased estimator $\hat{\mathbf{x}} (\mathbf{y}) \rmv=\rmv \mathbf{y}$. Since the CRB is thus achieved by an unbiased estimator, we conclude that the CRB is a \emph{maximally tight} lower bound for ${\|\xtrue\|}_0 \rmv<\rmv S$; no other lower bound can be tighter (higher). We 
%% can 
also conclude that for ${\|\xtrue\|}_0 \rmv<\rmv S$, the trivial estimator $\hat{\mathbf{x}} (\mathbf{y}) \rmv=\rmv \mathbf{y}$ is the LMVU estimator; no other unbiased estimator can have a smaller MSE.

For parameter values $\xtrue \!\in\! \mathcal{X}_{S}$ with maximal support, i.e., ${\|\xtrue\|}_0 \rmv=\rmv S$, we will see that the CRB is not maximally tight, and the trivial estimator $\hat{\mathbf{x}} (\mathbf{y}) \rmv=\rmv \mathbf{y}$ is not the  %%% AJ 
LMVU estimator. Indeed, one problem with the CRB in \eqref{equ crb} is that it is discontinuous in the transition between ${\|\xtrue\|}_0 \rmv=\rmv S$ and ${\|\xtrue\|}_0 \rmv<\rmv S$. Since the MSE of any estimator is continuous \cite{LC}, 
%% this 
this discontinuity implies that the CRB is not the tightest lower bound obtainable for 
%% (uniformly) 
unbiased estimators.
In order to obtain tighter bounds for ${\|\xtrue\|}_0 \rmv=\rmv S$, it is important to realize that the CRB is a local bound, which assumes unbiasedness only in a neighborhood of $\xtrue$.
Since we are interested in estimators that are unbiased for all $\x \!\in\! \XS$, which is a more restrictive constraint than local unbiasedness, 
tighter (i.e., higher) lower bounds can be expected for unbiased estimators in the case ${\|\xtrue\|}_0 \rmv=\rmv S$.

\subsection{Hammersley--Chapman--Robbins Bound}
\label{sec_hcrb}
%%%%%%%%%%%%%%%%%%%%%%%%%%%%%%%%%%%%%%%%%%%%%%%%%%%%%%%%%%

An alternative lower bound for unbiased % (not only in a local sense)  %%% AJ 
estimators is the Hammersley--Chapman--Robbins bound (HCRB) \cite{hammersley50,ChapmanRobbins51,GormanHero}, which can be stated, in our context, as follows.

\begin{proposition} \label{pr:hcrb}
Given a parameter value $\xtrue \!\in\! \XS$, consider a set of $p$ ``test points'' $\{ \mathbf{v}_{i} \}_{i=1}^p$ such that $\xtrue+\mathbf{v}_i \in \XS$ for all $i = 1,\ldots,p$.
%% $1 \le i \le p$.
Then, the covariance of any unbiased estimator $\hat{\x}(\cdot)$, 
$C(\x;\hx) \triangleq \mathsf{E}_{\xtrue} \big\{ 
\big[ \hat{\mathbf{x}}( \mathbf{y}) - \mathsf{E}_{\xtrue} \{\hat{\mathbf{x}}(\mathbf{y}) \} \big] 
\big[ \hat{\mathbf{x}}( \mathbf{y}) - \mathsf{E}_{\xtrue} \{\hat{\mathbf{x}}(\mathbf{y}) \} \big]^T \big\}$, satisfies %%% AJ 
%% satisfies
\beq \label{equ hcr vector}
C(\x;\hx) \ist\ist\succeq\ist\ist \V \J^\pinv \V^T
%% \vspace{-3.5mm}
\pagebreak %%%%%%%%%%
\eeq
where
\vspace{-2mm}
\beq \label{equ def V}
\mathbf{V} \,\triangleq\, ( \mathbf{v}_{1} \cdots \mathbf{v}_{p} ) \in \mathbb{R}^{N \times p}
\eeq
and the $(i,j)$th entry of the matrix $\mathbf{J} \rmv\in \mathbb{R}^{p \times p}$ is given by
\beq \label{equ def J}
(\mathbf{J})_{i,j} \ist\triangleq\, \expp{\frac{\v_i^T \v_j}{\sigma^{2}}} - 1 \,.
\vspace{-1mm}
\eeq
%$\mathbf{J}^{\dagger}$ denotes the Moore-Penrose pseudoinverse of $\mathbf{J}$, and $\mathbf{A} \succeq \mathbf{B}$ means that
%the matrix $\mathbf{A} \!-\! \mathbf{B}$ is positive semidefinite.
In particular, the MSE of $\hat{\x}(\cdot)$ 
\vspace{-2mm}
satisfies %%% AJ
%% by
\begin{equation}
\label{equ hcr}
\varepsilon(\xtrue;\hat{\mathbf{x}}) \ist\ist\geq\ist\ist {\rm tr}\big(\V \J^\dagger \V^T\big) \,.
\vspace{1mm}
\end{equation}
\end{proposition}

%% Note that the lower bounds in \eqref{equ hcr vector} and \eqref{equ hcr} depend on 
%% $\xtrue$ due to the condition that $\xtrue+\mathbf{v}_i \in \XS$.
The proof of Proposition~\ref{pr:hcrb}, which can be found in Appendix~\ref{ap:pr:hcrb}, involves an application of the multivariate HCRB of Gorman and Hero \cite{GormanHero} to the SSNM
%% present 
setting. Note that both the number of test points $p$ and their values $\mathbf{v}_i$ are arbitrary and can depend on $\xtrue$. 
%(recall that the $\mathbf{v}_i$ have to satisfy the condition $\xtrue+\mathbf{v}_i \in \XS$).  %%% AJ 
In general, including additional test points $\mathbf{v}_i$ 
will result in a tighter HCRB \cite{GormanHero}. Our goal in this section is to choose test points $\mathbf{v}_i$ which result in a tight but analytically tractable bound.

Before attempting to derive a bound which is tighter than the CRB, we
%% let us 
first observe that the CRB itself can be obtained as the limit of a sequence of HCRBs with appropriately chosen test points. Indeed, consider the specific test points given
by\footnote{Note %%%%%%%%%%%
that, with a slight abuse of notation, the index $i$ of the test points is 
now
allowed to take on non-sequential values from 
%% the set $\supp(\xtrue)$ or from 
the set $\{ 1, \ldots, N \}$.} %%%%%%%%%%%%%
\begin{subequations}\label{equ test pts crb}
\begin{align}
{\{ t \ist\mathbf{e}_i \}}_{i \in \supp(\xtrue)} \,, \quad {\|\xtrue\|}_0 & \rmv=\rmv S
\label{equ test pts =s} \\[-1mm]
{\{ t \ist\mathbf{e}_i \}}_{i \in \{1,\ldots,N\}} \,, \quad {\|\xtrue\|}_0 &\rmv<\rmv S
\label{equ test pts <s}
\end{align}
\end{subequations}
where $t \rmv\rmv>\rmv\rmv 0$ is a constant and %%% AJ 
$\mathbf{e}_i$ represents the $i$th column of the $N \times N$ identity matrix. Note that $p \rmv\rmv=\rmv\rmv S$ in \eqref{equ test pts =s} %%% AJ
%% the first case 
and $p \rmv\rmv=\rmv\rmv N$ in \eqref{equ test pts <s}.
%% the second case.
Each value of $t$ yields a different set of test points and, via Proposition \ref{pr:hcrb}, a different lower bound on the MSE of unbiased estimators. We show in Appendix~\ref{ap:t->0} that the CRB in \eqref{equ crb} is the limit of a sequence of such bounds as $t \rmv\rmv\rightarrow\rmv\rmv 0$,  
%% Suppose now that a sequence of such bounds is taken, in which $t \rightarrow 0$. Then, as shown in 
%% Appendix~\ref{ap:t->0}, the limit of this sequence of bounds is precisely the CRB of \eqref{equ crb}. 
%% Furthermore, it is shown in Appendix~\ref{ap:t->0} that \eqref{equ crb} 
and that it is tighter than any bound 
that can be 
obtained via Proposition \ref{pr:hcrb} using the test points 
%% in
%% of the form 
\eqref{equ test pts crb} for a fixed $t \rmv\rmv>\rmv\rmv 0$.

%% \pagebreak %%%%%%%%%

Can a set of test points different from \eqref{equ test pts crb} yield a lower bound that is tighter %%% AJ 
%% tighter bound 
than the CRB\@? 
As discussed above, this is only possible for parameter values $\xtrue$ having maximal support, i.e., ${\|\xtrue\|}_0 \rmv=\rmv S$, because for ${\|\xtrue\|}_0 \rmv<\rmv S$ the CRB is already maximally tight.
Therefore, let us consider a parameter $\xtrue$ with ${\|\xtrue\|}_0 \rmv=\rmv S$. Suppose one of the entries within the support, $x_j$ for some $j \in \supp(\xtrue)$, has a small magnitude. Such a parameter $\xtrue$
just  %%% AJ  
barely qualifies as having maximal support, so it makes sense to adapt the optimal test points \eqref{equ test pts <s} from the non-maximal support case. However, including a test point $t \ist\e_i$ 
with $i \notin \supp(\xtrue)$
is not allowed,
%%  if $i \notin \supp(\xtrue)$, 
since in this case $\xtrue + t \ist\e_i$ is not in $\XS$. Instead, one could include the test point $\v_i = t \ist\e_i - x_j \ist \e_j$,
%% $\v_i = -(\xtrue)_j \e_j + t \ist\e_i$, 
which 
satisfies the requirement $\xtrue + \v_i \in \XS$ and is still close to $t\ist\e_i$ if $x_j$ is small.
%% is close to $t\ist\e_i$ if $(\xtrue)_j$ is small, and satisfies the requirement $\xtrue + \v_i \in \XS$. 
More 
generally,
%% precisely, given a 
for any 
maximal-support parameter $\xtrue$, we propose the set of $N$ test points given 
\vspace{-1mm}
by
\begin{equation}
\label{equ_test_points}
\mathbf{v}_i =
\begin{cases}
t \ist\e_i \,, & i \in \supp(\xtrue) \\[-1.5mm]
 t \ist\e_i - \xitrue \ist \e_{(S)} \,,
%% -\xitrue \ist \e_{(S)} + t \ist\e_i , 
& i \notin \supp(\xtrue)
\end{cases}
\vspace{1mm}
\end{equation}
for $i = 1, \ldots, N$. Here, $\xitrue$ denotes the 
%% $S$-largest
%% %% smallest 
%% (in magnitude) component of $\xtrue$---i.e., the 
smallest (in magnitude) of the $S$ nonzero components of $\xtrue$ and $\e_{(S)}$ denotes the corresponding 
unit vector.
%% direction. 
These test points $\v_i$ satisfy the condition $\xtrue + \v_i \in \XS$. 
Note that the test points 
in
\eqref{equ test pts =s}, which yield the CRB, are a subset of the test points in \eqref{equ_test_points}. It can be shown \cite{GormanHero} that this implies that the bound induced by \eqref{equ_test_points} 
%% resulting bound 
will always be at least as tight as that obtained from \eqref{equ test pts =s}.
%% \eqref{equ test pts crb}. 
It is important to note that
%% Significantly, 
\eqref{equ_test_points} uses $N$ test points for parameter values 
%% $\xtrue \in \mathcal{X}_{S}$ 
with
%% having 
maximal support, 
just as \eqref{equ test pts <s} does for parameter values with non-maximal support. 
In fact, 
there is 
a smooth transition 
%% is formed 
between the optimal test points \eqref{equ test pts <s} for non-maximal support and the proposed test points \eqref{equ_test_points} for maximal support.
%%  points.

While an expression of the HCRB can be obtained by simply plugging \eqref{equ_test_points} into \eqref{equ hcr},
%% Proposition~\ref{pr:hcrb}, 
the resulting bound is extremely cumbersome and not very insightful. Instead, in analogy to the derivation of the CRB above, one can obtain a simple result by taking the limit 
%% of \eqref{equ hcr} using test points according to \eqref{equ_test_points} 
for $t \rmv\rmv\to\rmv\rmv 0$. This leads to 
%% The result is 
the following theorem, which combines the cases of maximal support (\eqref{equ hcr} using \eqref{equ_test_points} for $t \rmv\rmv\to\rmv\rmv 0$) and non-maximal support (\eqref{equ hcr} using \eqref{equ test pts <s} for $t \rmv\rmv\to\rmv\rmv 0$), and
whose proof can be found in Appendix~\ref{ap:th:hcrb}.

\begin{theorem}
\label{th:hcrb}
The MSE of any unbiased estimator $\hat{\mathbf{x}} \in \mathcal{U}$ for the SSNM satisfies  %%% AJ 
%% below by
\vspace*{1mm}
\begin{equation}
\label{equ_HCR}
\varepsilon(\xtrue;\hat{\mathbf{x}}) \,\geq\, \emph{HCRB}(\xtrue) \ist \triangleq
\begin{cases}
S\sigma^2 + (N \!-\! S \!-\! 1) \ist\ist e^{-\xitrue^{2}/\sigma^{2}} \sigma^2 , & {\|\xtrue\|}_0 \rmv=\rmv S \\[-1.5mm]
N\sigma^2 , & {\|\xtrue\|}_0 \rmv<\rmv S \,,
\end{cases}
%% \vspace{-1mm}
\end{equation}
where, in the case ${\|\xtrue\|}_0 \rmv=\rmv S$, $\xitrue$ is the smallest (in magnitude) of the $S$ nonzero entries of $\xtrue$.
\end{theorem}

For simplicity, we will continue to 
refer to %%% AJ 
%% refer to \eqref{equ_HCR} as 
\eqref{equ_HCR} as %%% AJ 
an HCRB\@, even though it was obtained as a limit of HCRBs. 
%% Strictly speaking, the bound in Theorem~\ref{th:hcrb} is obtained as a limit of HCRBs. 
%% However, for simplicity we will continue to refer to \eqref{equ_HCR} as an HCRB\@.
%% As expected, 
Note that when ${\|\xtrue\|}_0 \!<\! S$, the HCRB 
in \eqref{equ_HCR}
%% of 
%% Theorem~\ref{th:hcrb} 
is identical to the CRB 
in \eqref{equ crb}, 
since in that case the CRB is maximally tight and cannot be improved.
%% at those points. 
The HCRB also approaches the CRB when 
${\|\xtrue\|}_0 \rmv=\rmv S$ and
all components of $\xtrue$ are much larger than $\sigma$: here
%% in this case 
$e^{-\xitrue^2/\sigma^2}$ is negligible and the respective bound in \eqref{equ_HCR} converges to $S\sigma^2$, which is equal to the CRB in \eqref{equ crb}. This is due to the fact that the CRB is 
achieved by the ML
%% maximum likelihood 
estimator 
asymptotically\footnote{This %%%%%%%%%%%
can be explained by the fact that according to \eqref{equ_LS_ML_coincide}, the ML estimator for the SSNM retains the $S$ largest components in $\mathbf{y}$ and zeros out all other components. 
For noise variances $\sigma^{2}$ that are extremely small compared to the nonzero entries, i.e., for $\xitrue^2/\sigma^2 \to \infty$, the probability that the ML estimator selects the true components becomes 
very close to one. Therefore, for high 
%% values of 
$\xitrue^2/\sigma^2$, the ML estimator behaves like an oracle estimator which knows the support of $\mathbf{x}$ and whose MSE is equal to $S \sigma^{2}$. 
} %%%%%%%%%%% 
as $\xitrue^2/\sigma^2 \to \infty$, and is therefore also maximally tight when $\xitrue \rmv\gg\rmv \sigma$.
Furthermore, if we define the ``worst-case component SNR'' (briefly denoted as SNR) as $\xitrue^2/\sigma^2$,
%% (note the difference from the global SNR defined previously as $\mbox{SNR} = \| \xtrue \|^{2}_{2}/(N \sigma^{2})$), 
then Theorem \ref{th:hcrb} hints that the convergence to the high-SNR limit is exponential in the SNR.
%However, while the CRB predicts an abrupt transition between an MSE of $N\sigma^2$ at sub-maximal support points and $S \sigma^2$ at high SNR (i.e. for $\xitrue/\sigma \gg 1$), the HCRB forms a smoother transition, hinting that the rate of convergence to the high-SNR limit is exponential in the SNR\@ (the SNR could be defined by $\xitrue^{2}/\sigma^{2}$).

One of the motivations for improving the CRB \eqref{equ crb} was that \eqref{equ crb} 
is discontinuous in the transition between ${\|\xtrue\|}_0 \!=\! S$ and ${\|\xtrue\|}_0 \!<\! S$. 
%% Unfortunately, 
While the HCRB \eqref{equ_HCR} 
%% bound of Theorem~\ref{th:hcrb} 
is still discontinuous 
in this transition,
%% at these points, 
%% but 
the discontinuity
%% jump 
%% in this transition (whose size is $\sigma^{2}$) 
is much smaller than that of the CRB\@.
%%  (whose size is $(N-S)\sigma^{2}$). 
Indeed, the transition from ${\|\xtrue\|}_0 \rmv=\rmv S$ to ${\|\xtrue\|}_0 \rmv<\rmv S$ corresponds to $\xitrue \rmv\rmv\to\rmv\rmv 0$,
in which case the first bound in \eqref{equ_HCR} tends to $(N \!-\! 1) \ist \sigma^2$, whereas the second bound, valid for ${\|\xtrue\|}_0 \!<\! S$, is $N \sigma^{2}$; thus, the difference between the two bounds in \eqref{equ_HCR} is $\sigma^{2}$.
By contrast, the difference between the two bounds in \eqref{equ crb} is $(N \!-\! S) \ist \sigma^{2}$, which is typically much larger.
Again,
%% Since the MSE of any estimator is continuous \cite{LC}, 
%% this 
the discontinuity of \eqref{equ_HCR} 
implies that \eqref{equ_HCR} is not the tightest lower bound obtainable for unbiased estimators. In Section \ref{sec_sim}, we will demonstrate experimentally
%% show 
that this discontinuity can be eliminated altogether by 
%% adding more 
using a much larger number of test points. % compared to \eqref{equ_test_points}. %%% AJ 
However, in that case
%%  but 
the resulting bound no longer has a simple closed-form expression and can only be evaluated numerically. %%% AJ 

%% \pagebreak %%%%%%%%%

\vspace{1mm}

%%%%%%%%%%%%%%%%%%%%%%%%%%%%%%%%%%%%%%%%%%%%%%%%%%%%%%%%%%
\section{Upper Bound on the Minimum MSE}
\label{sec_upperbound}
%%%%%%%%%%%%%%%%%%%%%%%%%%%%%%%%%%%%%%%%%%%%%%%%%%%%%%%%%%

\vspace{1mm}

%% \subsection{Upper Bound on the MSE}
%% \label{sec_upperbound}
%% %%%%%%%%%%%%%%%%%%%%%%%%%%%%%%%%%%%%%%%%%%%%%%%%%%%%%%%%%%

%The lower bound $\mbox{HCRB}(\xtrue)$ on the MSE $\MSE$ of an unbiased estimator $\hat{\mathbf{x}} \in \mathcal{U}$ derived in the previous section for the SSNM is tight for all parameter vector values $\xtrue \in \mathcal{X}_{S}$ that either have less than $S$ non-zeros, i.e., $\| \xtrue \|_{0} < S$ or which
%have exactly $S$ non-zeros and additionally the magnitudes of the non-zeros are all much larger than the noise variance $\sigma^{2}$. However, for all other parameter vectors $\xtrue \in \mathcal{X}_{S}$

As pointed out in the previous section, the lower bound 
%% given by 
$\mbox{HCRB}(\xtrue)$ on the BB 
%% in \eqref{equ_HCR} 
is not maximally tight since it is discontinuous in the transition between parameter vectors with
%% having 
${\|\xtrue\|}_0 \!=\! S$ and those with
%% having 
${\|\xtrue\|}_0 \rmv<\rmv S$. In other words, there is a gap between the HCRB and the BB\@. 
How large is this gap?
%% The question arises as to how large 
%% this gap is.
%% %% the gap actually is. 
We will address this issue by deriving an \emph{upper} bound on the BB.
%% , i.e., 
%% %% an upper bound 
%% on the minimum MSE/variance $\min_{\hat{\mathbf{x}}(\cdot) \in\, \mathcal{U}} V(\xtrue; \hat{\mathbf{x}})$ achieved by an LMVU estimator 
%% as defined by the optimization problem \eqref{equ_opt_var_proof}.
%%  at the parameter vector $\xtrue$.
This will be done by finding a
%% n approximate ???
constrained solution of \eqref{equ_opt_var_proof}.
If this upper bound is close to the lower bound $\mbox{HCRB}(\xtrue)$, we can conclude that both bounds are fairly tight and thus provide a fairly accurate
characterization of the BB. As before,
%% in the later part of Section \ref{sec_lowerbound}, 
we consider the nontrivial case ${\| \xtrue \|}_{0} \rmv=\rmv S$.

%% It will be convenient to use an equivalent form of 
%% %% the optimization problem 
%% \eqref{equ_opt_var_proof}. First, as 
We first note (cf.\ \eqref{equ_mean_power}) that \eqref{equ_opt_var_proof} is equivalent to the optimization problem
$\argmin_{\hat{\mathbf{x}}(\cdot) \ist\in\, \mathcal{U}} \mathsf{E}_{\xtrue} \big\{ \| \hat{\mathbf{x}}(\mathbf{y}) \|^{2}_{2} \big\}$\linebreak %%%%%%%%
$= \argmin_{\hat{\mathbf{x}}(\cdot) \ist\in\, \mathcal{U}}  \sum_{k=1}^N  \mathsf{E}_{\xtrue} \big\{ ( \hat{x}_{k} ( \mathbf{y} ) )^{2} \big\}$,
%%  which consists of $N$ independent scalar optimization problems. 
where $\hat{x}_{k}$ denotes the $k$th entry of $\hat{\mathbf{x}}$.
This, in turn, is equivalent to 
%% solving for $k=1...N$ 
the $N$ individual scalar optimization problems
\begin{equation}
\label{equ_scalar_opt}
\argmin_{\hat{x}_{k}(\cdot) \ist\in\, \mathcal{U}^k} \ist \mathsf{E}_{\xtrue} \big\{ ( \hat{x}_{k} ( \mathbf{y} ) )^{2} \big\} \,, \qquad k = 1,\ldots,N 
%% \quad \quad \mbox{s.t.} \quad \quad \hat{x}_{k} \in \mathcal{U}^{k}
\vspace{1mm}
\end{equation}
%% for $k = 1,\ldots,N$,
where $\mathcal{U}^{k}$ denotes the set of unbiased estimators of the $k$th entry of $\xtrue$, i.e., 
\[
\mathcal{U}^{k} \ist\triangleq\, \big\{ \hat{x}_k(\cdot) \,\big|\, \mathsf{E}_{\mathbf{x}} \{ \hat{x}_k(\mathbf{y}) \} = x_k  \;\, \text{for all}\; \mathbf{x} \!\in\! \mathcal{X}_{S} \big\} \,.
\]
%% the set of estimators $\hat{x}_{k}(\mathbf{y})$ satisfying $\mathsf{E}_{\xtrue} \{ \hat{x}_{k} ( \mathbf{y} ) \} = x_k$ for all $\xtrue \in \mathcal{X}_{S}$.
By combining the unbiased estimators $\hat{x}_{k}(\cdot)$ for $k = 1,\ldots,N$ into a vector, we obtain an unbiased estimator of the parameter $\xtrue$.
%In particular if $\hat{x}^{(\xtrue)}_{k}$ denotes the solution of
%\eqref{equ_scalar_opt} then the solution of \eqref{equ_opt_var_proof} is given by the vector valued estimator $\hat{\mathbf{x}}^{(\xtrue)} (\mathbf{y}) = \begin{pmatrix} \hat{x}^{(\xtrue)}_{1}(\mathbf{y}) &\hat{x}^{(\xtrue)}_{2}(\mathbf{y}) & ... & \hat{x}^{(\xtrue)}_{N}(\mathbf{y}) \end{pmatrix}^{T}$.

It will be convenient to 
%% We can always 
write the $k$th scalar estimator as
%% now decompose the optimization variable in \eqref{equ_scalar_opt} (which is a scalar estimator $\hat{x}_{k}(\mathbf{y})$) into
\begin{equation}
\label{equ_scalar_est_additive}
\hat{x}_{k}( \mathbf{y}) \ist\ist=\ist\ist y_{k} + \hat{x}'_{k}(\mathbf{y}) 
\end{equation}
with $\hat{x}'_{k}(\mathbf{y}) \triangleq \hat{x}_{k}( \mathbf{y}) - y_{k}$.
Since for any $\hat{x}_{k}(\cdot) \!\in\rmv \mathcal{U}^{k}$ we have $\mathsf{E}_{\xtrue} \{ \hat{x}_{k} (\mathbf{y}) \} = \mathsf{E}_{\xtrue} \{ y_{k} \} + \mathsf{E}_{\xtrue} \{ \hat{x}'_{k} (\mathbf{y}) \} = x_{k} + \mathsf{E}_{\xtrue} \{ \hat{x}'_{k} (\mathbf{y}) \}$, 
%% we have that 
the unbiasedness condition $\hat{x}_{k}(\cdot) \rmv\in \mathcal{U}^{k}$ is equivalent to 
%% requiring that the mean of $\hat{x}'_{k}$ is zero for all parameter vectors, i.e.,
\[
\mathsf{E}_{\xtrue} \{ \hat{x}'_{k} (\mathbf{y}) \} = 0 \quad \quad  \text{for all }\, \xtrue \!\in\! \mathcal{X}_{S} \,.
\]
For $k \in \supp(\xtrue)$,
%%  the optimization problem \eqref{equ_scalar_opt} is readily solved as shown in 
the solution of the optimization problem \eqref{equ_scalar_opt} is stated in the following lemma,
%% theorem, 
which is proved in Appendix \ref{app_proof_opt_est_k_in_supp}.
In what follows, it will be convenient to denote by $\hat{\mathbf{x}}^{(\xtrue)}(\mathbf{y})$ a solution of the optimization problem \eqref{equ_opt_var_proof} 
for a given parameter vector $\xtrue \!\in\! \mathcal{X}_{S}$.
We recall that the estimator $\hat{\mathbf{x}}^{(\xtrue)}(\mathbf{y})$ is an
%% the 
LMVU at the parameter value
%% vector 
$\xtrue$, and 
%%  for the SSNM. 
%% Therefore 
its MSE, $\varepsilon( \xtrue; \hat{\mathbf{x}}^{(\xtrue)}) = \min_{\hat{\mathbf{x}}(\cdot) \in\, \mathcal{U}} V(\xtrue; \hat{\mathbf{x}})$,
%% ---which 
equals the BB at $\x$.
%% ---is the tightest lower bound on the MSE/variance at $\xtrue$ of any unbiased estimator for the SSNM.

\begin{lemma}
\label{equ_opt_k_in_supp}
%% Given a fixed 
Consider a parameter vector $\xtrue \!\in\! \mathcal{X}_{S}$ with maximal support, i.e., ${\| \xtrue \|}_{0} \!=\! S$. Then, for any $k \in \supp(\xtrue)$, 
the solution of the optimization problem \eqref{equ_scalar_opt} is given by 
\[
\hat{x}^{(\xtrue)}_{k}(\mathbf{y}) \ist=\ist y_{k} \,, \qquad k \in \supp(\xtrue) \,.
\]
Moreover, this is the LMVU for $k \in \supp(\xtrue)$. The MSE of this estimator is $\sigma^2$.
\end{lemma}

%% \pagebreak %%%%%%%%%

Because 
Lemma
%% Theorem 
\ref{equ_opt_k_in_supp}
%% we have found 
describes the scalar LMVU estimators for all
%% the 
indices $k \in \supp(\xtrue)$, it remains to  %%% AJ 
%% in the following 
consider the scalar problem \eqref{equ_scalar_opt} for $k \notin \supp (\xtrue)$.
Since 
%% the MSE 
$\varepsilon( \xtrue; \hat{\mathbf{x}}^{(\xtrue)})$ is the minimum 
of $\varepsilon(\xtrue; \hat{\mathbf{x}})$ as defined by the 
%% value of the minimization 
optimization problem \eqref{equ_opt_var_proof}, we can obtain an upper bound on $\varepsilon( \xtrue; \hat{\mathbf{x}}^{(\xtrue)})$ by
placing further constraints on the estimator $\hat{\mathbf{x}}(\cdot)$ 
%% optimization variable, 
to be optimized.
%% , i.e., the estimator $\hat{\mathbf{x}}(\cdot)$ in \eqref{equ_opt_var_proof}. 
We will thus consider the modified optimization 
\vspace{-1mm}
problem
\begin{equation}
\label{equ_opt_approx_bb}
\argmin_{\hat{\mathbf{x}}(\cdot) \ist\in\, \mathcal{U} \ist\cap\ist \mathcal{A}_{\xtrue}} \! V(\xtrue; \hat{\mathbf{x}})  
%% \quad \quad \mbox{s.t.} \quad \quad \hat{\mathbf{x}} \in \mathcal{U} \cap \mathcal{A}_{\xtrue}
\vspace{.7mm}
\end{equation}
%where the additional constraint is that the estimator should be also a member of a given class $\mathcal{A}_{\xtrue}$ of estimators. The notation $\mathcal{A}_{\xtrue}$ should
%indicate that the class may depend on the fixed value $\xtrue$ of the parameter since the optimzation problem \eqref{equ_opt_var_proof} also depends on it.
where the set $\mathcal{A}_{\xtrue}$ is chosen 
%% below so as to obtain 
such that a simpler 
%% optimization 
problem 
is obtained.
%% 
%% Obviously, this is an ad-hoc method, and the gap between the solution of the modified optimization problem 
%% \eqref{equ_opt_approx_bb} and the original one \eqref{equ_opt_var_proof} can be large.
%% However, our goal of obtaining an upper bound on the minimum MSE of unbiased estimators can be achieved by this method.
%% %we are interested in an upper bound on the minimum value of \eqref{equ_opt_var_proof}, i.e. 
%% the MSE $\varepsilon( \xtrue; \hat{\mathbf{x}}^{(\xtrue)})$ on the LMVU $\hat{\mathbf{x}}^{(\xtrue)}$ for the
%% %parameter vector $\xtrue \in \mathcal{X}_{S}$ and the minimum value of \eqref{equ_opt_approx_bb} is always an upper bound.
%% If, moreover, this upper bound is close to the $\mbox{HCRB}(\xtrue)$ of Section \ref{sec_lowerbound}, 
%% we can conclude that both bounds are fairly tight.
We will define $\mathcal{A}_{\xtrue}$ in a componentwise fashion. More specifically, 
%% we define $\mathcal{A}_{\xtrue}$ according to properties required of each of the components of $\hx(\y)$. The 
the $k$th component $\hat{x}_k(\y)$ of $\hx(\y)$, where $k \notin \supp (\xtrue)$, is said to belong to the set $\mathcal{A}_{\xtrue}^{k}$ if
%% , when written in the form
%% \be
%% \label{equ_scalar_est_additive}
%% \hat{x}_{k} (\mathbf{y}) = y_k + \hat{x}'_{k} (\mathbf{y}),
%% \ee
the correction term $\hat{x}'_{k}(\mathbf{y}) = \hat{x}_{k}( \mathbf{y}) - y_{k}$ (see \eqref{equ_scalar_est_additive}) satisfies the following two properties.
%to a class $\mathcal{A}_{\xtrue}^{k}$ of scalar estimators. This class in turn consists of estimators that are
%composed of the $k$th component of the observation, i.e., $y_k$ and an additive term $\hat{x}'_{k}$ that belongs to the class $\mathcal{B}_{\xtrue}^{k}$, i.e.,
%\be
%\mathcal{A}_{\xtrue}^{k} = \left\{ \hat{x}_{k} | \hat{x}_{k} (\mathbf{y}) = y_k + \hat{x}'_{k} (\mathbf{y}) \quad \mbox{,}  \quad \hat{x}'_{k} (\mathbf{y})  \in \mathcal{B}_{\xtrue}^{k} \right\}
%\ee
%The elements $\hat{x}'_{k}$ of $\mathcal{B}_{\xtrue}^{k}$  are defined via requiring the following requirements:
\begin{itemize}
\item {\em Odd symmetry} with respect to $k$ and all indices in $\supp(\xtrue)$:
\be
\hat{x}_k^\prime (\ldots,- y_{l}, \ldots) \,=\, -\, \hat{x}_k^\prime ( \ldots, y_{l}, \ldots) \,, \qquad \text{for all }\, l \in \{k\} \cup \supp(\xtrue) \,.
\label{equ_constr_odd_symm}
\ee
\item {\em Independence}
%% The function $\hat{x}_k^\prime (\mathbf{y})$
%%  \!\in\! \mathcal{A}_k^\prime$
with respect to
%% from 
all other indices:
%% indices not in the set $\{k\} \cup \supp(\xtrue)$:
\be
\hat{x}_k^\prime (\ldots, y_{l}, \ldots) \,=\, \hat{x}_k^\prime ( \ldots, 0, \ldots) \,, \qquad \text{for all }\, l \notin \{k\} \cup \supp(\xtrue) \,.
\label{equ_constr_invariance}
\ee
\end{itemize}
We then define 
$\mathcal{A}_{\xtrue}$ 
%% is now defined 
as the set of estimators $\hx(\y)$ such that $\hat{x}_k(\y) \rmv\rmv \in\rmv\rmv \mathcal{A}_{\xtrue}^k$ for all $k \notin \supp (\xtrue)$.
%% $k=1,\ldots,N$.
Note that any function $\hat{x} (\mathbf{y}) \rmv\rmv \in\rmv\rmv \mathcal{A}^{k}_{\xtrue}$
is fully specified by its values for all arguments $\mathbf{y}$ such that $\supp(\mathbf{y}) = \{k\} \ist\cup\, \supp(\xtrue)$ and all entries of $\mathbf{y}$ are nonnegative. The values
of $\hat{x}(\mathbf{y})$ for all other 
%%  argument vectors 
$\mathbf{y}$ follow by the decomposition \eqref{equ_scalar_est_additive} and the properties \eqref{equ_constr_odd_symm} and \eqref{equ_constr_invariance}.

%% \pagebreak %%%%%%%%%

To solve the modified optimization problem \eqref{equ_opt_approx_bb}, we consider the equivalent scalar 
\vspace{-1mm}
form
\begin{equation}
\label{equ_scalar_opt_modified}
\argmin_{\hat{x}_{k}(\cdot) \ist\in\, \mathcal{U}^{k} \ist\cap\ist \mathcal{A}_{\xtrue}^{k}} \!\mathsf{E}_{\xtrue} \big\{ ( \hat{x}_{k} ( \mathbf{y} ) )^{2} \big\}  
  \,, \qquad k \notin \supp(\xtrue) \,.
%% \quad \quad \mbox{s.t.} \quad \quad \hat{x}_{k} \in \mathcal{U}^{k} \cap  \mathcal{A}_{\xtrue}^{k}
\vspace{.7mm}
\end{equation}
%% for $k \notin \supp(\xtrue)$.
The resulting minimum MSE is stated by the following lemma,
%% theorem, 
whose proof can be found in Appendix \ref{app_proof_equ_upper_bound_comp}.

\begin{lemma}
\label{equ_upper_bound_comp}
Consider a parameter vector $\xtrue \!\in\! \mathcal{X}_{S}$ with maximal support, i.e., ${\| \xtrue \|}_{0} \!=\! S$. Then, for any $k \notin \supp(\xtrue)$, 
the minimum MSE of any estimator $\hat{x}_{k}(\cdot) \ist\in\, \mathcal{U}^{k} \cap \mathcal{A}_{\xtrue}^{k}$, denoted by \emph{$\mbox{BB}_{\text{c}}^{k}(\xtrue)$}, is given by 
\emph{
\be
\label{equ_expr_bbc_comp1}
\mbox{BB}^{k}_{\text{c}}(\xtrue) 
%% \,\triangleq\, \mathsf{E}_{\xtrue} \{ \hat{x}_{k, \xtrue} ( \mathbf{y} ) |^{2} \} 
  \eq \bigg[ 1 -\! \prod_{l \ist\in\ist \supp(\xtrue) } \!\! g( x_l; \sigma^{2}) \bigg] \ist \sigma^{2} 
\vspace{-2mm}
\ee
\emph{with}}
\vspace{1mm}
\be
\label{equ_expr_bbc_g}
g(x;\sigma^{2}) \eq \frac{1}{\sqrt{2 \pi \sigma^{2}} } \int_0^\infty \! e^{-(x^2 + y^2 )/(2 \sigma^{2}) } \, {\rm sinh} \bigg( \frac{x y}{\sigma^{2}} \bigg) \, 
  {\rm tanh} \bigg( \frac{x y}{\sigma^{2}} \bigg) \ist dy \,.
\vspace{3.5mm}
\ee
\end{lemma}

Lemma \ref{equ_upper_bound_comp} identifies the minimum MSE of any unbiased estimator  %%% AJ 
of the $k$th component of $\mathbf{x}$ (where $k \notin \supp(\xtrue)$) that is also constrained to be an element of $\mathcal{A}^{k}_{\xtrue}$. Note that $\mbox{BB}_{\text{c}}^{k}(\xtrue)$ does not depend on $k$.
%% , apart from the constraint $k \notin \supp(\xtrue)$. 
It provides an upper bound on the minimum MSE of any unbiased estimator of the $k$th component of $\mathbf{x}$, for 
any $k \notin \supp(\xtrue)$.

The total MSE of a vector estimator $\hat{\mathbf{x}}(\cdot)$ can be decomposed as 
$\varepsilon(\xtrue;\hat{\mathbf{x}}) = \sum_{k \ist\in\ist \supp(\xtrue)}\varepsilon(\xtrue;\hat{x}_k) + \sum_{k \ist\notin\ist \supp(\xtrue)}\varepsilon(\xtrue;\hat{x}_k)$ with the component MSE
$\varepsilon(\xtrue;\hat{x}_k) \triangleq \mathsf{E}_{\xtrue} \big\{ ( \hat{x}_k( \mathbf{y}) - x_k )^{2} \big\}$.
Inserting the minimum component MSE for $k \in \supp(\xtrue)$ (which is $\sigma^2$ according to Lemma
%% Theorem 
\ref{equ_opt_k_in_supp}) in the first sum and the upper bound $\mbox{BB}^{k}_{\text{c}}(\xtrue)$ on the minimum component MSE for $k \notin \supp(\xtrue)$ 
%% (as expressed by \eqref{equ_expr_bbc_comp1}) 
in the second sum, we obtain the following upper bound on the minimum total MSE of any unbiased vector estimator.

%% \pagebreak %%%%%%%%%

\begin{theorem}
\label{equ_upper_bound}
The minimum MSE achievable by any unbiased estimator for the SSNM at a parameter vector $\xtrue \!\in\! \mathcal{X}_{S}$ with ${\| \xtrue \|}_{0} \rmv=\rmv S$ satisfies
\emph{
\be
\label{equ_upper_bound_bbc}
\varepsilon( \xtrue; \hat{\mathbf{x}}^{(\xtrue)}) \,\leq\, \mbox{BB}_{\text{c}}(\xtrue) \,\triangleq\, S \ist \sigma^{2} +\ist (N \!-\! S) \,\mbox{BB}^{k}_{\text{c}}(\xtrue)
\ee
}
with \emph{$\mbox{BB}^{k}_{\text{c}}(\xtrue)$} given by \eqref{equ_expr_bbc_comp1}.
%% where $k$ is an arbitrary element of $\{1,\ldots,N\} \setminus \supp(\xtrue)$. 
%The first component of the right hand side in \eqref{equ_upper_bound_bbc} is due to the MSE of estimating the components $\left( \mathbf{x} \right)_{l}$ with $l \in \supp(\xtrue)$ and the
%second component is due to the MSE of estimating the components $ \left( \mathbf{x} \right)_{l}$ with $l \notin \supp(\xtrue)$.\footnote{ALEX: It seems to me that the last sentence is an interpretation, not something that should appear as a theorem.}
\end{theorem}

Depending on the parameter vector $\mathbf{x}$, the upper bound $\mbox{BB}_{\text{c}}(\xtrue)$ varies between two extreme values. For decreasing SNR $\xi^{2}/ \sigma^2$, it converges to 
the low-SNR value $N \sigma^{2}$ (because the factor $g(\xi,\sigma^{2})$ in \eqref{equ_expr_bbc_comp1} vanishes for $\xi^{2}/ \sigma^2 \rightarrow 0$). On the other hand, we will show below that for increasing SNR, $\mbox{BB}_{\text{c}}(\xtrue)$ converges to its high-SNR value, which is given by $S \sigma^{2}$. 

The lower bound $\mbox{HCRB}(\xtrue)$ in \eqref{equ_HCR} for the case $ {\|\xtrue\|}_0 \rmv=\rmv S$, i.e., $S \ist \sigma^2 +\ist (N - S - 1) \ist\ist e^{-\xitrue^{2}/\sigma^{2}} \sigma^2 \!$,
%% derived in the previous section 
exhibits an exponential transition between the low-SNR and high-SNR regimes. More specifically, when considering a sequence 
%% ${\{ \mathbf{x}_{i} \}}_{i = 1,2,\ldots}$ 
%%% BEGIN ALEX 08042010
of parameter vectors $\mathbf{x} \!\in\! \mathcal{X}_{S}$ with increasing SNR $\xitrue^2/\sigma^2$,
the bound transitions from the low-SNR value $(N \!-\! 1) \ist \sigma^2$ (obtained for $\xitrue^2/\sigma^2 = 0$) to the high-SNR value $S \ist \sigma^{2}$ (obtained for $\xitrue^2/\sigma^2 \to \infty$); this transition is %via the factor $e^{-\xitrue^{2}/\sigma^{2}}\!$ and thus  %%% AJ 
exponential in the SNR. The upper bound $\mbox{BB}_{\text{c}}(\xtrue)$ in \eqref{equ_upper_bound_bbc} also exhibits a transition that is exponential in $\xitrue^2/\sigma^2$. 
%%% END ALEX 08042010
In fact, it is shown in Appendix \ref{app_proof_exponential} that
\be
\label{equ_exp_conv_bb}
%%% BEGIN ALEX 08042010
\mbox{BB}_{\text{c}}(\xtrue) \,\leq\, S \ist \sigma^2 +\ist (N \!-\! S) \, 3^{S} \ist e^{- \xi^{2}/(2 \sigma^{2})} \ist \sigma^2 .
%%% END ALEX 08042010
 % \qquad \text{with} \;\, c \,\triangleq\, \sqrt{\frac{2}{\pi \sigma^{2}}} + \frac{1}{2} \,.
\ee
This 
%% explicitly 
%%% BEGIN ALEX 08042010
shows that for increasing $\xitrue^2/\sigma^2$, the upper bound $\mbox{BB}_{\text{c}}(\xtrue)$---just like the lower bound $\mbox{HCRB}(\xtrue)$---decays  %%% AJ 
exponentially to its asymptotic value $S \ist \sigma^{2}$, which is also the asymptotic value of $\mbox{HCRB}(\xtrue)$.
%% the HCRB\@. 
%%% END ALEX 08042010
It follows
%% This implies 
that the BB itself also converges exponentially to   %%% AJ 
%% a value of 
%%% BEGIN ALEX 08042010
$S \ist \sigma^2$ as $\xitrue^2/\sigma^2$ increases. This result  will be further explored in Section \ref{sec.simu.role}.
%% The implications of this result will be explored in the next section. 
%%% END ALEX 08042010

\vspace{1mm}

%%%%%%%%%%%%%%%%%%%%%%%%%%%%%%%%%%%%%%%%%%%%%%%%%%%%%%%%%%
\section{Numerical Results}
\label{sec_sim}
%%%%%%%%%%%%%%%%%%%%%%%%%%%%%%%%%%%%%%%%%%%%%%%%%%%%%%%%%%

\vspace{1mm}

In this section, we describe several numerical studies which  %%% AJ 
%several experiments designed to further 
explore 
and extend %%% AJ
the theoretical bounds developed above. These include a numerical improvement of the bounds,
a comparison with practical (biased) estimation techniques, an analysis of the performance at high SNR, and an examination of the ability to 
%identify 
estimate %%% AJ 
the threshold region in which the transition from low to high SNR occurs.

We will first show that it is possible to obtain significantly tighter versions of the lower and upper bounds  developed in Sections \ref{sec_lowerbound} and \ref{sec_upperbound}. %%% AJ 
% section
%consider the tightness of the upper and lower bounds developed in the previous section. As we will see, it is possible to obtain significantly
%% somewhat 
%tighter versions of both the upper and lower bounds. 
These tightened versions can only be computed numerically and no longer have a simple form; consequently, they are less convenient for theoretical analyses. %%% AJ 
Nevertheless, they characterize the BB very accurately and therefore also provide an indication of the accuracy of the simpler, closed-form bounds.  %%% AJ 

\subsection{Numerical Lower Bound}
%%%%%%%%%%%%%%%%%%%%%%%%%%%%%%%%%%%%%%%%%%%%%%%%%%%%%%%%%%

\begin{figure}
\vspace{-1mm}
\centering
\psfrag{SNR}[c][c][1]{\uput{2.5mm}[270]{0}{SNR (dB)}}
\psfrag{Bounds}[c][c][1]{\uput{2.5mm}[270]{0}{}}
\psfrag{x_0}[c][c][1]{\uput{0.1mm}[270]{0}{$0$}}
\psfrag{x_0_001}[c][c][1]{\uput{0.1mm}[270]{0}{$-30$}}
\psfrag{x_0_01}[c][c][1]{\uput{0.1mm}[270]{0}{$-20$}}
\psfrag{x_0_1}[c][c][1]{\uput{0.1mm}[270]{0}{$-10$}}
\psfrag{x_1}[c][c][1]{\uput{0.1mm}[270]{0}{$0$}}
\psfrag{x_10}[c][c][1]{\uput{0.1mm}[270]{0}{$10$}}
\psfrag{y_1}[c][c][1]{\uput{0.1mm}[180]{0}{$1$}}
\psfrag{y_2}[c][c][1]{\uput{0.1mm}[180]{0}{$2$}}
\psfrag{y_3}[c][c][1]{\uput{0.1mm}[180]{0}{$3$}}
\psfrag{y_4}[c][c][1]{\uput{0.1mm}[180]{0}{$4$}}
\psfrag{y_5}[c][c][1]{\uput{0.1mm}[180]{0}{$5$}}
\psfrag{data2}[l][l][0.8]{$\mbox{BB}'_{\text{c}}(\xtrue)$}
\psfrag{data4}[l][l][0.8]{$\mbox{HCRB}(\xtrue)$}
\psfrag{data1}[l][l][0.8]{$\mbox{BB}_{\text{c}}(\xtrue)$}
\psfrag{data3}[l][l][0.8]{$\mbox{HCRB}_{\mathcal{V}} (\xtrue)$}
\centering
\includegraphics[height=6cm,width=12.3cm]{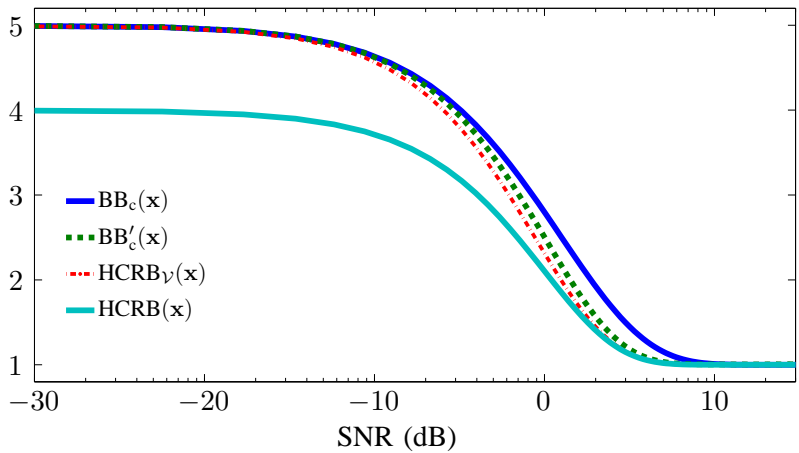}
  \caption{Lower bounds $\mbox{HCRB}(\xtrue)$, $\mbox{HCRB}_{\mathcal{V}}(\xtrue)$ and upper bounds $\mbox{BB}_{\text{c}}(\xtrue)$, $\mbox{BB}'_{\text{c}}(\xtrue)$
  on the MSE $\varepsilon(\xtrue;\hat{\mathbf{x}}^{(\xtrue)})$ of the LMVU estimator at $\xtrue = c \,(1 \;\, 0 \;\, 0 \;\, 0 \;\, 0)^T$, with $c$ varied to obtain different
values of $\mbox{SNR}(\xtrue) = \xi^{2}/\sigma^{2}$. The SSNM parameters are $N \!=\! 5$, $S \!=\! 1$, and $\sigma^{2} \!=\! 1$.}
\label{fig_bounds_1}
\end{figure}

For a parameter vector $\mathbf{x}$ with ${\| \mathbf{x} \|}_{0} \!=\! S$, let us reconsider the HCRB in \eqref{equ hcr}.
%%  in Proposition \ref{pr:hcrb}. %%% AJ 
We will show that by using an increased number of appropriately chosen test points, we can obtain a lower bound that is higher (thus, tighter) than \eqref{equ_HCR}.
Specifically, assume without loss of generality that $\supp(\xtrue) = \{1,\ldots,S \}$, and consider the set of test points 
\[
\mathcal{V} \,\triangleq\, \mathcal{V}_{0}  \cup\, \bigcup_{k=1}^S (\mathcal{V}_{k} \cup \mathcal{W}_{k})
\vspace{-2mm}
\]
with the component sets
\begin{align*}
\mathcal{V}_{0} & \,\triangleq\! \bigcup_{l \in \, \supp(\xtrue)}  \!\{ \alpha \ist\ist\mathbf{e}_{l} \}  \\
\mathcal{V}_{k} & \,\triangleq\! \bigcup_{l \in \{S+1,...,N\}} \!\{  \alpha \ist\ist\mathbf{e}_{l} \rmv-\rmv x_k \ist \mathbf{e}_{k}  \} \,, \qquad k=1,\ldots,S \\
\mathcal{W}_{k} & \,\triangleq\! \bigcup_{l \in \{S+1,...,N\}} \!\{  x_k \ist \ist\mathbf{e}_{l} \rmv-\rmv x_k \ist \mathbf{e}_{k}  \} \,, \qquad k=1,\ldots,S\\[-9mm]
&
\end{align*}
where $\alpha = 0.02 \ist\ist \sigma$.
In Fig.~\ref{fig_bounds_1}, the HCRB \eqref{equ hcr} for the new set $\mathcal{V}$ of test points---denoted $\mbox{HCRB}_{\mathcal{V}}(\xtrue)$---is displayed versus the SNR and compared with $\mbox{HCRB}(\mathbf{x})$.
%%  in \eqref{equ_HCR}.
%%  of Theorem~\ref{th:hcrb}. 
For this figure, we chose $NÊ\! =\! 5$, $S \! =Ê\! 1$, $\sigma^{2} \!=\! 1$, and $\mathbf{x} = c \,(1 \;\, 0 \;\, 0 \;\, 0 \;\, 0)^T$, where the parameter $c \! \in \! \mathbb{R}$ is varied to obtain different SNR  %%% AJ 
values.\footnote{The %%%%%%%%%%
use of a low-dimensional model is mandated by the complexity of the numerical approximation to the upper bound on the BB which will be described 
in Section \ref{sec_num_upper_bound}.} %%%%%%%%%%  AJ 
As before, the SNR is defined as
%% The horizontal axes of Fig. \ref{fig_bounds_1} indicates the SNR which we define as
%% \label{equ_SNR_def}
$\mbox{SNR}(\xtrue) = \xi^{2}/\sigma^{2}$, where $\xi$ is the $S$-largest (in magnitude) component of $\x$ (in our example with $S \! = \! 1$, $\xi$ is simply  %%% AJ 
the single nonzero component).
%% The justification for this definition of the SNR will be given presently.
It can be seen from Fig.~\ref{fig_bounds_1} that the numerical lower bound $\mbox{HCRB}_{\mathcal{V}}(\xtrue)$ computed from the above test points 
%% \eqref{equ_V'} 
is indeed tighter than the closed-form lower bound $\mbox{HCRB}(\mathbf{x})$ in \eqref{equ_HCR}.
%% obtained in the previous section.

\subsection{Numerical Upper Bound}
%%%%%%%%%%%%%%%%%%%%%%%%%%%%%%%%%%%%%%%%%%%%%%%%%%%%%%%%%%
\label{sec_num_upper_bound}

It is also possible to find upper bounds on the BB that are tighter (lower) than the upper bound $\mbox{BB}_{\text{c}}(\xtrue)$ in \eqref{equ_upper_bound_bbc}. 
Consider a parameter vector %%% AJ 
$\xtrue$ with ${\|\xtrue\|}_{0} \! = \! S$. We recall that $\mbox{BB}_{\text{c}}(\xtrue)$ was derived by constructing, for all %%% AJ 
%% components $x_{k}$ with 
$k \notin \supp(\xtrue)$, unbiased estimators $\hat{x}_{k}(\mathbf{y}) = y_{k} + \hat{x}'_{k}(\mathbf{y})$ with $\hat{x}'_{k} (\mathbf{y})$ constrained by \eqref{equ_constr_odd_symm} and \eqref{equ_constr_invariance}. We will now investigate how much we can improve on $\mbox{BB}_{\text{c}}(\xtrue)$ if we remove
%% gives up 
the constraint \eqref{equ_constr_odd_symm}. Thus, in the optimization problem \eqref{equ_opt_approx_bb}, the constraint set $\mathcal{A}_{\xtrue}$ is hereafter considered to correspond only to the constraint \eqref{equ_constr_invariance}. 

In order to numerically solve this modified optimization problem \eqref{equ_opt_approx_bb},
%% compute the optimum unbiased estimator satisfying only the constraint \eqref{equ_constr_invariance}, 
a discrete approximation for $\hat{x}'_k(\y)$ was used.
%% constructed.
More specifically, we defined $\hat{x}'_k(\y)$ to be piecewise constant in each of the components $y_l$ with $l \in \{k\}  \cup \supp(\x) $, and constant in the remaining components $y_l$ (the latter being required by \eqref{equ_constr_invariance}). 
%% In particular, we 
We used $Q$ piecewise constant segments for each
%% in each axis 
$l \in \{k\} \cup \supp(\x)$, with each segment of length $\Delta \! = \! 10 \, \sigma/Q$. These arrays of constant segments were centered about $\y \! = \! \x$. %%% AJ 
%% These constant segments were centered on $\y=\x$. 
The remaining values of $\hat{x}'_k(\y)$ were set to $0$. Thus, we obtained a function $\hat{x}'_k(\y)$ 
%% which linearly depends 
with linear dependence on a finite number $Q^{S+1}$ of parameters. For functions of this form, the optimization problem \eqref{equ_opt_approx_bb} becomes a finite-dimensional quadratic program with linear constraints, 
%% a well-studied problem 
which can be solved efficiently \cite{BoydConvexBook}. The MSE of the resulting estimator, denoted by $\mbox{BB}'_{\text{c}}(\xtrue)$, is an upper bound on the BB. This bound is tighter than the closed-form upper bound $\mbox{BB}_{\text{c}}(\xtrue)$ in \eqref{equ_upper_bound_bbc} if $Q$ is large enough.  %%% AJ 
In Fig.~\ref{fig_bounds_1}, we compare $\mbox{BB}'_{\text{c}}(\xtrue)$ for $Q \! =\! 20$ with $\mbox{BB}_{\text{c}}(\xtrue)$ as a function of the SNR. The improved accuracy of $\mbox{BB}'_{\text{c}}(\xtrue)$ relative to $\mbox{BB}_{\text{c}}(\xtrue)$ is evident, particularly at high SNR values. Moreover, the proximity of the numerical upper bound $\mbox{BB}'_{\text{c}}(\xtrue)$ to the numerical lower bound $\mbox{HCRB}_{\mathcal{V}}(\xtrue)$ indicates that these two bounds achieve an accurate characterization of the BB, since the BB lies between them. %%% AJ 

%% \pagebreak %%%%%%%%%

\subsection{The Role of $\xi$}
\label{sec.simu.role}
%%%%%%%%%%%%%%%%%%%%%%%%%%%%%%%%%%%%%%%%%%%%%%%%%%%%%%%%%%

We have seen in Section \ref{sec_upperbound} that for ${\|\xtrue\|}_{0} \! = \! S$, the MSE %%% AJ 
of the LMVU estimator at high SNR is given by $S \sigma^2$, and furthermore, convergence to this value is exponential in the quantity $\xi^2/\sigma^2$. A
%% The most 
remarkable aspect of this conclusion is the fact that convergence to the high-SNR regime depends solely on $\xi$, the smallest nonzero component of $\x$, rather than having a more complex dependency on all the $S$ nonzero components of $\x$. For example, one might imagine the behavior of an estimator to be rather different when all nonzero components have the same value $\xi$, as opposed to the situation in which one component equals $\xi$ and the others are much larger. However, our analysis shows that when $\xi \rmv\rmv\gg\rmv\rmv \sigma$, the remaining components of $\x$ have no effect on the performance of the LMVU estimator. We will next investigate whether practical estimators also exhibit such an effect. 
%% In other words, when the noise level is low, is practical performance effectively 
%% determined by $\xi$, or is the behavior of actual estimators more complex?

\begin{figure}
\centering
\psfrag{x}[c][c][1]{\uput{3.5mm}[270]{0}{$r$}}
\psfrag{y}[c][c][1][270]{\uput{0.8mm}[180]{0}{$\varepsilon(\mathbf{x}_{r}; \ML_est)$}}
\psfrag{x_10}[c][c][1]{\uput{0.2mm}[270]{0}{$10$}}
\psfrag{x_50}[c][c][1]{\uput{0.2mm}[270]{0}{$50$}}
\psfrag{x_100}[c][c][1]{\uput{0.2mm}[270]{0}{$100$}}
\psfrag{y_4}[c][c][1]{\uput{0.1mm}[180]{0}{$4$}}
\psfrag{y_6}[c][c][1]{\uput{0.1mm}[180]{0}{$6$}}
\psfrag{y_8}[c][c][1]{\uput{0.1mm}[180]{0}{$8$}}
\psfrag{y_10}[c][c][1]{\uput{0.1mm}[180]{0}{$10$}}
\psfrag{SNR=4}[l][l][0.7]{$\mbox{SNR}=4$}
\psfrag{SNR=9}[l][l][0.7]{$\mbox{SNR}=9$}
\psfrag{SNR=16}[l][l][0.7]{$\mbox{SNR}=16$}
\psfrag{SNR=25}[l][l][0.7]{$\mbox{SNR}=25$}
\centering
\hspace*{20mm}\includegraphics[height=6cm,width=15cm]{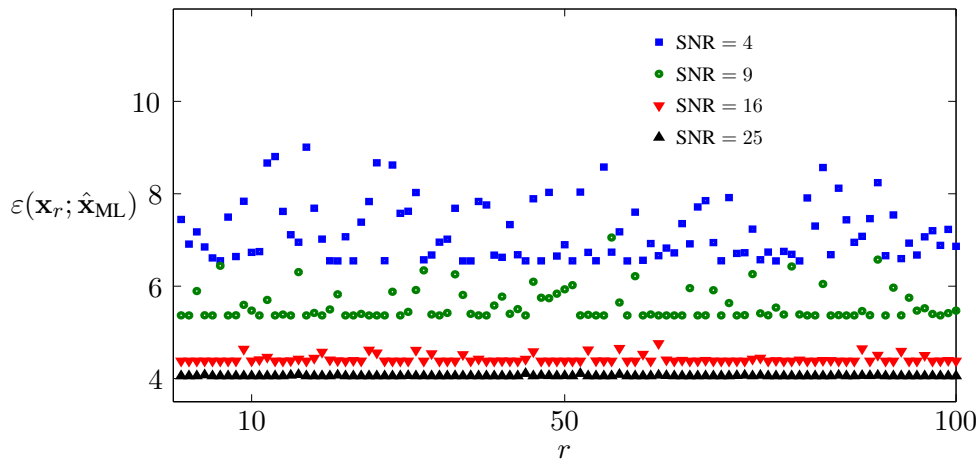}
\caption{MSE $\varepsilon(\mathbf{x}_{r}; \ML_est)$ of the ML estimator for randomly generated parameter vectors $\mathbf{x}_{r}$ at four different SNRs $\xi^2/\sigma^2$,
for SSNM parameters $N \!=\! 10$, $S \!=\! 4$, and $\sigma^2 \!=\! 1$.} %%% AJ
\label{fig_ml_monte_carlo}
\end{figure}

%% \pagebreak %%%%%%%%%

To answer this question, we examined the MSE of the ML estimator \eqref{equ_LS_ML_coincide}
%% \eqref{equ_ML_est} 
for a wide range of parameter vectors $\x$ having a predetermined smallest component $\xi$. More specifically, for a given value of $\xi$, we randomly generated $100$ parameter vectors $\mathbf{x}_{r}$, $r=1,\ldots,100$, with $\mathbf{x}_{r} \!\!\in\! \mathcal{X}_{S}$ and ${\| \mathbf{x}_{r} \|}_{0}\!=\!S$, whose minimum nonzero component was equal to $\xi$. The other nonzero components were generated 
%% randomly 
as independent, identically distributed realizations of the random variable $x = \xi \ist (1 + 3 \ist \sigma |q|)$, where $q \sim \mathcal{N}(0,1)$ is a standard Gaussian random variable and $\sigma$ is the standard deviation of the noise. The MSE $\varepsilon(\mathbf{x}_{r}; \ML_est)$ of the ML estimator is shown in Fig.~\ref{fig_ml_monte_carlo} for $N \!=\! 10$, $S \!=\! 4$, and four different SNRs $\xi^2/\sigma^2$, with the horizontal axis representing the different choices of $\mathbf{x}_{r}$ in arbitrary order. %%% AJ
It is seen that for large $\xi$, the performance of the ML estimator, like that of the LMVU, depends almost exclusively on $\xi$. This suggests that the performance guarantees of Sections \ref{sec_lowerbound} and \ref{sec_upperbound}, %%% AJ 
while formally valid only for unbiased estimators, can still provide general conclusions which are relevant to biased techniques such as the ML estimator. Moreover, this result also justifies our definition 
%% \eqref{equ_SNR_def} 
of the SNR as the ratio $\xi^2/\sigma^2$, since this 
%% measure 
is the most significant factor determining estimation performance for the SSNM\@.

\subsection{Threshold Region Identification}
%%%%%%%%%%%%%%%%%%%%%%%%%%%%%%%%%%%%%%%%%%%%%%%%%%%%%%%%%%

In Sections~\ref{sec_lowerbound} and ~\ref{sec_upperbound}, we characterized the performance of unbiased estimators as a means of quantifying the difficulty of estimation for the SSNM\@. A common use of  %%% AJ 
this analysis is in the identification of the threshold region, a range of SNR values which constitutes %%% AJ
a transition between low-SNR and high-SNR behavior \cite{LowerBoundAbel,ForsterICASSP2002,TodrosICASSP2008}. %%% AJ 
Specifically, in many cases the performance of estimators can be calculated analytically when the SNR is either very low or very high. It is then important to identify the threshold 
region which separates these two regimes. Although the analysis is based on bounds for unbiased estimators, the result is often heuristically assumed to approximate the threshold region for biased techniques as well \cite{LowerBoundAbel,TodrosICASSP2008}.  %%% AJ

For ${\| \xtrue \|}_{0} \! = \! S$, the lower and upper bounds on the BB 
%% which were derived above 
($\text{HCRB}(\xtrue)$ in \eqref{equ_HCR}, $\text{BB}_{\text{c}}(\xtrue)$ in \eqref{equ_upper_bound_bbc}) 
%% also 
exhibit a transition between a low-SNR region, where both bounds are on the order of $N \sigma^2\rmv$, and a high-SNR region, for which both bounds converge to $S \sigma^2\rmv$. The BB %%%AJ
therefore also displays such a transition. One can define the threshold region of the SSNM (for unbiased estimation) as the range of values of $\mathbf{\xi}^{2}/\sigma^{2}$ in which this transition takes place. Since the BB is itself a lower bound on the performance of unbiased estimators, one would expect the transition region of 
actual estimators to occur at slightly higher SNR values than that of the BB. 
%% of the BB to occur at slightly lower SNR values than that of actual estimators.

%% \pagebreak %%%%%%%%%

To test this hypothesis, we compared the bounds of Sections~\ref{sec_lowerbound} and ~\ref{sec_upperbound} with the MSE of two well-known estimation schemes, namely, the ML estimator in \eqref{equ_LS_ML_coincide} and the hard-thresholding (HT) estimator 
$\hat{\mathbf{x}}_{\text{HT}}(\mathbf{y})$, which is given componentwise as
\[
%% \label{equ_threshold_est}
\hat{x}_{\text{HT},k}(\mathbf{y}) \eq \begin{cases}
 y_k ,  & 
 %% \text{if } 
 |y_k| \geq T \\[-1mm]
 0, & \text{else}
\end{cases}
\]
for a given threshold $T \!>\!0$. In our simulations, we chose the commonly used value $T = \sigma \sqrt{2 \log{N}}$ \cite{Mallat98}.
%\footnote{ALEX: Isn't it strange that hard thresholding, which is asymptotically optimal, is asymptotically worse than the ML estimator, which is also hard thresholding, presumably with a different threshold?}
Note that since the ML and HT estimators are biased, their MSE is not bounded by $\mbox{BB}_{\text{c}}(\xtrue)$, $\mbox{HCRB}(\xtrue)$, and the CRB. Assuming SSNM parameters $N \!=\! 10$ and $S \!= 4$, we generated a number of parameter vectors $\x$ from the set 
$\mathcal{R} \triangleq \big\{ c \, (1\;\, 1\;\, 1\;\, 1\;\, 0\;\, 0\;\, 0\;\, 0\;\, 0\;\, 0)^{T} \big\}_{c \in \mathbb{R}}$, where $c$ was varied to obtain a range of SNR values. For these $\x$, we calculated the MSE of the 
two estimators $\ML_est$ and $\hat{\mathbf{x}}_{\text{HT}}$  %%% AJ 
by means of numerical integration (see Appendix \ref{app_MSE_ML} for a discussion of the 
%% somewhat more involved 
computation of $\varepsilon(\mathbf{x};\ML_est)$).
%% The somewhat more involved computation of $\varepsilon(\x;\ML_est)$ is discussed in Appendix \ref{app_MSE_ML}.

\begin{figure*}
\vspace{-1mm}
\centering
\psfrag{SNR}[c][c][1]{\uput{2.5mm}[270]{0}{SNR (dB)}}
\psfrag{MSE}[c][c][1.2][270]{\uput{0.8mm}[180]{0}{$\varepsilon$}}
\psfrag{x_0_1}[c][c][1]{\uput{0.1mm}[270]{0}{$-10$}}
\psfrag{x_0_01}[c][c][1]{\uput{0.1mm}[270]{0}{$-20$}}
\psfrag{x_1}[c][c][1]{\uput{0.1mm}[270]{0}{$0$}}
\psfrag{x_10}[c][c][1]{\uput{0.1mm}[270]{0}{$10$}}
\psfrag{x_100}[c][c][1]{\uput{0.1mm}[270]{0}{$20$}}
\psfrag{y_2}[c][c][1]{\uput{0.1mm}[180]{0}{$2$}}
\psfrag{y_12}[c][c][1]{\uput{0.1mm}[180]{0}{$12$}}
\psfrag{y_14}[c][c][1]{\uput{0.1mm}[180]{0}{$14$}}
\psfrag{y_4}[c][c][1]{\uput{0.1mm}[180]{0}{$4$}}
\psfrag{y_6}[c][c][1]{\uput{0.1mm}[180]{0}{$6$}}
\psfrag{y_8}[c][c][1]{\uput{0.1mm}[180]{0}{$8$}}
\psfrag{y_10}[c][c][1]{\uput{0.1mm}[180]{0}{$10$}}
\psfrag{data4}[l][l][0.8]{$\varepsilon_{\mathrm{ML}}$}
\psfrag{data2}[l][l][0.8]{HCRB}
\psfrag{data1}[l][l][0.8]{$\mbox{BB}_{\text{c}}$}
\psfrag{data3}[l][l][0.8]{$\varepsilon_{\mathrm{HT}}$}
\psfrag{data5}[l][l][0.8]{CRB}
\centering
\includegraphics[height=6cm,width=12.3cm]{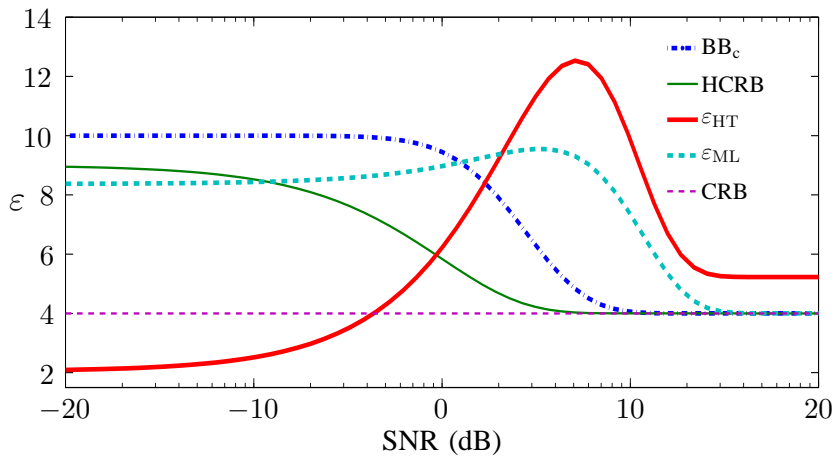}
  \caption{MSE of the ML and HT estimators compared with the 
  performance bounds $\mbox{BB}_{\text{c}}(\xtrue)$, $\mbox{HCRB}(\xtrue)$, and CRB ($\equiv S \sigma^{2}$), as a function of the SNR $\xi^2/\sigma^2$, for SSNM parameters $N \!=\! 10$, $S \!=\! 4$, and $\sigma^{2} \!=\! 1$.} %%% AJ 
\label{fig_all_equal}
\end{figure*}

The results are displayed in Fig.~\ref{fig_all_equal} as a function of the SNR $\xi^2/\sigma^2$. %corresponding to the respective $\x$. %%% AJ
Although there is some gap between the lower bound (HCRB) and the upper bound ($\mbox{BB}_{\text{c}}$), a rough indication of the behavior of the BB is conveyed.
As expected, the threshold region exhibited by the ML and HT estimators is somewhat higher than that predicted by the bounds.
%% the threshold predicted by these bounds is somewhat lower than that obtained by practical estimators. 
Specifically, the threshold region of the BB (as indicated by the bounds) can be seen to occur at SNR values between $-5$ and $5$ dB, while the threshold 
region of the ML and HT estimators is at SNR values between $5$ and $12$ dB\@. %%% AJ
Another effect which is visible in Fig.~\ref{fig_all_equal} is the convergence of the ML estimator to the BB at high SNR; this is a manifestation %%% AJ 
of the well-known fact that the ML estimator %%% AJ 
is asymptotically unbiased and asymptotically optimal. 
Finally, at low SNR, both the ML and HT estimators are better than the best unbiased approach.  %%% AJ 
This is because unbiased methods generally perform poorly at low SNR, so that even the best unbiased technique %%% AJ 
is outperformed by the biased ML and HT estimators.
%\footnote{ALEX: In my opinion the figure showing the asymptotic rate of convergence, and its discussion (both commented, below this footnote), can be removed. I think that Fig.~\ref{fig_ml_monte_carlo} demonstrates everything we want to say about asymptotic behavior, and the discussion doesn't have too much of a punchline. Also the results are not very good in terms of our bounds' performance.}
On the other hand, for medium %%%AJ
SNR, the MSE of the ML and HT estimators is significantly higher than the BB. Thus, there is a potential for unbiased estimators to perform better than biased estimators in the medium-SNR regime. 

\begin{figure*}[t]
\vspace{-1mm}
\centering
\psfrag{SNR}[c][c][1]{\uput{3.5mm}[270]{0}{SNR (dB)}}
\psfrag{y}[c][c][1][270]{\uput{0.8mm}[180]{0}{}}
\psfrag{x_0_01}[c][c][1]{\uput{0.2mm}[270]{0}{$-20$}}
\psfrag{x_0_1}[c][c][1]{\uput{0.2mm}[270]{0}{$-10$}}
\psfrag{x_1}[c][c][1]{\uput{0.2mm}[270]{0}{$0$}}
\psfrag{x_10}[c][c][1]{\uput{0.2mm}[270]{0}{$10$}}
\psfrag{x_100}[c][c][1]{\uput{0.2mm}[270]{0}{$20$}}
\psfrag{x_1000}[c][c][1]{\uput{0.2mm}[270]{0}{$30$}}
\psfrag{y_1}[c][c][1]{\uput{0.1mm}[180]{0}{$1$}}
\psfrag{y_1_2}[c][c][1]{\uput{0.1mm}[180]{0}{$1.2$}}
\psfrag{y_1_4}[c][c][1]{\uput{0.1mm}[180]{0}{$1.4$}}
\psfrag{y_1_6}[c][c][1]{\uput{0.1mm}[180]{0}{$1.6$}}
\psfrag{y_1_8}[c][c][1]{\uput{0.1mm}[180]{0}{$1.8$}}
\psfrag{data2}[l][l][0.7]{$\mathcal{R}_{2}$}
\psfrag{data1}[l][l][0.7]{$\mathcal{R}$}
\psfrag{data3}[l][l][0.7]{$\mathcal{R}_{3}$}
%% \psfrag{(a)}[c][c][1]{\uput{0.1mm}[270]{0}{(a)}}
%% \psfrag{(b)}[c][c][1]{\uput{0.1mm}[270]{0}{(b)}}
\centering
\hspace*{5mm} \includegraphics[height=6cm,width=12.3cm]{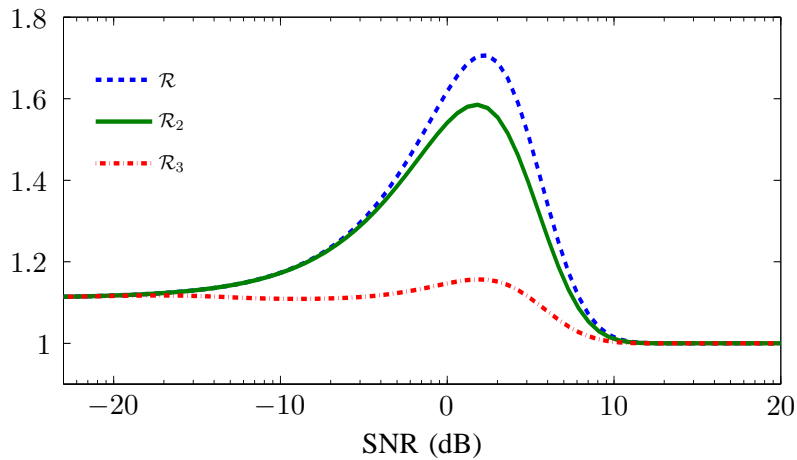}
%%     \hspace*{-2mm}\resizebox{9cm}{!}{\includegraphics{bild-3.eps}}
  \caption{Ratio $\mbox{BB}_{\text{c}}(\x) / \mbox{HCRB}(\x)$ versus the SNR $\xi^{2}/\sigma^{2}$ for different sets of parameter vectors $\x$.}  %%% AJ
%%  for CS-based channel estimation and least-squares channel estimation.}
\label{fig_tightness}
%% \vspace{-1mm}
\end{figure*}

One may argue that considering only parameter vectors $\x$ in the set $\mathcal{R}$ is
not representative, since $\mathcal{R}$ covers only a small part of the parameter space $\mathcal{X}_S$. However,
the choice of $\mathcal{R}$ is conservative in that the maximum deviation between $\mbox{HCRB}(\xtrue)$ and $\mbox{BB}_{\text{c}}(\xtrue)$ is largest
when the nonzero entries of $\x$ have approximately the same magnitude, which is the case for each element of $\mathcal{R}$.
This is illustrated in Fig.~\ref{fig_tightness}, which
shows the ratio between the two bounds versus the SNR $\xi^2/\sigma^2$ for three different configurations of the nonzero entries in the parameter vector. Specifically, we considered the two additional sets $\mathcal{R}_{2} \triangleq \big\{ c \, (10\;\, 1\;\, 1\;\, 1\;\, 0\;\, 0\;\, 0\;\, 0\;\, 0\;\, 0)^{T} \big\}_{c \in \mathbb{R}}$ 
%% \{ c \, (10,1,1,1,0,0,0,0,0,0)^{T} \}_{c \in \mathbb{R}_{+}}$ 
and $\mathcal{R}_{3} \triangleq \big\{ c \, (0.1\;\, 1\;\, 1\;\, 1\;\, 0\;\, 0\;\, 0\;\, 0\;\, 0\;\, 0)^{T} \big\}_{c \in \mathbb{R}}$, 
%% \{ c \, (0.1,1,1,1,0,0,0,0,0,0)^{T} \}_{c \in \mathbb{R}_{+}}$,
in which the nonzero entries have different magnitudes. It can be seen from Fig.~\ref{fig_tightness} that the ratio $\mbox{BB}_{\text{c}}(\xtrue) / \mbox{HCRB}(\xtrue)$
is indeed highest when $\x$ is in $\mathcal{R}$.

%%%%%%%%%%%%%%%%%%%%%%%%%%%%%%%%%%%%%%%%%%%%%%%%%%%%%%%%%%
\section{Conclusion}
%%%%%%%%%%%%%%%%%%%%%%%%%%%%%%%%%%%%%%%%%%%%%%%%%%%%%%%%%%

In this paper, %%% AJ 
we have studied unbiased estimation of a sparse vector in white Gaussian noise within a frequentist setting. As we have seen, without the assumption of sparsity, there exists only a single unbiased estimator.
%% , so that unbiased bounds in this case have little practical use. 
However, the addition of a sparsity assumption yields a rich family of unbiased estimators. The analysis of the performance of these estimators has been the primary goal of this paper. We first demonstrated that there exists no uniformly minimum variance unbiased estimator, i.e., no single unbiased estimator is optimal for all parameter values. Consequently, we focused on analyzing the Barankin bound (BB), i.e., the MSE of the locally minimum variance unbiased estimator, or equivalently, the smallest MSE achievable by an unbiased estimator for each 
%% specific 
value of the sparse vector.
%%  to be estimated.

For the sparse estimation problem considered, as for most estimation problems, the BB cannot be computed precisely. However, we demonstrated that it
%% the BB 
can be characterized quite accurately using numerical lower and upper bounds. Furthermore, we derived simple closed-form lower and upper bounds 
%% on the BB 
which are somewhat looser than the numerical bounds. These closed-form bounds allow an estimation
%% identification 
of the threshold region separating the low-SNR and high-SNR regimes, and they indicate the asymptotic behavior of the BB at high SNR.
In particular, a notable conclusion is that the high-SNR behavior of the BB depends solely on the value of the smallest nonzero component of the sparse vector.

While the unbiasedness property is intuitively appealing and 
%% unbiasedness is an intuitively appealing property which is 
related to several desirable asymptotic features of an estimator \cite{LC},
%% . However, there is some degree of arbitrariness in the examination of unbiased estimators, especially in light of the fact that 
one can often obtain biased estimators which outperform any unbiased estimator \cite{1956_Stein,BlindMinimaxZvika,RethinkingBiasedEldar}. Thus, it is 
%% particularly 
interesting to note that some of the 
%% general 
conclusions obtained from our analysis of unbiased estimators 
%% for the sparse model 
appear to provide insight into the behavior of standard biased estimators. In particular, we saw that the behavior of two commonly used biased estimators at high SNR corresponds to the predictions of our unbiased bounds, not only in terms of the asymptotically achievable MSE but also in certain finer details, such as the SNR range of the threshold region and the fact that the convergence to the high-SNR regime depends primarily on the value of the smallest nonzero component of the sparse vector,
%%  $\xi$, 
rather than on the entire vector. This gives additional merit to the analysis of achievable estimation performance within the unbiased setting.

%%%%%%%%%%%%%%%%%%%%%%%%%%%%%%%%%%%%%%%%
\appendices

%% \pagebreak %%%%%%%%%

%%%%%%%%%%%%%%%%%%%%%%%%%%%%%%%%%%%%%%%%
\section{Proof of Theorem \ref{thm_exit_uniqu_unbiasedness_S_N}}
\label{app_proof_ thm_exit_uniqu_unbiasedness_S_N}
%%%%%%%%%%%%%%%%%%%%%%%%%%%%%%%%%%%%%%%%

\vspace{1mm}

We wish to show that for $S \!=\! N$, the only unbiased estimator with bounded MSE is the trivial estimator $\hat{\mathbf{x}} (\mathbf{y}) \rmv=\rmv \mathbf{y}$.
We will first show that a bounded MSE implies that $\hat{\mathbf{x}} (\mathbf{y})$ is equivalent to a tempered distribution. This will allow
us to reformulate the unbiasedness condition in the Fourier transform domain.  

Using \eqref{equ_pdf}, the unbiasedness condition in \eqref{equ_unbiased_cond1} for $S \!=\! N$ 
\vspace{1mm}
reads
%% can be written as
\be
\label{cond_conv_unbiased}
%% \mathbf{b}(\xtrue) \triangleq 
\normconstgauss \int_{\mathbb{R}^N} \!\hat{\mathbf{x}}(\mathbf{y}) \, \exp\rmv\rmv \bigg( \!\!-\frac{1}{2 \sigma^{2}} \| \mathbf{y} \!-\! \xtrue \|_2^{2} \bigg)
%% e^{-\frac{1}{2 \sigma^{2}} \| \mathbf{y} - \xtrue \|^{2} } 
\, d \mathbf{y} \,=\,
%% \stackrel{!}{=} 0 
\xtrue
 \qquad \text{for all }\, \xtrue \!\in\! \mathbb{R}^{N} \,.
\vspace{1mm}
\ee
The integral in \eqref{cond_conv_unbiased} is
%% can be viewed as a 
the convolution of $\hat{\mathbf{x}}(\mathbf{y})$ with 
$\exp\rmv\rmv \big( \!\! -\!\rmv \frac{1}{2 \sigma^{2}} \| \mathbf{y} \|_2^{2} \big)$.
%% $e^{-\frac{1}{2 \sigma^{2}} \| \mathbf{y} \|^{2} }$. 
%% For \eqref{cond_conv_unbiased} to hold, the 
The result of this convolution, viewed as a function of $\xtrue$, must equal $(2\pi\sigma^2)^{N/2}\, \xtrue$ for all parameter vectors $\xtrue$.
For absolutely integrable functions, the Fourier transform maps a convolution onto a pointwise product, and consequently it seems natural to consider the Fourier transform of condition \eqref{cond_conv_unbiased} in order to simplify the analysis.
However, typically, the estimator function $\hat{\mathbf{x}}(\mathbf{y})$ will be neither absolutely integrable nor square  %%% AJ
integrable, and thus its Fourier transform can only exist in the sense of a tempered distribution \cite{Strichartz}.
From a practical point of view, the class of tempered distributions is large enough 
so
%% in 
that it does not exclude reasonable estimators such as the LS  estimator \eqref{equ_LS_ML_coincide}. The following lemma
states that $\hat{\mathbf{x}}(\mathbf{y})$ can be viewed as a tempered distribution if it has a bounded MSE. 
%% A sufficient condition ensuring that an estimator can be viewed as a tempered distribution is provided in the following lemma.

\begin{lemma}
\label{thm_bounded_tempered}
Consider an estimator $\hat{\mathbf{x}}$ for the SSNM \eqref{equ_ssnm} with $S \!=\! N$. If $\hat{\mathbf{x}}$ has a bounded MSE, i.e., $\MSE \leq C$ for all $\xtrue \!\in\! \mathbb{R}^{N}$ (where $C$ is a constant which may depend on 
%% the parameters
$N$, $S$, and $\sigma^{2}$),
%% ($N$,$S$,$\sigma^{2}$), 
then $\hat{\mathbf{x}}$ is equivalent to a tempered distribution.
\end{lemma}

\begin{bproof}
%% [Proof of Lemma \ref{thm_bounded_tempered}]
%% The proof of Lemma \ref{thm_bounded_tempered} is provided in Appendix \ref{app_proof_bounded_mse_temp}.
%% We have to show that an estimator $\hat{\mathbf{x}}$ with bounded MSE is equivalent to a tempered distribution. %%% AJ 
The proof of Lemma \ref{thm_bounded_tempered} is based on the following result which gives a sufficient condition for a function $\hat{\mathbf{x}}(\mathbf{y})$ to be (equivalent to) a tempered distribution. 

\begin{proposition}[\hspace*{-1.5mm}\cite{Strichartz}]
\label{prop_tempered}
%% Given a function $\hat{\mathbf{x}}(\mathbf{y})$, if 
If there exist 
%% some fixed 
constants $B,n,R_{0} \!\in\rmv\rmv \mathbb{R}_{+}$ such that
\be
\label{equ_suff_cond_bounded}
\int_{ {\| \mathbf{y} \|}_2 \leq R} \| \hat{\mathbf{x}}(\mathbf{y}) \|_2^{2} \, d \mathbf{y} \,\leq\,  B R^{n} \quad \quad \text{for all }\, R \rmv \geq \! R_{0}
\ee
then 
%% the estimator function  
$\hat{\mathbf{x}}(\mathbf{y})$ is equivalent to a tempered distribution.
\end{proposition}

%% The proof shows that a bounded MSE implies that  $\hat{\mathbf{x}}(\mathbf{y})$ satisfies the 
%% condition \eqref{equ_suff_cond_bounded} which (according to Proposition \ref{prop_tempered}) is sufficient 
%% for being equivalent to a tempered distribution. 

%% \begin{proof}[Proof of Lemma \ref{thm_bounded_tempered}]
%%% BEGIN ALEX 08042010
Let $\hx(\y)$ be an estimator function with
%%% END ALEX 08042010
%% having a 
bounded MSE,
%%  function, 
i.e., there exists a constant $C$ such that
%% $\mathsf{E}_{\xtrue} \{ \|\hat{\mathbf{x}}(\mathbf{y}) - \xtrue \|^{2}  \}  \leq C$ for all $\xtrue \in \mathcal{X}_{S}$,
\be
\label{equ_bounded_MSE}
\mathsf{E}_{\xtrue} \{ \|\hat{\mathbf{x}}(\mathbf{y}) \rmv - \rmv \xtrue \|_2^{2}  \}  \,\leq\, C  \quad \quad \text{for all }\, \xtrue \! \in \! \mathcal{X}_{S} \,.
\ee
%% which implies that
%% $\sqrt{\mathsf{E}_{\xtrue} \{ \|\hat{\mathbf{x}}(\mathbf{y}) - \xtrue \|^{2}  \} }  \leq \sqrt{C }$.
%% %% \sqrt{\mathsf{E}_{\xtrue} \{ \|\hat{\mathbf{x}}(\mathbf{y}) - \xtrue \|^{2}  \} }  \leq \sqrt{C }.
Defining the usual norm ${\| \cdot \|}_{\rm RV}$ on the space of of random vectors by $\| \mathbf{y} \|_{\rm RV} \triangleq  \sqrt{\mathsf{E}_{\xtrue} \{ \| \mathbf{y} \|^{2}_{2} \}}$, we can
use the (reverse) triangle inequality ${\| \hat{\mathbf{x}}(\mathbf{y}) \!  - \! \xtrue \|}_{\rm RV} \ge {\| \hat{\mathbf{x}}(\mathbf{y}) \|}_{\rm RV}- {\| \xtrue \|}_{\rm RV}$
%% for ${\| \cdot \|}_{\rm RV}$ 
to 
\vspace{1mm}
obtain
\[
\sqrt{\mathsf{E}_{\xtrue} \{ \|\hat{\mathbf{x}}(\mathbf{y}) \rmv - \rmv \xtrue \|_2^{2}  \} } \,\geq\, \sqrt{\mathsf{E}_{\xtrue} \{ \|\hat{\mathbf{x}}(\mathbf{y}) \|_2^{2}  \} } -\sqrt{\mathsf{E}_{\xtrue} \{ \| \xtrue \|_2^{2}  \} }   \eq \sqrt{\mathsf{E}_{\xtrue} \{ \|\hat{\mathbf{x}}(\mathbf{y}) \|_2^{2}  \} } - {\| \xtrue \|}_2 \, .
\]
From this, it follows 
\vspace{1mm}
that
\[
\sqrt{\mathsf{E}_{\xtrue} \{ \|\hat{\mathbf{x}}(\mathbf{y}) \|_2^{2}  \} } \,\leq\, \sqrt{\mathsf{E}_{\xtrue} \{ \|\hat{\mathbf{x}}(\mathbf{y}) - \xtrue \|_2^{2}  \} } \ist+ {\|Ê\xtrue \|}_2
 \,\leq\, \sqrt{C} + {\|Ê\xtrue \|}_2  \qquad \text{for all }\, \xtrue \! \in \! \mathcal{X}_{S} \, ,
\]
where \eqref{equ_bounded_MSE} has been used. Squaring both sides %%% AJ 
and using the inequality $(x + y)^2 \leq 2 (x^2 + y^2)$, we obtain
\[
\mathsf{E}_{\xtrue}\{ \|\hat{\mathbf{x}}(\mathbf{y}) \|_2^{2}  \} \,\leq\, ( \sqrt{C} + {\|Ê\xtrue \|}_2 )^2
  \,\le\, 2 \, ( C + \|Ê\xtrue \|_2^{2} )  \qquad \text{for all }\,  \xtrue \! \in \! \mathcal{X}_{S}
\vspace{-1mm}
\]
or equivalently
\be
\label{equ_proof_bounded_temp_int_1}
%% \mathsf{E}_{\xtrue}\{ \|\hat{\mathbf{x}}(\mathbf{y}) \|^{2}  \}    = 
\normconstgauss  \int_{\mathbb{R}^N} \!\| \hat{\mathbf{x}}(\mathbf{y}) \|_2^{2} \,e^{-\| \mathbf{y} - \xtrue \|_2^{2}/(2\sigma^{2}) } \ist d \mathbf{y} 
\,\leq\, 2 \ist ( C + \|Ê\xtrue \|_2^{2} ) \quad \quad \text{for all }\, \xtrue \! \in \! \mathcal{X}_{S}.
\vspace{1mm}
\ee
%For a given radius $R \in \mathbb{R}_{++}$ we define the following rectangular set of points in $\mathbb{R}^{N}$

We will now show that \eqref{equ_suff_cond_bounded} holds for $R_{0}  \rmv\rmv=\! 1$, i.e., $R \rmv\geq\! 1$.  
%% Let $\Delta$ be arbitrary except that $0 \le \Delta \le R$, and define 
We define the $N$-dimensional grid
\[
\mathcal{G} \,\triangleq\, \{ - m \Delta, -(m \!-\! 1) \Delta, \ldots,-\Delta, 0, \Delta , \ldots , m \Delta \}^{\!N}
\]
%%%%% BEGIN ALEX 08042010
where $0 < \Delta \leq R$ (hence, $R / \Delta \ge 1$) and $m = \lfloor R / \Delta \rfloor \le R/\Delta$. %%% AJ 
The number of grid points in any single dimension
%% in $\mathcal{G}$ 
satisfies 
\be
\label{equ_m_bound}
2m \rmv+\! 1 \,\leq\, \frac{2R}{\Delta} +1  %%% AJ 
\vspace{-2mm}
\ee
so that  %%% AJ 
%. The cardinality of $\mathcal{G}$ obeys:
\vspace{1mm}
\be
\label{equ_proof_bounded_temp_number_est}
| \mathcal{G} | \ist=\ist (2m \rmv+\! 1)^{N} \leq\ist \left( \frac{2R}{\Delta} +1 \right)^{\!\rmv N} \!.
\ee
We thus 
\vspace*{1mm}
have
\begin{align}
\sum_{\mathbf{x} \in \mathcal{G}} \|\mathbf{x}\|_2^{2} 
& \ist= \sum_{\mathbf{x} \in \mathcal{G}} \sum_{k=1}^{N} x_{k}^{2} 
  \,=\ist \sum_{k=1}^{N}  \sum_{\mathbf{x} \in \mathcal{G}} \rmv x_{k}^{2}  
  \,=\ist \sum_{k=1}^{N} \Bigg[ (2m \rmv+\! 1)^{N-1} \!\sum_{l= -m}^{m } \!(l \Delta)^{2} \Bigg] 
  =\, N (2m \rmv+\! 1)^{N-1} \!\sum_{l= -m}^{m } \!(l \Delta)^{2} \nonumber\\[2.5mm]Ê
& \ist\leq\, N \ist (2m \rmv+\! 1)^{N-1} \Delta^{2} \!\int_{x=-R/ \Delta}^{R / \Delta} \!x^{2}dx  
  \,\leq\, N \!\left( \frac{2R}{\Delta}+1 \right)^{\!\rmv N-1}  \frac{2}{3} \frac{R^3}{\Delta}  
\label{equ_sum_R} \\[-9mm]
& \nonumber
\end{align}
where \eqref{equ_m_bound} was used in the last step. Furthermore, for $c \ist\triangleq \normconstgauss e^{- N \Delta^{2}/(2\sigma^{2})}$, we 
\vspace{.7mm}
have
\be
\label{equ_c_larger_1}
\frac{1}{c} \, \normconstgauss \sum\limits_{\mathbf{x} \in \mathcal{G}} e^{-\| \mathbf{y} - \mathbf{x} \|_2^{2}/(2\sigma^{2}) } \geq 1 \,, \qquad
  \text{for all $\mathbf{y}$ with ${\| \mathbf{y} \|}_2 \rmv\leq\rmv R$}
\vspace{1mm}
\ee
In order to verify this inequality, consider an arbitrary $\mathbf{y} \!\in\! \mathbb{R}^{N}$ with ${\| \mathbf{y} \|}_2 \rmv\leq\rmv R$. 
Since $0 < \Delta \leq R$, and since ${\| \mathbf{y} \|}_2 \rmv\leq\rmv R$ implies that no component $y_k$ of $\mathbf{y}$ can be larger than $R$, 
there always exists a grid point $\tilde{\mathbf{x}} \!\in\! \mathcal{G}$ (dependent on $\mathbf{y}$) 
%% , denoted by $\mathbf{g}(\mathbf{y})$, 
such that $|y_k - \tilde{x}_{k} | \leq \Delta$ for all $k \in \{1,\ldots,N\}$. It follows that $\| \mathbf{y} - \tilde{\mathbf{x}} \|^{2}_{2} \leq N \Delta ^{2}$ and, in 
\vspace{-1mm}
turn, 
\[
e^{- N \Delta^{2}/(2\sigma^{2})} \ist\leq\, e^{-\| \mathbf{y} - \tilde{\mathbf{x}} \|_2^{2}/(2\sigma^{2}) } \ist\leq\ist \sum_{\mathbf{x} \in \mathcal{G}} e^{-\| \mathbf{y} - \mathbf{x} \|_2^{2}/(2\sigma^{2}) } \,, 
  \qquad {\| \mathbf{y} \|}_2 \rmv\leq\rmv R
\vspace{-1mm}
\]
which is equivalent to \eqref{equ_c_larger_1}. 

Successively using \eqref{equ_c_larger_1}, \eqref{equ_proof_bounded_temp_int_1}, \eqref{equ_proof_bounded_temp_number_est}, \eqref{equ_sum_R}, and $1 \le 2R/\Delta$,
we obtain the following sequence of inequalities:
\begin{align}
\int_{ {\| \mathbf{y} \|}_2 \leq R } \| \hat{\mathbf{x}}(\mathbf{y})\|_2^{2} \, d \mathbf{y} & \,\leq\, \int_{ {\| \mathbf{y} \|}_2 \leq R } \| \hat{\mathbf{x}}(\mathbf{y})\|_2^{2} 
\ist \bigg[ \frac{1}{c} \, \normconstgauss \sum\limits_{\mathbf{x} \in \mathcal{G}} e^{-\| \mathbf{y} - \mathbf{x} \|_2^{2}/(2\sigma^{2}) } \bigg] \, d \mathbf{y} \nonumber \\[1mm]
& \,\leq\, \frac{1}{c} \sum_{\mathbf{x} \in \mathcal{G}} \normconstgauss  \int_{\mathbb{R}^{N}} \! \| \hat{\mathbf{x}}(\mathbf{y}) \|_2^{2} \, e^{- \| \mathbf{y} - \mathbf{x} \|_2^{2}/(2\sigma^{2}) } \ist d \mathbf{y} 
\nonumber \\[1mm]
                                        & \,\leq\,  \frac{1}{c}\sum_{\mathbf{x} \in \mathcal{G}} 2 \ist (C + \|\mathbf{x}\|_2^{2}) \nonumber \\[1mm]Ê
                                 %       & \,\leq\,  \frac{2}{c}  \left(  \left( \frac{2R+2}{\Delta} \right)^N  C +  \sum_{\mathbf{x} \in \mathcal{G}} \|\mathbf{x} \|^{2} \right) \label{equ_proof_bounded_temp_central_3}\\
                                        & \,\leq\, \frac{2}{c} \bigg[  \left( \frac{2R}{\Delta} +1 \right)^{\!\rmv N} \!\rmv C + N \!\left( \frac{2R}{\Delta}+1 \right)^{\!\rmv N-1}  \frac{2}{3} \frac{R^3}{\Delta} \bigg] \nonumber \\[1.5mm]Ê
                                         &\,\leq\, \frac{2}{c} \bigg[  \left( \frac{4R}{\Delta} \right)^{\!\rmv N} \!\rmv C + N \!\left( \frac{4R}{\Delta}\right)^{\!\rmv N-1} \frac{2}{3} \frac{R^3}{\Delta} \bigg] \,.\label{equ_proof_bounded_temp_central_3} \\[-9mm]
& \nonumber
\end{align}
%The final step is to partition the $N$-dimensional ball of radius $R$, i.e., the set of $\mathbf{y}$ vectors $\{ \mathbf{y} : \| \mathbf{y} \| \leq R \}$, into a component where $\| \hat{\mathbf{x}}(\mathbf{y}) \| \leq 1$ and
%a second component where $\| \hat{\mathbf{x}}(\mathbf{y}) \| > 1$. On the latter component we can then use the inequality $\| \hat{\mathbf{x}}(\mathbf{y}) \| \leq \| \hat{\mathbf{x}}(\mathbf{y}) \|^{2}$ to obtain:
%\begin{align}
%\int_{ \| \mathbf{y} \| \leq R }  \| \hat{\mathbf{x}}(\mathbf{y}) \|^{2} d \mathbf{y} & = \int_{ \| \mathbf{y} \| \leq R \cap\| \hat{\mathbf{x}}(\mathbf{y}) \| \leq 1}  \| \hat{\mathbf{x}}(\mathbf{y})\|^{2} d \mathbf{y}  + \int_{ \| \mathbf{y} \| \leq R \cap\| \hat{\mathbf{x}}(\mathbf{y}) \| > 1}  \| \hat{\mathbf{x}}(\mathbf{y})\|^{2} d \mathbf{y} \\Ê
%& \leq \int_{ \| \mathbf{y} \| \leq R}  d \mathbf{y} + \int_{ \| \mathbf{y} \| \leq R }  \| \hat{\mathbf{x}}(\mathbf{y})\|^{2} d \mathbf{y} \\Ê
%& \leq  R^{N} \frac{\pi^{N/2}}{\Gamma(1+N/2)}+ \frac{2}{c}  \left(  \left( \frac{2R+2}{\Delta} \right)^N  C +    \left( \frac{2R+2}{\Delta} \right)^{N-1} \frac{2}{3}(R+1)^3 \right) \label{equ_equ_partition_bounded_implies_temp}
%\end{align}
%where we have used that the volume of a $N$-dimensional ball of radius $R$ is equal to $R^{N} \frac{\pi^{N/2}}{\Gamma(1+N/2)}.$\footnote{$\Gamma(x) \triangleq \int\limits_{0}^{\infty} t^{x-1} e^{-t} dt$}
It then follows from \eqref{equ_proof_bounded_temp_central_3} that for 
\vspace{1mm}
$R \geq 1$
\begin{align*}
\int_{ {\| \mathbf{y} \|}_2 \leq R }  \| \hat{\mathbf{x}}(\mathbf{y}) \|_2^{2} \, d \mathbf{y}    
&\,\leq\, \frac{2}{c} \bigg[  \left( \frac{4}{\Delta} \right)^{\!\rmv N} \!\rmv R^{N+2}  \ist\ist C + N \!\left( \frac{4}{\Delta}\right)^{\!\rmv N-1} \frac{2}{3} \frac{R^{N+2}}{\Delta} \bigg]\\[1.5mm]
&\,\leq\, \frac{2}{c} \frac{R^{N+2}}{\Delta^{N}} \bigg( \rmv 4^N C + N 4^{N} \frac{2}{3} \bigg)\\[1.5mm]
%% &\,\leq\, \frac{2}{c} \normconstgauss \bigg[  \left( \frac{4}{\Delta} \right)^{\!N} R^{N+2}  \! C + N 4^{N} \left( \frac{1}{\Delta}\right)^{\!N-1} \frac{2}{3} \frac{R^{N+2}}{\Delta} \bigg]\\[.5mm]
&\,=\, \frac{2^{2N+1}}{c  \,\Delta^{N}} \, \bigg(   \rmv C +  \frac{2N}{3}  \bigg) \, R^{N+2} \,.  \\[-9mm]
& \nonumber
\end{align*}
Thus, we have established that under the conditions of Lemma \ref{thm_bounded_tempered} (bounded MSE), the bound \eqref{equ_suff_cond_bounded} holds with $R_{0} \rmv=\! 1$, 
$B= \frac{2^{2N+1}}{c  \,\Delta^{N}} \, ( C + 2N/3 )$, and $n=N \rmv+\rmv 2$.
Therefore, it follows from Proposition \ref{prop_tempered} that an estimator with bounded MSE is 
%END ALEX 08042010
%% necessarily 
equivalent to a tempered distribution. This concludes the proof of Lemma \ref{thm_bounded_tempered}. \hfill $\Box$
\vspace{4mm}
\end{bproof}

We now continue our proof of Theorem \ref{thm_exit_uniqu_unbiasedness_S_N}. 
%% We know from Lemma \ref{thm_bounded_tempered} that,
%% unless the MSE tends to infinity, any estimator for the SSNM is equivalent to a tempered distribution 
%% and thus has a Fourier transform in the distributional sense. Therefore, in the sequel, 
%% we will restrict our attention to estimators %satisfying the conditions of Lemma \ref{thm_bounded_tempered}. %%% AJ
%% having bounded MSE. %%% AJ
%Using the Fourier transform calculus for tempered distributions \cite{Strichartz}, we can deduce the following result, which is proved in Appendix \ref{app_proof_ thm_exit_uniqu_unbiasedness}. %%% AJ
Any estimator $\hat{\mathbf{x}}(\mathbf{y})$ for the SSNM \eqref{equ_ssnm} can be written as
\be
\label{equ_correction}
\hat{\mathbf{x}}(\mathbf{y}) \ist=\ist \mathbf{y} + \hat{\mathbf{x}}'(\mathbf{y})
\ee
with the correction term $\hat{\mathbf{x}}'(\mathbf{y}) \triangleq \hat{\mathbf{x}}(\mathbf{y}) - \mathbf{y}$.
Because $\mathsf{E}_{\xtrue} \{ \hat{\mathbf{x}}(\mathbf{y}) \} = \mathsf{E}_{\xtrue} \{ \mathbf{y} \} + \mathsf{E}_{\xtrue} \{ \hat{\mathbf{x}}'(\mathbf{y}) \} = \mathbf{x} + \mathsf{E}_{\xtrue} \{ \hat{\mathbf{x}}'(\mathbf{y}) \}$, 
%% the linearity of expectations, a necessary and sufficient condition for 
$\hat{\mathbf{x}}(\mathbf{y})$ is unbiased if and only if
\be
\label{equ_cond_conv_ident_zero}
\mathbf{b}(\mathbf{x}; \hat{\mathbf{x}}) \,=\, \mathsf{E}_{\xtrue} \{ \hat{\mathbf{x}}'(\mathbf{y}) \}
  \,\equiv\, \normconstgauss \int_{\mathbb{R}^N} \!\hat{\mathbf{x}}'(\mathbf{y}) \, e^{-Ê\| \mathbf{y} - \mathbf{x} \|^{2}_{2}/(2 \sigma^{2})} \ist d \mathbf{y}
  \,=\, \mathbf{0} \qquad \text{for all }\, \xtrue \!\in\! \mathcal{X}_{S} \,.
\vspace{1mm}
\ee

%We first consider the case $S=N$ and prove the first part of Theorem \ref{thm_exit_uniqu_unbiasedness}.
%% show that for $S=N$, $\hat{\mathbf{x}} (\mathbf{y}) = \mathbf{y}$ is the unique unbiased estimator.
%%%% BEGIN ALEX 08042010
Remember that we assume that $\hat{\mathbf{x}}$ has a bounded MSE, so that according to our above proof of Lemma \ref{thm_bounded_tempered}, the estimator function $\hat{\mathbf{x}}(\mathbf{y})$ satisfies condition \eqref{equ_suff_cond_bounded} with $n=N+2$, i.e.,
\be
\label{equ_suff_cond_bounded_1}
\int_{ {\| \mathbf{y} \|}_2 \,\leq\, R} {\| \hat{\mathbf{x}}(\mathbf{y}) \|}_2^{2} \, d \mathbf{y} \,\leq\,  B R^{N+2} \qquad \text{for all} \;\, R \rmv\geq\rmv\rmv 1
\vspace{.7mm}
\ee
with $B$ as given at the end of the proof of Lemma \ref{thm_bounded_tempered}. 
We will also need the following bound, in which $\mathcal{R} \triangleq [-R,R]^N$:
\be
\label{equ_aux_bound}
\int_{ {\| \mathbf{y} \|}_2 \leq R} \| \mathbf{y} \|_2^{2} \, d \mathbf{y} 
\,\le \int_{ \mathcal{R} } \| \mathbf{y} \|_2^{2} \, d \mathbf{y}  
\,=\ist  \sum_{k=1}^{N} \int_{ \mathcal{R} } y_{k}^{2} \, d \mathbf{y} 
\,=\ist \sum_{k=1}^{N} (2R)^{N-1} \, \frac{2}{3} \ist\ist R^{3}  
\eq \frac{N}{3} \ist 2^{N} \rmv R^{N+2} \,.
\vspace{1mm}
\ee
We then have for the correction term $\hat{\mathbf{x}}'(\mathbf{y})$, for all $R \rmv\ge\rmv 1$, 
\begin{align*}
\int_{ {\| \mathbf{y} \|}_2 \leq R} \| \hat{\mathbf{x}}'(\mathbf{y}) \|_2^{2} \, d \mathbf{y} 
& \,= \int_{ {\| \mathbf{y} \|}_2 \leq R} \| \hat{\mathbf{x}}(\mathbf{y}) - \mathbf{y} \|_2^{2} \, d \mathbf{y} \\[1mm] 
& \,\leq\,  \int_{ {\| \mathbf{y} \|}_2 \leq R} 2 \ist \big( \| \hat{\mathbf{x}}(\mathbf{y}) \|_2^{2} + \|\mathbf{y} \|_2^{2} \big) \ist d \mathbf{y} \\[1mm] 
& \eq 2 \ist\ist \Bigg( \int_{ {\| \mathbf{y} \|}_2 \leq R} \| \hat{\mathbf{x}}(\mathbf{y}) \|_2^{2} \, d \mathbf{y} \, + \int_{ {\| \mathbf{y} \|}_2 \leq R}  \| \mathbf{y} \|_2^{2} \, d \mathbf{y} \Bigg) \\[1mm]Ê
 %%% BEGIN ALE X08042010
& \,\leq\, 2 \ist\ist \bigg( B R^{N+2} \,+\,  \frac{N}{3} \ist 2^{N} \rmv R^{N+2} \bigg) \\[1mm]Ê
 %%% BEGIN ALE X08042010
& \eq \bigg( 2B + \frac{N}{3} 2^{N+1} \bigg) \ist R^{N+2} 
%%% END ALEX 08042010
\end{align*}
%for $R \geq 1$,  %%% AJ 
where \eqref{equ_suff_cond_bounded_1} and \eqref{equ_aux_bound} have been used.
Therefore, the correction term $\hat{\mathbf{x}}'(\mathbf{y})$ also satisfies \eqref{equ_suff_cond_bounded} and thus, according to Proposition \ref{prop_tempered}, it is equivalent to a tempered distribution.
%\footnote{ALEX: I don't understand the structure of this sentence.}

The bias function $\mathbf{b}(\mathbf{x},\hat{\mathbf{x}})$ in \eqref{equ_cond_conv_ident_zero}  %%% AJ 
is the convolution of $\hat{\mathbf{x}}'(\mathbf{y})$ with the Gaussian function $(2\pi\sigma^2)^{-N/2} \, e^{- \| \mathbf{y} \|^{2}_{2}/(2 \sigma^{2})}$.  %%% AJ 
Because $S \! = \!N$, we have $\mathcal{X}_{S} \!  = \! \mathbb{R}^N$, and thus \eqref{equ_cond_conv_ident_zero} holds for all $\xtrue \!  \in \! \mathbb{R}^N$. Since $\hat{\mathbf{x}}'(\mathbf{y})$ is a tempered distribution and the Gaussian function is in the Schwartz class, it follows that the Fourier transform of the convolution product \eqref{equ_cond_conv_ident_zero} is a smooth function which can be calculated 
%% by converting the convolution to a 
as the pointwise product $\bar{\mathbf{x}}' (\bar{\mathbf{y}}) \, e^{- \| \bar{ \mathbf{y}} \|^{2}_{2}/(2 \sigma^{2})}$, where $\bar{\mathbf{x}}' (\bar{\mathbf{y}})$ denotes the Fourier transform of $\hat{\mathbf{x}}'(\mathbf{y})$ \cite{Strichartz}. Therefore, \eqref{equ_cond_conv_ident_zero} is equivalent to $\bar{\mathbf{x}}' (\bar{\mathbf{y}}) \, e^{- \| \bar{ \mathbf{y}} \|^{2}_{2}/(2 \sigma^{2})} = \mathbf{0}$ for all $\bar{\mathbf{y}} \! \in \! \mathbb{R}^N$. This can only be satisfied if $\bar{\mathbf{x}}' (\bar{\mathbf{y}}) \equiv \mathbf{0}$, which in turn implies that $\hat{\mathbf{x}}'(\mathbf{y}) \equiv \mathbf{0}$ (up to deviations of zero measure) and further, %%% AJ
by \eqref{equ_correction}, that $\hat{\mathbf{x}} (\mathbf{y}) = \mathbf{y}$. Recalling that $\mathcal{X}_{S} \!  = \! \mathbb{R}^N$, %%% AJ 
it is clear from \eqref{equ_LS_est} that $\hat{\mathbf{x}} (\mathbf{y}) = \mathbf{y}$ is the LS estimator. Thus, we have shown that $\LS_est(\mathbf{y}) = \mathbf{y}$ is the unique unbiased estimator for the SSNM with $S \! = \!N$.

%% \pagebreak %%%%%%%%%

%%%%%%%%%%%%%%%%%%%%%%%%%%%%%%%%%%%%%%%%
\section{Proof of Theorem \ref{thm_no_umvu}}
\label{app_proof_thm_no_umvu}
%%%%%%%%%%%%%%%%%%%%%%%%%%%%%%%%%%%%%%%%

\vspace{1mm}

%% \begin{proof}[Proof of Theorem \ref{thm_no_umvu}]
%%
We must show that there %%% AJ 
exists no UMVU estimator for the SSNM with $S \! < \! N$. The outline of our proof is as follows. We first demonstrate 
that the unique solution of the optimization problem \eqref{equ_opt_var_proof} at the parameter value $\xtrue \!=\! \zero$, i.e., $\argmin_{\hat{\mathbf{x}}(\cdot) \ist\in\, \mathcal{U}} V(\zero; \hat{\mathbf{x}})$, is 
%% given by 
the estimator $\hx^{(\zero)}(\y)=\y$. We then show that there exist unbiased estimators which have lower variance than $\hx^{(\zero)}$ at other points $\xtrue$. This implies that neither $\hx^{(\zero)}$ nor any other estimator uniformly minimizes the variance for all $\xtrue$ among all unbiased estimators.
%\footnote{ALEX: I have removed the notation $\hat{\mathbf{x}}^{(\xtrue)}$ for ``the solution to \eqref{equ_opt_var_proof} for a specific parameter vector $\xtrue$'' because this solution is not necessarily unique.}

The estimator $\hx^{(\zero)}(\y)=\y$ is a solution of \eqref{equ_opt_var_proof} when $\xtrue \! = \! \zero$ because the minimum variance at $\xtrue=\mathbf{0}$ of any unbiased estimator is bounded below by $N \sigma^2$ and $\hat{\mathbf{x}}^{(\mathbf{0})}(\y) = \mathbf{y}$ achieves this lower bound \cite{ZvikaCRB}.
To show that $\hat{\mathbf{x}}^{(\mathbf{0})}$ is the unique solution of \eqref{equ_opt_var_proof} for $\xtrue \! = \! \mathbf{0}$, suppose by contradiction that there exists a second
unbiased estimator $\hat{\mathbf{x}}_{a}$ different from $\hat{\mathbf{x}}^{(\mathbf{0})}$, also having variance $N \sigma^{2}$ at $\xtrue \!=\! \mathbf{0}$. Consider the estimator $\hx_{\text{new}} \triangleq (\hx^{(\zero)} + \hat{\mathbf{x}}_{a})/2$. Since $\hx^{(\zero)}$ and $\hat{\mathbf{x}}_{a}$ are unbiased, $\hx_{\text{new}}$
%% This estimator is a convex combination of two unbiased estimators, and therefore 
is unbiased as well. Thus, its variance is (see \eqref{equ_mean_power}) $V(\xtrue;\hx_{\text{new}}) = P(\xtrue;\hx_{\text{new}}) - \| \xtrue \|^{2}_{2}$. In particular, we obtain for 
\vspace{1mm}
$\xtrue \!=\! \mathbf{0}$
\begin{align*}
 V(\mathbf{0};\hx_{\text{new}}) &\eq P(\mathbf{0};\hx_{\text{new}}) \eq \mathsf{E}_{\xtrue = \zero} \bigg\{ \bigg\| \frac{1}{2} \ist (\hx^{(\zero)} \rmv+ \hat{\mathbf{x}}_{a}) \bigg\|^{2}_{2} \bigg\} \\ 
                                                      &\eq \frac{1}{4} \ist \big[ \,\mathsf{E}_{\xtrue = \zero} \big\{ \|\hx^{(\zero)} \|^{2}_{2} \big\}  + \mathsf{E}_{\xtrue = \zero} \big\{ \|\hat{\mathbf{x}}_{a}\|^{2}_{2} \big\} + 2 \, \mathsf{E}_{\xtrue = \zero} \big\{ (\hx^{(\zero)})^{T} \hat{\mathbf{x}}_{a} \big\} \ist \big] \\[1mm]
                                                   & \stackrel{(*)}{\,<\,}  \frac{1}{4} \ist \Big[ \,\mathsf{E}_{\xtrue = \zero} \big\{ \|\hx^{(\zero)} \|^{2}_{2} \big\}  + \mathsf{E}_{\xtrue = \zero} \big\{ \|\hat{\mathbf{x}}_{a}\|^{2}_{2} \big\} + 2 \sqrt{ \mathsf{E}_{\xtrue = \zero} \big\{ \|\hx^{(\zero)} \|^{2}_{2} \big\} \ist \mathsf{E}_{\xtrue = \zero} \big\{ \|\hat{\mathbf{x}}_{a}\|^{2}_{2} \big\}} \,\ist\Big] \\[1.5mm] 
                                                      &\eq \frac{1}{4} \rmv\cdot 4 \ist N \sigma^{2} \rmv\eq N \sigma^{2} 
\end{align*}
where the strict inequality $(*)$ follows from the Cauchy-Schwarz inequality applied to the inner product $\mathsf{E}_{\xtrue = \zero} \big\{ (\hx^{(\zero)})^{T} \hat{\mathbf{x}}_{a} \big\}$, combined with the fact that $\hat{\mathbf{x}}^{(\mathbf{0})}$ and $\hat{\mathbf{x}}_{a}$ are not linearly dependent  (indeed, $\hat{\mathbf{x}}_{a} \not= c \ist \hat{\mathbf{x}}^{(\mathbf{0})}$ %%% AJ 
since $\hat{\mathbf{x}}^{(\mathbf{0})}$ and $\hat{\mathbf{x}}_{a}$ were assumed to be different unbiased estimators). This inequality means that the variance of $\hx_{\text{new}}$ at $\xtrue \!=\! \zero$ is lower than $N \sigma^2$. But this is impossible, as $N \sigma^2$ is the minimum variance at $\xtrue \!=\! \mathbf{0}$ achieved by any unbiased estimator. Thus, we have shown that $\hat{\mathbf{x}}^{(\mathbf{0})}$ is the unique solution of \eqref{equ_opt_var_proof} for $\xtrue \!=\! \mathbf{0}$.

Next, still for $S \!<\! N$, we consider the specific parameter value $\mathbf{x}' \!\rmv \in \! \mathcal{X}_{S}$ whose components are given 
\vspace{-2mm}
by
\[
x'_{k} \,=\ist \begin{cases} 1 \ist , & k = 2,\ldots,S \!+\! 1 \ist,\\[-1mm]
0 \ist , &\text{else} \ist .
\end{cases}
\]
The estimator $\hat{\mathbf{x}}^{(\mathbf{0})}$ 
%% which has minimum variance for $\xtrue = \mathbf{0}$ yields 
has variance $V(\mathbf{x}';\hat{\mathbf{x}}^{(\mathbf{0})}) \! =Ê\! N \sigma^{2}$ at 
%% the parameter value 
$\mathbf{x}'$ (and at all other 
%% parameter vectors 
$\mathbf{x} \! \in \! \mathcal{X}_{S}$).
We will now construct an unbiased estimator $\hat{\mathbf{x}}_{b}(\mathbf{y})$ whose variance at $\mathbf{x}'$ is smaller than $N \sigma^{2}$. The components of this estimator are defined as
\be
\label{equ_second_est_proof_no_umvu}
\hat{x}_{b,k}(\mathbf{y}) \,\triangleqÊ\begin{cases} y_1 + A \ist y_1 \prod_{l=2}^{S+1} h( y_l) \ist, & k=1 \\[-1mm]
%%% y_1 + y_1 A \prod_{k=2}^{S+1} h( y_k)
        y_k \,, & k=2,\ldots,N \end{cases}
\vspace{-1mm}
\ee
where
\vspace{1mm}
\[ 
%% \label{equ_def_app_h}
h(y) \,\triangleq\ist \begin{cases} \mbox{sgn}(y) \ist, & |y| \in [0.4,0.6] \\[-1mm] 
%% [\frac{4}{10},\frac{6}{10}]
0 \ist, & \mbox{else} \end{cases}
\]
and $A \! \in \! \mathbb{R}$ is a parameter to be determined 
shortly.\footnote{The %%%%%%%%%%%%
interval $[0.4,0.6]$ in the definition of $h(y)$ is chosen rather arbitrarily. 
%% However this choice ensures that $\beta$ in \eqref{equ_P_1_quadr_form} is nonzero.} 
Any interval which ensures that $\beta$ in \eqref{equ_P_1_quadr_form} is nonzero can be 
used.} %%%%%%%%%%%%%
A direct calculation shows %%% AJ 
that $\hat{\mathbf{x}}_{b}(\mathbf{y})$ is unbiased for all $\mathbf{x} \!\in\! \mathcal{X}_S$. Note that $\hat{\mathbf{x}}_{b}(\mathbf{y})$ 
is identical to $\hat{\mathbf{x}}^{(\mathbf{0})}(\y)=\y$ except for the first component, $\hat{x}_{b,1}(\mathbf{y})$.

%% \pagebreak %%%%%%%%%

We recall that for unbiased estimators, minimizing the variance $V(\xtrue; \hat{\mathbf{x}})$ is equivalent to minimizing the mean power
$P(\xtrue; \hat{\mathbf{x}}) = \mathsf{E}_{\xtrue} \big\{ \| \hat{\mathbf{x}}(\mathbf{y}) \|^{2}_{2} \big\}$ (see \eqref{equ_mean_power}); furthermore, 
$P(\xtrue; \hat{\mathbf{x}}) = \sum_{k=1}^{N} P(\xtrue; \hat{x}_{k})$ with $P(\xtrue;\hat{x}_{k}) \triangleq \mathsf{E}_{\xtrue} \big\{ ( \hat{x}_k(\mathbf{y}) )^{2} \big\}$.  
%= \normconstgauss \int_{\mathbf{y}} | \left( \hat{\mathbf{x}}(\mathbf{y}) \right)_{k} |^2 e^{- \frac{1}{2 \sigma^{2}} \| \xtrue - \mathbf{y} \|^{2}} d \mathbf{y}.
For the proposed estimator $\hat{\mathbf{x}}_{b}$, $P(\mathbf{x}'; \hat{x}_{b,k}) = P\big(\mathbf{x}'; \hat{x}_{k}^{(\mathbf{0})} \big)$ 
%% for $k = 2,\ldots,N$.
except for $  k \!=\! 1$. Therefore, our goal is to choose $A$ such that $P(\mathbf{x}';\hat{x}_{b,1})$ is smaller than
%% that of the estimator $\hat{\mathbf{x}}^{(\mathbf{0})}$, which is given by 
$P\big(\mathbf{x}';\hat{x}_{1}^{(\mathbf{0})}\big) = \sigma^2 + (x'_{1})^{2} = \sigma^2$.
%% To this end, we develop $P_{1}(\mathbf{x}_{1};\hat{\mathbf{x}}'')$ as
We 
\vspace{1mm}
have
\be
P(\mathbf{x}';\hat{x}_{b,1}) \eq \mathsf{E}_{\mathbf{x}'} \rmv\Bigg\{ \rmv\Bigg( \ist y_1 + A \ist y_1 \prod_{l=2}^{S+1} h(y_l) \Bigg)^{\!\! 2} \ist \Bigg\} %\nonumber\\[1mm]Ê
%% &\eq  \mathsf{E}_{\mathbf{x}_{1}} \Bigg\{ ( | y_1 |^{2} + 2 \big| y_1 \big|^{2}  A  \prod_{k=2}^{S+1} h(y_k) + A^{2} | y_1 |^{2}  \prod_{k=2}^{S+1} h^{2}(y_k) \big) \bigg \} \\
                                                                 %   & \eq \mathsf{E}_{\mathbf{x}_{1}} \rmv\big\{  y_1^{2}  \big\}  \,+\, \mathsf{E}_{\mathbf{x}_{1}} \!\Bigg\{  2 \ist y_1^2 \ist A  \prod_{k=2}^{S+1} h(y_k) \Bigg\}
                                                                   %  \ist+\, \mathsf{E}_{\mathbf{x}_{1}} \!\Bigg\{  y_1^{2}  A^{2}  \prod_{k=2}^{S+1} h^{2}(y_k)  \Bigg\}  \nonumber \\[1mm]           
                                                                     %% \label{equ_proof_non_umvu_quad_form}.
 \eq \alpha \ist A^2 +\ist  \beta \ist A \ist +\ist  \gamma %\,, 
\label{equ_P_1_quadr_form}
\pagebreak %%%%%%%%%%%%
\ee
with
%% the coefficients are given by
\vspace{1mm}
\[
\alpha \ist=\ist \mathsf{E}_{\mathbf{x}'} \rmv\Bigg\{ y_1^{2} \prod_{l=2}^{S+1} h^{2}(y_l) \Bigg\} \,, \qquad
\beta \ist=\ist \mathsf{E}_{\mathbf{x}'} \rmv\Bigg\{  2 \ist y_1^{2} \prod_{l=2}^{S+1} h(y_l) \Bigg\} \,, \qquad
   %& = \frac{2 \sigma^{2}}{(2 \pi \sigma^{2})^{S/2}}  \left( \int_{y} h(y) e^{- \frac{1}{2 \sigma^{2}} (y-1)^{2}} d y \right)^{S}
\gamma \ist=\ist \mathsf{E}_{\mathbf{x}'} \rmv\big\{  y_1^{2}  \big\} \ist=\ist \sigma^{2} .
\vspace{1mm}
\]
Note that $\gamma = P\big(\mathbf{x}'; \hat{x}_{1}^{(\mathbf{0})} \big)$.
%%  and that $a,b > 0$.
From \eqref{equ_P_1_quadr_form}, the $A$ minimizing $P(\mathbf{x}';\hat{x}_{b,1})$ is obtained as $-\beta/(2 \alpha)$; the associated minimum $P(\mathbf{x}';\hat{x}_{b,1})$ is given by $\gamma - \beta^2/(4 \alpha^2)$. It can be shown that $\beta$ is nonzero due to the construction of $h(y)$. It follows that $\beta$ is positive, and therefore $P(\mathbf{x}';\hat{x}_{b,1})$ is smaller than $\gamma = P\big(\mathbf{x}';  \hat{x}_{1}^{(\mathbf{0})} \big)$. %%% AJ
Thus, using $A = -\beta/(2 \alpha)$ in \eqref{equ_second_est_proof_no_umvu}, we obtain an estimator $\hat{\mathbf{x}}_{b}$ which has a smaller component power $P(\mathbf{x}'; \hat{x}_{b,1})$ than  %%% AJ 
%% the estimator 
$\hat{\mathbf{x}}^{(\mathbf{0})}$. Since $P(\mathbf{x}'; \hat{x}_{b,k}) = P\big(\mathbf{x}'; \hat{x}_{k}^{(\mathbf{0})} \big)$ for $k = 2,\ldots,N$,
%%  the other components of $\hx''$ and $\x^{(\zero)}$ are equal, 
it follows that the overall mean power of $\hat{\mathbf{x}}_{b}$ at $\mathbf{x}'$ is smaller than that of 
%% the estimator 
$\hat{\mathbf{x}}^{(\mathbf{0})}$, i.e., $P(\mathbf{x}'; \hat{\mathbf{x}}_{b}) < P(\mathbf{x}'; \hat{\mathbf{x}}^{(\mathbf{0})})$. 
Since both estimators are unbiased, this moreover implies that at $\mathbf{x}'$, the variance of $\hat{\mathbf{x}}_{b}$ is smaller than that of $\hat{\mathbf{x}}^{(\mathbf{0})}$.
Thus, $\hat{\mathbf{x}}^{(\mathbf{0})}$ cannot be the LMVU estimator at $\xtrue \!=\! \mathbf{x}'$.
%% minimum variance unbiased estimator. 
On the other hand, as we have seen, $\hat{\mathbf{x}}^{(\mathbf{0})}$ is the unique LMVU
%% minimum variance unbiased 
estimator at $\xtrue \!=\! \mathbf{0}$. We conclude that there does not exist a single unbiased estimator which simultaneously minimizes the variance for all parameters $\mathbf{x} \! \in \! \mathcal{X}_S$.
%%  when $S < N$.
%% \end{proof}

%% \pagebreak %%%%%%%%%

%%%%%%%%%%%%%%%%%%%%%%%%%%%%%%%%%%%%%%%%
\section{Proof of Proposition~\ref{pr:hcrb}}
\label{ap:pr:hcrb}
%%%%%%%%%%%%%%%%%%%%%%%%%%%%%%%%%%%%%%%%

\vspace{1mm}

We begin by stating the multivariate HCRB.

\begin{proposition}[Gorman and Hero \cite{GormanHero}]
\label{pr:gorman hero}
Let $f(\y;\x)$ be a family of pdf's of $\y$ indexed by $\x \! \in \! \XS$, and let $\x + \v_1, \ldots, \x + \v_p$ be a set of points in $\XS$. 
Given an estimator $\hx$, define 
\begin{align*}
\m_\x & \,\triangleq\, \mathsf{E}_{\xtrue} \{\hx\} \notag \\
\delta_i \m_\x & \,\triangleq\, \m_{\x+\v_i} - \m_\x \notag\\
\bdel \m_\x & \,\triangleq\, (\delta_1 \m_\x \, \cdots \, \delta_p \m_\x )^T \notag\\[-12.5mm]
&
\end{align*}
and 
\vspace{-1.5mm}
\begin{align}
\delta_i f &\,\triangleq\, f(\y;\x+\v_i) - f(\y;\x) \notag\\
\bdel f &\,\triangleq\, (\delta_1 f \, \cdots \, \delta_p f )^T  \notag\\
\Q &\,\triangleq\, \mathsf{E}_{\xtrue} \Bigg\{ \frac{\bdel f}{f} \, \frac{\bdel f^T}{f} \Bigg\} \,. \label{eq:Q_def}
\end{align}
Then, the covariance matrix 
%% $\Cov(\hx)$ 
of $\hx$ satisfies
\beq 
\label{eq:HCR}
C(\x;\hx) \,\succeq\,
\bdel\m_\x^T \ist\Q^{\pinv} \ist \bdel\m_x.
\vspace{3mm}
\eeq
%where $\A \succeq \B$ indicates that the matrix $\A-\B$ is positive semidefinite and $^{\dagger}$ denotes the pseudoinverse.
\end{proposition}

We will now prove Proposition~\ref{pr:hcrb} by applying the multivariate HCRB \eqref{eq:HCR} to the case of unbiased estimation under Gaussian noise. 
For an unbiased estimator $\hx$, we have $\m_\x \!=\rmv\rmv \x$, so $\delta_i \m_\x \!=\rmv\rmv \v_i$ and further 
\be
\bdel \m_\x \ist=\ist \V \ist\ist\triangleq\, ( \mathbf{v}_{1} \cdots \mathbf{v}_{p} )
\label{eq:V}
\ee
%with $\V = ( \mathbf{v}_{1} \cdots \mathbf{v}_{p} )$ 
(see \eqref{equ def V}). We next show that the matrix $\Q$ in \eqref{eq:Q_def} coincides with $\J$ in \eqref{equ def J}. Because of the Gaussian noise,  
$f(\y;\x) = (2\pi\sigma^2)^{-N/2} \expp{ - \|\y \rmv-\rmv \x\|_2^2/(2\sigma^2) }$, and thus we obtain by direct 
\vspace{2mm}
calculation
\[
   \frac{\delta_i f}{f}
%&= \frac{\expp{-\frac{\|\y-\x-\v_i\|^2}{2\sigma^2}}
%             - \expp{-\frac{\|\y-\x\|^2}{2\sigma^2}}}
%        {\expp{-\frac{\|\y-\x\|^2}{2\sigma^2}}} \notag\\
\eq \expp{ \frac{2\v_i^T(\y \!-\! \x) - \|\v_i\|_2^2}{2\sigma^2} } - 1 
\vspace*{-2mm}
\]
and 
\vspace{1.5mm}
consequently
\begin{align*}
   {( \Q )}_{i,j}
&\eq \mathsf{E}_{\xtrue} \Bigg\{  \frac{\delta_i f}{f} \frac{\delta_j f}{f} \Bigg\} \notag\\[1mm]
&\eq 1 -\, \expp{ \rmv-\frac{\|\v_i\|_2^2}{2\sigma^2}} \mathsf{E}_{\xtrue} \bigg\{ \!\expp{\frac{\v_i^T(\y\!-\!\x)}{\sigma^2}} \!\rmv\bigg\} 
 \ist-\, \expp{-\frac{\|\v_j\|_2^2}{2\sigma^2}} \mathsf{E}_{\xtrue} \bigg\{ \!\exp \! \bigg( \frac{\v_j^T(\y\!-\!\x)}{\sigma^2} \bigg) \rmv \bigg\} \notag\\[1.5mm]
&\hspace*{57mm}+\, \expp{ \rmv-\frac{\|\v_i\|_2^2+\|\v_j\|_2^2}{2\sigma^2}} \mathsf{E}_{\xtrue} \bigg\{ \!\expp{\frac{(\v_i \rmv+\rmv \v_j)^T(\y \!-\! \x)}{\sigma^2}} \!\rmv \bigg\} \,. \notag \\[-9mm]
& \nonumber
\end{align*}
Now $\mathsf{E}_{\xtrue} \big\{  \!\expp{\a^T(\y\!-\! \x)} \!\big\}$ is the moment-generating function of the zero-mean Gaussian random vector $\y\!-\! \x$, which equals $\expp{\|\a\|_2^2 \, \sigma^2/2}$. 
We thus 
\vspace{1mm}
have  %%% AJ 
\begin{align}
   {(\Q )}_{i,j}
%&\eq -1 \ist+\, \expp{ \rmv-\frac{\|\v_i\|_2^2 + \|\v_j\|_2^2}{2\sigma^2}} \expp{\frac{\|\v_i+\v_j\|_2^2}{2\sigma^2}} \notag\\
& \eq 1 -\, \expp{ \rmv-\frac{\|\v_i\|_2^2}{2\sigma^2}}\expp{ \rmv\frac{\|\v_i\|_2^2}{2\sigma^2}} -\, \expp{ \rmv-\frac{\|\v_j\|_2^2}{2\sigma^2}}\expp{ \rmv\frac{\|\v_j\|_2^2}{2\sigma^2}}\notag\\[1mm]
& \hspace*{45mm} +\, \expp{ \rmv-\frac{\|\v_i\|_2^2 + \|\v_j\|_2^2}{2\sigma^2}} \expp{\frac{\|\v_i+\v_j\|_2^2}{2\sigma^2}} \notag\\[-1mm]
&\eq -1 \ist+\, \expp{\frac{\v_i^T\v_j}{\sigma^2}}
\label{eq:Q} \\[-9.5mm]
& \nonumber
\end{align}
which equals ${(\mathbf{J})}_{i,j}$ in \eqref{equ def J}. Inserting \eqref{eq:V} and \eqref{eq:Q} into \eqref{eq:HCR}, we obtain \eqref{equ hcr vector}. Finally, taking the trace of both sides of \eqref{equ hcr vector} yields \eqref{equ hcr}.
%% \end{proof}

%%%%%%%%%%%%%%%%%%%%%%%%%%%%%%%%%%%%%%%%
\section{Obtaining the CRB from the HCRB}
\label{ap:t->0}
%%%%%%%%%%%%%%%%%%%%%%%%%%%%%%%%%%%%%%%%

\vspace{1mm}

We will demonstrate that the CRB \eqref{equ crb} can be obtained as a limit of HCRBs \eqref{equ hcr} by choosing the test points $\mathbf{v}_{i}$ according to \eqref{equ test pts crb} and letting 
$t \rmv\rmv\to\rmv\rmv 0$. Since the test points \eqref{equ test pts crb} are orthogonal vectors, it follows from \eqref{equ def J} that the matrix $\J$ is diagonal. More specifically, we have
\[
\J \ist=
\begin{cases}
\big[ \exp (t^2/\sigma^2) - 1 \big] \ist \I_S \,, & {\|\x\|}_0  \!=\! S \\
 \big[ \exp (t^2/\sigma^2) - 1 \big] \ist \I_N \,, & {\|\x\|}_0 \!<\! S \,.
\end{cases}
\]
Thus, both for ${\|\x\|}_0  \!=\! S$ and for $ {\|\x\|}_0 \!<\! S$, the pseudoinverse of $\J$ is obtained simply by inverting the diagonal entries of $\J$. 
%% Applying Proposition~\ref{pr:hcrb}, we have
From \eqref{equ hcr}, we then obtain
\beq \label{equ t->0 1}
\MSE \,\ge\ist
\begin{cases}
\displaystyle \frac{S t^2}{\exp(t^2/\sigma^2) - 1} \,, & {\|\x\|}_0 \!=\! S  \\[3mm]
\displaystyle \frac{N t^2}{\exp(t^2/\sigma^2) - 1} \,, & {\|\x\|}_0 \!<\! S\,.
\end{cases}
\eeq
We now use the third-order Taylor series expansion 
\beq \label{equ t->0 taylor}
\expp{\frac{t^2}{\sigma^2}} =\, 1 + \frac{t^2}{\sigma^2} + \frac{\tau^4}{2\sigma^4} \,, \quad\; \mbox{where } \tau \!\in\! [0,t] \,.
\eeq
Substituting \eqref{equ t->0 taylor} into \eqref{equ t->0 1} 
\vspace{-1mm}
yields
\beq \label{equ t->0 2}
\MSE \,\ge\ist
\begin{cases}
\displaystyle \frac{S t^2}{t^2/\sigma^2 + \tau^4/(2 \sigma^4)} \,, & {\|\x\|}_0 \!=\! S \\[3mm]
\displaystyle \frac{N t^2}{t^2/\sigma^2 + \tau^4/(2 \sigma^4)} \,, & {\|\x\|}_0 \!<\! S \,.
\end{cases}
\eeq
In the limit as $t \rmv\rmv\rightarrow\rmv\rmv 0$, $\tau^4 \!\in\rmv\rmv [0,t^4]$ decays faster than $t^2$, and thus the bound \eqref{equ t->0 2} converges to the CRB \eqref{equ crb}. 

%% Furthermore, from 
The CRB can also be obtained by formally replacing $\expp{t^2/\sigma^2}$ 
with $1 + t^2/\sigma^2$ in \eqref{equ t->0 1}. From \eqref{equ t->0 taylor}, we have $\expp{t^2/\sigma^2} \geq 1 + t^2/\sigma^2$ for all $t \rmv\rmv>\rmv\rmv 0$. This shows that for any 
%% finite value of 
$t \rmv\rmv>\rmv\rmv 0$, the bound \eqref{equ t->0 1} is lower than the CRB \eqref{equ crb}. Thus, the CRB (which, as shown above, is obtained using the test points \eqref{equ test pts crb} in the limit $t \rmv\rmv\rightarrow\rmv\rmv 0$) is tighter than any bound that is obtained using the test points \eqref{equ test pts crb} for any fixed $t \rmv\rmv>\rmv\rmv 0$.

%% \pagebreak %%%%%%%%%

%%%%%%%%%%%%%%%%%%%%%%%%%%%%%%%%%%%%%%%%
\section{Proof of Theorem~\ref{th:hcrb}}
\label{ap:th:hcrb}
%%%%%%%%%%%%%%%%%%%%%%%%%%%%%%%%%%%%%%%%

\vspace{1mm}

We will prove the HCRB-type bound
%%  for the SSNM 
in \eqref{equ_HCR}. 
For ${\|\xtrue\|}_0 \!<\! S$, \eqref{equ_HCR} was already demonstrated by the %%% AJ 
CRB \eqref{equ crb}, and thus it remains to show
%% demonstrate the validity of the bound 
\eqref{equ_HCR} for ${\|\xtrue\|}_0 \!=\! S$. This will be done by plugging the test points \eqref{equ_test_points} into the HCRB \eqref{equ hcr},
%%  in Proposition~\ref{pr:hcrb}, 
calculating the resulting bound for an arbitrary constant $t \!>\! 0$, and then taking the limit as $t \!\rightarrow\! 0$.
We will use the following lemma, whose proof is provided at the end of this appendix.

\begin{lemma} \label{le:alg}
Let $\P$ be an $(r+1) \times (r+1)$ matrix with the following structure:
\beq
\label{eq:P_def}
\P \,=
\begin{pmatrix}
 a      & b \one^T      \\
 b\one      & \bM
%%% BEGIN ALEX 08042010
\end{pmatrix}
= \begin{pmatrix}
 a      & b      & b & b & \cdots & b      \cr
 b      & d      & c  & c       &  \cdots        & c      \cr
 b & c      & d    &  c        & \cdots & c \cr
 b & c      &  c   &  \ddots        &  \ddots & \vdots \cr
 \vdots & \vdots &  \vdots & \ddots    &  \ddots        & c      \cr
 b      & c      & c & \cdots    &  c             & d
\end{pmatrix}
\eeq
%%% END ALEX 08042010
where $\one$
%%  \in \RR^r$ is a 
is the column vector of dimension $r$ whose entries all equal $1$, and 
%% $\bM$ is given by
\beq \label{eq:bM}
\bM \eq (d \!-\! c) \ist \I_r + c \ist \one \one^T.
\vspace{-3mm}
\eeq
Let
\vspace{-1.5mm}
\beq \label{eq:def q}
q \,\triangleq\, rb^2 - ad - (r\!-\!1) \ist ac
\vspace{-2.5mm}
\eeq
and assume that
\be\label{eq:inv cond}
d \!-\! c \,\ne\, 0\,, \qquad
d + (r\!-\!1) \ist c \,\ne\, 0\,, \qquad
q \,\ne\, 0 \,.
\ee
%% \begin{subequations} \label{eq:inv cond}
%% \begin{align}
%% d - c &\ne 0, \label{eq:inv cond 1} \\
%% d + (r-1)c &\ne 0,\label{eq:inv cond 2} \\
%% q &\ne 0.\label{eq:inv cond 3}
%% \end{align}
%% \end{subequations}
Then, $\P$ is nonsingular
%% invertible 
and its inverse is given by
\beq \label{eq:P inv}
\P^{-1} =
\begin{pmatrix}
 a'      & b' \one^T      \\
 b'\one      & \bM'
\end{pmatrix}
= 
%%% BEGIN ALEX 08042010
%\begin{pmatrix}
% a'      & b'      & \multicolumn{2}{c}{\cdots}   & b'      \cr
% b'      & d'      & c'         &  \cdots         & c'      \cr
% \vdots  & c'      & \ddots     &  \ddots         & \vdots  \cr
% \vdots  & \vdots  & \ddots     &  \ddots         & c'      \cr
% b'      & c'      & \cdots     &  c'             & d'
%\end{pmatrix}
\begin{pmatrix}
 a'      & b'      & b' & b' & \cdots & b'      \cr
 b'      & d'      & c'  & c'       &  \cdots        & c'      \cr
 b'      & c'      & d'    &  c'        & \cdots & c' \cr
 b'     & c'      &  c'   &  \ddots        &  \ddots & \vdots \cr
 \vdots & \vdots &  \vdots & \ddots    &  \ddots        & c'      \cr
 b'      & c'      & c' & \cdots    &  c'             & d'
\end{pmatrix}
%%% END ALEX 08042010
\vspace{-2mm}
\eeq
where $\bM' = (d' \!\rmv\rmv -\rmv\rmv c') \ist \I_r + c' \one \one^T$ 
\vspace{1mm}
and
\be\label{eq:abcd tag}
a' \rmv\rmv\eq -\frac{ d + (r\!-\!1) \ist c}{q}\,, \qquad\!
b' \rmv\rmv \eq \frac{b}{q}\,, \qquad\!
c' \rmv\rmv\eq \frac{ac \rmv -\rmv b^2}{(d \!-\! c) \ist q}\,, \qquad\!
d' \rmv\rmv\eq \frac{(r\!-\! 1) \ist b^2 - (r\!-\! 2) \ist ac - ad}{(d \!-\! c) \ist q} \,. \rule{2mm}{0mm}
\vspace{4mm}
\ee
\end{lemma}

%% \pagebreak %%%%%%%%%

Let 
%% $\xtrue$ be a maximal-support parameter vector (
${\|\xtrue\|}_0 \rmv=\rmv S$, and assume for concreteness and without loss of generality that $\supp(\xtrue) = \{ 1, \ldots, S \}$ and that $\xitrue$, the smallest (in magnitude) nonzero component of $\xtrue$, is the $S$th entry. A direct calculation of the matrix $\J$ in \eqref{equ def J} based on the test points \eqref{equ_test_points} then yields
\[
%% \label{eq:ap:J}
\J \ist=
\begin{pmatrix}
a \I_{S-1}       & \zero_{(S-1)\times(r+1)} \cr
\zero_{(r+1)\times(S-1)} & \P
\end{pmatrix} .
%% \begin{pmatrix}
%% a &        & 0 &       &        &            &                &        \cr
%%   & \ddots &   & \multicolumn{5}{c}{0}                                 \cr
%% 0 &        & a &        &        &           &                &        \cr
%%   &        &   & a      & b      & \multicolumn{2}{c}{\cdots} & b      \cr
%%   &        &   & b      & d      & c         &  \cdots        & c      \cr
%%   &    0   &   & \vdots & c      & \ddots    &  \ddots        & \vdots \cr
%%   &        &   & \vdots & \vdots & \ddots    &  \ddots        & c      \cr
%%   &        &   & b      & c      & \cdots    &  c             & d
%% \end{pmatrix}
\]
Here, $\P$ is an $(r+1) \times (r+1)$ matrix, where $r = N \!-\rmv S$, having the structure \eqref{eq:P_def} with entries
\be 
\label{eq:abcd}
a \eq e^{t^2/\sigma^2} \!\!-\rmv 1 \,, \qquad
b \eq e^{-t \ist \xitrue/\sigma^2} \!\!-\rmv 1 \,, \qquad
c \eq e^{\xitrue^2/\sigma^2} \!\!-\rmv 1  \,, \qquad
d \eq e^{(t^2+\xitrue^2)/\sigma^2} \!\!-\rmv 1 \,.
\ee

%% \pagebreak %%%%%%%%%

We now 
%% wish to 
apply Lemma~\ref{le:alg} in order to show
%% demonstrate 
that $\J$ is nonsingular
%% invertible 
and to calculate its inverse. More precisely, it suffices to calculate the inverse for all but a finite number of values of $t$,
since any finite set of values
%% such values 
can simply be excluded from consideration when $t$ tends %%% AJ 
to $0$. When applying Lemma~\ref{le:alg},
%% To this end, 
we first have to verify that the conditions \eqref{eq:inv cond} hold for all but a finite number of values of $t$.
By substituting \eqref{eq:abcd}, it is seen that the left-hand sides of \eqref{eq:inv cond} are nonconstant entire functions of $t$, and thus have a finite number of roots on any compact set of values of $t$.
%% values of $t$. 
By Lemma~\ref{le:alg}, this implies that $\J$ is nonsingular
%% invertible 
for all but a finite number of values of values of $t$,
and that the inverse (if it exists) is given 
\vspace{.5mm}
by
\beq \label{eq:J inv}
\J^{-1} =
\begin{pmatrix}
\frac{1}{a} \I_{S-1}       & \zero_{(S-1)\times(r+1)} \\[.5mm]
\zero_{(r+1)\times(S-1)} & \P^{-1}
\end{pmatrix}
\vspace{.5mm}
\eeq
where $\P^{-1}$ is given by \eqref{eq:P inv} and \eqref{eq:abcd tag}, again with $r = N \!-\rmv S$.
Next, we observe that for our choice of test points 
\vspace{-.5mm}
\eqref{equ_test_points},
\beq \label{eq:V^T V}
\V^T\V \ist=
\begin{pmatrix}
t^2 \I_{S-1}       & \zero_{(S-1)\times(r+1)} \\[.5mm]
\zero_{(r+1)\times(S-1)} & \widetilde\P
\end{pmatrix} 
%% \begin{pmatrix}
%% t^2 &        & 0   &       &             &            &                &             \cr
%%     & \ddots &     & \multicolumn{5}{c}{0}                                           \cr
%% 0   &        & t^2 &        &             &           &                &             \cr
%%     &        &     & t^2    & -t\xitrue     & \multicolumn{2}{c}{\cdots} & -t\xitrue     \cr
%%     &        &     & -t\xitrue& t^2+\xitrue^2 & \xitrue^2   &  \cdots        & \xitrue^2     \cr
%%     &    0   &     & \vdots & \xitrue^2     & \ddots    &  \ddots        & \vdots      \cr
%%     &        &     & \vdots & \vdots      & \ddots    &  \ddots        & \xitrue^2     \cr
%%     &        &     & -t\xitrue& \xitrue^2     & \cdots    &  \xitrue^2       & t^2+\xitrue^2
%% \end{pmatrix}.
\vspace{.5mm}
\eeq
where $\widetilde\P$ is an $(r+1) \times (r+1)$ matrix having the structure \eqref{eq:P_def} with
entries
\[ 
%% \label{eq:abcd_V}
\tilde{a} \eq t^2 , \qquad
\tilde{b} \eq -t \ist \xitrue \,, \qquad
\tilde{c} \eq \xitrue^2 , \qquad
\tilde{d} \eq t^2 \rmv+\xitrue^2 .
\]
Using \eqref{equ hcr}
%% Proposition~\ref{pr:hcrb} 
together with \eqref{eq:J inv} and \eqref{eq:V^T V}, a direct calculation 
\vspace{.5mm}
yields
%% we obtain
\begin{align} 
\eps(\xtrue;\hx)
& \,\geq\, {\rm tr}\big(\V \J^\dagger \V^T\big) 
  \eq \tr{\V^T\V\J^{-1}} 
  \eq \sum_{i=1}^{N} \sum_{j=1}^N {(\V^T\V)}_{i,j} {(\J^{-1})}_{i,j} \notag\\
&\hspace*{20mm}\eq (S \!-\! 1) \ist \frac{t^2}{a} +\ist t^2 a' - 2 \ist rt\xitrue b' \rmv+\ist r \ist (r\!-\!1) \ist \xitrue^2 c' +\ist r \ist (t^2 \!+\xitrue^2) \ist d'.
\label{eq:hcr wo lim}
\end{align}
%% For (almost) any value of $t$, the above expression is a lower bound on the MSE of unbiased estimators. 
%% However, this expression is not very insightful. A much simpler lower bound can be obtained by taking 
We now take the limit $t \rmv\rmv\rightarrow\rmv\rmv 0$ in \eqref{eq:hcr wo lim}.
For the first term, we obtain
\be \label{eq:term1}
(S \!-\! 1) \ist \frac{t^2}{a}
\eq (S \!-\! 1) \ist \frac{t^2}{e^{t^2/\sigma^2} \!-\rmv 1} 
\eq (S \!-\! 1) \ist \frac{t^2}{t^2/\sigma^2 + o(t^2)}
 \;\longrightarrow\; (S \!-\! 1) \ist \sigma^2
\ee
where we have expanded $e^{t^2/\sigma^2}\!$ into a second-oder Taylor series. Here, $o(f(t))$ indicates terms which are negligible compared with $f(t)$ when $t \!\rightarrow\! 0$, i.e., 
$\lim_{t\to 0} \ist o(f(t))/f(t) = 0$.
To find the limit of the second term in \eqref{eq:hcr wo lim}, $t^2 a' \rmv= -(t^2/q) \ist [ d + (r \!-\!1) \ist c \ist ]$, we first consider the reciprocal of the first factor, $t^2/q$. We have
%% must first perform some preliminary calculations. Specifically, note that
\vspace*{-3mm}
\[
%% \label{eq:pre q/t^2}
\frac{q}{t^2}
\eq \frac{1}{t^2} \big[ r \ist \big( e^{-t \xitrue/\sigma^2} \!\!-\rmv 1 \big)^2 \rmv- \big(e^{t^2/\sigma^2} \!\!-\rmv 1\big) \big(e^{(t^2+\xitrue^2)/\sigma^2} \!\!-\rmv 1 \big) 
  - (r \!-\! 1)\big(e^{t^2/\sigma^2} \!\!-\rmv 1\big) \big(e^{\xitrue^2/\sigma^2} \!\!-\rmv 1 \big) \big] \,.
\]
%% We now wish to find the limit of $q/t^2$ as $t \rightarrow 0$. To this end, expand the exponents in \eqref{eq:pre q/t^2} to series, drop high-order terms, 
Expanding some of the $t$-dependent exponentials 
%% in \eqref{eq:pre q/t^2} to 
into Taylor series, dropping higher-order terms, and simplifying, we obtain
% Details of this boring calculation can be found in notes_091013, p.1, Aside.
\begin{align} \label{eq:q/t^2}
\frac{q}{t^2}
&\eq \frac{1}{t^2} \bigg[ r \bigg( \frac{-t\xitrue}{\sigma^2} + o(t) \rmv\bigg)^{\!2} \!- \bigg(\frac{t^2}{\sigma^2} + o(t^2) \rmv\bigg) \big(e^{(t^2+\xitrue^2)/\sigma^2} \!\!-\rmv 1 \big)
                         - (r\!-\! 1) \bigg(\frac{t^2}{\sigma^2} + o(t^2)\rmv\bigg)\big(e^{\xitrue^2/\sigma^2} \!\!-\rmv 1 \big)
                 \bigg] \rule{7mm}{0mm}\notag\\[1mm]
&\hspace*{20mm}\longrightarrow\; r  \, \frac{\xitrue^2}{\sigma^4} - \frac{1}{\sigma^2} \big(e^{\xitrue^2/\sigma^2} \!\!-\rmv 1 \big) - (r\!-\!1) \,  \frac{1}{\sigma^2} \big(e^{\xitrue^2/\sigma^2} \!\!-\rmv 1 \big) 
\eq \frac{r}{\sigma^4} \ist\big[ \xitrue^2 \rmv- \sigma^2 \big(e^{\xitrue^2/\sigma^2} \!\!-\rmv 1 \big) \big] \,.
\end{align}
For the second factor, we obtain
\beq \label{eq:d+(r-1)c}
d + (r \!-\!1)\ist c 
\eq e^{(t^2+\xitrue^2)/\sigma^2} \!\!-\rmv 1 \,+\, (r \!-\!1)\ist \big(e^{\xitrue^2/\sigma^2} \!\!-\rmv 1 \big)
\;\longrightarrow\; r \ist \big(e^{\xitrue^2/\sigma^2} \!\!-\rmv 1 \big) \,.
\eeq
Then, using \eqref{eq:q/t^2} and \eqref{eq:d+(r-1)c}, it is seen that the second term in \eqref{eq:hcr wo lim} converges 
\vspace{1mm}
to
% notes_091013, p.1, bottom
\be \label{eq:term2}
t^2 a'
\eq -\frac{t^2}{q} \ist [ d + (r \!-\!1) \ist c \ist ] 
\;\longrightarrow\; - \frac{ r \ist ( e^{\xitrue^2/\sigma^2} \!\!-\rmv 1 ) }
                      { \frac{r}{\sigma^4} \big[ \xitrue^2 - \sigma^2 \ist ( e^{\xitrue^2/\sigma^2} \!-\rmv 1) \big] } 
\eq \sigma^2 \bigg[ 1 + \frac{\xitrue^2}{\sigma^2 \ist (e^{\xitrue^2/\sigma^2} \!-\rmv 1) - \xitrue^2} \bigg] \,.
\vspace{1mm}
\ee
Next, we 
%% determine the limit of 
consider the third term in \eqref{eq:hcr wo lim}, $-2 r t \xitrue b'$, which can be written as $-2 r \xitrue \frac{b/t}{q/t^2}$.
We have
\[
\frac{b}{t} 
\eq \frac{1}{t} \big( e^{-t \ist \xitrue/\sigma^2} \!\!-\rmv 1 \big)
\eq \frac{1}{t} \bigg( \frac{-t \xitrue}{\sigma^2} + o(t) \bigg) 
\;\longrightarrow\; -\frac{\xitrue}{\sigma^2} \,.
\]
Combining with \eqref{eq:q/t^2}, we obtain
% notes_091013, p.2, top
\be \label{eq:term3}
-2 r t \xitrue b'
%% \eq -2 r \xitrue \frac{b/t}{q/t^2} 
\,\longrightarrow\; 2 r \xitrue \, \frac{\xitrue/\sigma^2}{\frac{r}{\sigma^4} \big[ \xitrue^2 \rmv- \sigma^2 \ist (e^{\xitrue^2/\sigma^2} \!-\rmv 1 ) \big]} 
\eq \frac{2 \ist \sigma^2\xitrue^2}{\xitrue^2 \rmv- \sigma^2 \ist (e^{\xitrue^2/\sigma^2} \!-\rmv 1)} \,.
\vspace{1mm}
\ee

The fourth and fifth terms in \eqref{eq:hcr wo lim} have to be calculated together because each of them by itself diverges. The sum of these terms is 
\begin{align}
r \ist (r\!-\!1) \ist \xitrue^2 c' +\ist r \ist (t^2 \!+\xitrue^2) \ist d'
&\eq \frac{r}{(d \!-\! c) \ist q} \,\big[ (r\!-\! 1) \ist \xitrue^2 \ist (ac \rmv -\rmv b^2) \ist+\ist (t^2 \!+\xitrue^2) \ist [(r\!-\!1)\ist b^2 \rmv- (r\!-\! 2) \ist ac - ad \ist\ist ] \ist \big] \notag\\[1mm]
&\eq \frac{r}{(d \!-\! c) \ist q} \,\big[ \!-\!\xitrue^2 a \ist (d \!-\! c) + t^2 \ist [(r\!-\!1)\ist b^2 \rmv- (r\!-\! 2) \ist ac - ad \ist\ist ] \ist \big] \notag\\[1mm]
%%% BEGIN ALEX 08042010
&\eq -\frac{r\xitrue^2 a}{q} \ist+\ist \frac{rt^2}{(d \!-\! c) \ist q} \, (q + ac \rmv-\rmv b^2) \notag\\[1mm]
&\eq \underbrace{-\frac{r\xitrue^2 a}{q}}_{\displaystyle z_1} \,+\, \underbrace{\frac{rt^2}{d \!-\! c}}_{\displaystyle z_2} \,+\, \underbrace{\frac{rt^2}{(d \!-\! c)\ist q}(ac \rmv-\rmv b^2)}_{\displaystyle z_3}.
%%% END ALEX 08042010
\label{eq:z123}
\end{align}
%% Let us separately calculate the limits of the terms $z_1$, $z_2$, and $z_3$. 
Using \eqref{eq:q/t^2}, $z_1$ in \eqref{eq:z123} becomes
% notes_091013, p.3, top
\be \label{eq:z1}
z_1
%%% BEGIN ALEX 08042010
\eq -\frac{r \xitrue^2 a/t^2}{q/t^2} 
\eq -r \xitrue^2 \frac{(e^{t^2/\sigma^2} \!\!-\rmv 1)/t^2}{q/t^2}
\;\longrightarrow\; -r \xitrue^2 \frac{1/\sigma^2}{\frac{r}{\sigma^4} \ist [ \xitrue^2 \rmv- \sigma^2 (e^{\xitrue^2/\sigma^2} \!\!-\rmv 1 ) ]}
%% &= -r \xitrue^2 \frac{ 1/\sigma^2 + o(1) }{ \frac{r}{\sigma^4} \left( \xitrue^2 - \sigma^2 \left(e^{\xitrue^2/\sigma^2} - 1 \right) \right) + o(1) }\notag\\
\eq -\frac{\sigma^2 \xitrue^2}{\xitrue^2 \rmv- \sigma^2 (e^{\xitrue^2/\sigma^2} \!\!-\rmv 1 )} \,.
%%% END ALEX 08042010
\vspace{-2mm}
\ee
Furthermore, a direct calculation yields
\be \label{eq:z2}
z_2
\eq \frac{rt^2}{e^{ (t^2+\xitrue^2)/\sigma^2} \!- e^{\xitrue^2/\sigma^2}}
\eq r \ist e^{-\xitrue^2/\sigma^2} \rmv\rmv\frac{t^2}{e^{t^2/\sigma^2}\! -\rmv\rmv 1} 
\;\longrightarrow\; r \sigma^2 e^{-\xitrue^2/\sigma^2}.
\ee
To take the limit of $z_3$, first note 
\vspace{1mm}
that
%notes_091013, p.3, middle
\begin{align*}
\frac{ac \rmv-\rmv b^2}{d \!-\! c}
&\eq \frac{(e^{t^2/\sigma^2} \!\!-\rmv 1)(e^{\xitrue^2/\sigma^2} \!\!-\rmv 1)\rmv- (e^{-t \ist \xitrue/\sigma^2} \!\!-\rmv 1)^2}{e^{(t^2+\xitrue^2)/\sigma^2} \!- e^{\xitrue^2/\sigma^2}}\\[1mm]
&\hspace*{30mm}\longrightarrow\; \frac{(t^2/\sigma^2)(e^{\xitrue^2/\sigma^2} \!\!-\rmv 1)\rmv- (-t \ist \xitrue/\sigma^2)^2}{e^{\xitrue^2/\sigma^2} \, t^2/\sigma^2}
\eq \frac{ \sigma^2 (e^{\xitrue^2/\sigma^2} \!\!-\rmv 1) - \xitrue^2 }{ \sigma^2 \ist e^{\xitrue^2/\sigma^2} } \;. \\[-12mm]
& \nonumber
\end{align*}
Together with \eqref{eq:q/t^2}, we thus 
\vspace{1mm}
have
%notes_091013, p.3, middle and bottom
\be \label{eq:z3}
z_3
\eq r \ist \frac{t^2}{q} \, \frac{ac \rmv-\rmv b^2}{d \!-\! c}
\;\longrightarrow\; r \ist \frac{1}{\frac{r}{\sigma^4} \ist [ \xitrue^2 \rmv- \sigma^2 (e^{\xitrue^2/\sigma^2} \!\rmv-\rmv 1 ) ]} \,
  \frac{ \sigma^2 (e^{\xitrue^2/\sigma^2} \!\!-\rmv 1) - \xitrue^2 }{ \sigma^2 \ist e^{\xitrue^2/\sigma^2} } 
\eq - \sigma^2 e^{-\xitrue^2/\sigma^2}.
\vspace{1mm}
\ee
Adding the limits of $z_1$, $z_2$, and $z_3$ in \eqref{eq:z1}--\eqref{eq:z3}, we find that the sum of the fourth and fifth terms in \eqref{eq:hcr wo lim} converges 
\vspace{.5mm}
to
% notes_091013 p.4 top
\vspace{-2mm}
\be \label{eq:term45}
z_1+z_2+z_3
\;\longrightarrow\; \frac{-\ist\sigma^2 \xitrue^2}{\xitrue^2 \rmv- \sigma^2 (e^{\xitrue^2/\sigma^2} \!\rmv-\rmv 1 )} \ist+\ist (r\!-\!1) \ist \sigma^2 \ist e^{-\xitrue^2/\sigma^2}.
\vspace{1mm}
\ee

Finally, adding the limits of all terms in \eqref{eq:hcr wo lim} as given by \eqref{eq:term1}, \eqref{eq:term2}, \eqref{eq:term3}, and \eqref{eq:term45} and simplifying, we obtain the following result for the limit of the bound \eqref{eq:hcr wo lim} for $t \rmv\rmv\to\rmv\rmv 0$: 
%% \begin{align}
%% \eps(\xtrue;\hx)
%% &\ge (S-1)\sigma^2
%%    + \sigma^2 - \frac{\sigma^2 \xitrue^2}{\xitrue^2 - \sigma^2\left(e^{\xitrue^2/\sigma^2}-1\right)} \notag\\
%% &\quad  + \frac{2\sigma^2 \xitrue^2}{\xitrue^2 - \sigma^2\left(e^{\xitrue^2/\sigma^2}-1\right)} + (r-1)\sigma^2 e^{-\xitrue^2/\sigma^2} \notag\\
%% &\quad  - \frac{\sigma^2 \xitrue^2}{\xitrue^2 - \sigma^2\left(e^{\xitrue^2/\sigma^2}-1\right)}
%% \end{align}
%% which simplifies to
\[
\eps(\xtrue;\hx) \,\ge\, S\sigma^2 + (r \!-\! 1) \ist \sigma^2 \ist e^{-\xitrue^{2}/\sigma^{2}} .
\vspace{-2mm}
\]
This equals \eqref{equ_HCR}, as 
\vspace{4mm}
claimed.

\emph{Proof of Lemma~\ref{le:alg}:}
We first calculate the inverse of $\bM$ in \eqref{eq:bM}. Applying the Sherman--Morrison--Woodbury formula \cite[\S 2.8]{bernstein05}
\[
\left( \A + c \ist \u \v^T \right)^{-1} =\, \A^{-1} - \frac{c}{1 + c \ist \v^T \!\A^{-1} \u} \,\A^{-1} \u\v^T \!\A^{-1}
\]
to \eqref{eq:bM} and simplifying 
\vspace{-1mm}
yields
\beq \label{eq:bM inv}
\bM^{-1} =\,\frac{1}{d \!-\! c} \ist \I_r \ist-\ist \frac{c}{(d \!-\! c) \ist [d+(r \!-\! 1) \ist c]} \, \one \one^T.
\vspace{1mm}
\eeq
Next, we invoke the block inversion 
\pagebreak %%%%%%%%
lemma \cite[\S 2.8]{bernstein05}
\[
%% \label{eq:block inv_0}
\begin{pmatrix}
\A & \B^T \cr
\B & \bM
\end{pmatrix}^{\!\!\!-1}
\!=
\begin{pmatrix}
\bE^{-1}             & -\bE^{-1} \B^T \bM^{-1}  \cr
-\bM^{-1} \B \bE^{-1} & \bM^{-1} \rmv+ \bM^{-1} \B \bE^{-1} \B^T \bM^{-1}
\end{pmatrix} , \quad\; \text{with} \;\, \bE \ist\triangleq \A \rmv- \B^T\bM^{-1}\B \,.
\vspace{1mm}
\]
Specializing to $\A \rmv= a$ and $\B = b \one$ as is appropriate for $\P$ in \eqref{eq:P_def},
we obtain for the inverse 
\vspace{1mm}
of $\P$
\beq \label{eq:block inv}
\P^{-1} \ist=
\begin{pmatrix}
1/e             &\; -(b/e) \ist \one^T \bM^{-1}  \cr
-(b/e) \ist \bM^{-1}\one &\; \bM^{-1} \rmv+ (b^2/e) \ist \bM^{-1} \one\one^T \bM^{-1}
\end{pmatrix} , 
\quad\; \text{with} \;\, e \ist\triangleq a - b^2 \one^T\bM^{-1}\one \,.
\vspace{1mm}
\eeq
We now develop the various blocks of $\P^{-1}$ by using the expression of $\bM^{-1}$ in \eqref{eq:bM inv}. We first consider the upper-left block, $1/e$. We 
\vspace{.5mm}
have
\[
e \eq a \ist-\ist \frac{b^2}{d \!-\! c} \,\one^T \!\left[ \I_r - \frac{c}{d+(r \!-\! 1) \ist c} \ist \one\one^T \right] \!\one
\eq a - \frac{b^2}{d \!-\! c} \left[ r - \frac{cr^2}{d+(r \!-\! 1) \ist c} \right] 
\eq \frac{ad + (r \!-\! 1) \ist ac -rb^2}{d + (r \!-\! 1) \ist c} \, .
\vspace{.7mm}
\]
%% \begin{align}
%% e
%% %&= a - b^2 \one^T \left( \frac{1}{d-c}\I - \frac{c}{(d-c)(d+(r-1)c)} \one \one^T \right) \one \notag\\
%% &\eq a \ist-\ist \frac{b^2}{d \!-\! c} \,\one^T \!\left[ \I_r - \frac{c}{d+(r \!-\! 1) \ist c} \ist \one\one^T \right] \!\one \notag\\[1mm]
%% &\eq a - \frac{b^2}{d \!-\! c} \left[ r - \frac{cr^2}{d+(r \!-\! 1) \ist c} \right] \notag\\[1mm]
%% &\eq \frac{ad + (r \!-\! 1) \ist ac -rb^2}{d + (r \!-\! 1) \ist c} \, .
%% \end{align}
Thus, using the definitions in \eqref{eq:def q} and \eqref{eq:abcd tag} yields
\beq \label{eq:bE inv}
\frac{1}{e} \eq -\frac{d + (r\!-\!1)c}{q} \eq a'
\eeq
which proves
%% demonstrates 
the validity of the upper-left entry of $\P^{-1}$ in \eqref{eq:P inv}. Next, using \eqref{eq:bM inv} and \eqref{eq:bE inv} and simplifying, the upper-right block
%% term $-\bE^{-1}\B^T\bM^{-1}$ 
%% $-(b/e) \ist \one^T \bM^{-1}$ 
in \eqref{eq:block inv} 
\vspace{.5mm}
becomes
\[
-\frac{b}{e} \, \one^T \bM^{-1}
\eq -ba' \! \left[ \frac{1}{d \!-\! c} - \frac{rc}{(d \!-\! c) \ist [d+(r \!-\! 1) \ist c]} \right] \!\one^T 
  \!\eq\ist - \frac{ba'}{d+(r \!-\! 1) \ist c} \,\one^T
  \!\eq\ist \frac{b}{q} \ist \one^T \!\eq\ist b' \ist \one^T .
\vspace{1mm}
\]
Thus, we have shown the validity of the first row and first column of $\P^{-1}$ in \eqref{eq:P inv}.
Finally, to develop the remaining block $\bM^{-1} \rmv+ (b^2/e) \ist \bM^{-1} \one\one^T \bM^{-1}$ in \eqref{eq:block inv}, we first 
\vspace{.5mm}
calculate
% Note: details of this boring calculation are in notes_091012.pdf, page 3
\be
\label{eq:u_def}
\u \,\triangleq\, \bM^{-1} \one
\eq\ist \frac{1}{d \!-\! c} \rmv \left[ 1 - \frac{rc}{d+(r \!-\! 1) \ist c} \right] \!\one
%&= \frac{1}{d-c} \left( \frac{d-c}{d+(r-1)c} \right) \one \notag\\
\rmv\eq\ist \frac{1}{d+(r \!-\! 1) \ist c} \,\one \,.
\vspace{-2mm}
\ee
We then 
\vspace{.5mm}
have
\be  \label{eq:bM inv dev}
\bM^{-1} \rmv+ \frac{b^2}{e} \ist \bM^{-1} \one\one^T \bM^{-1} \rmv 
\eq\ist \bM^{-1} \rmv+ b^2 a' \u \u^T 
%&\eq \bM^{-1} - \frac{b^2}{q(d + (r-1)c)} \one\one^T \notag\\
%&\quad= \frac{1}{d-c}\I - \frac{c}{(d-c)(d+(r-1)c)} \one \one^T - \frac{b^2}{q(d + (r-1)c)} \one\one^T \notag\\
\!\eq\ist \frac{1}{d \!-\! c} \,\I_r - \frac{1}{d+(r \!-\! 1) \ist c} \left[ \frac{c}{d \!-\! c} + \frac{b^2}{q} \right] \!\one\one^T
\vspace{.5mm}
\ee
where \eqref{eq:bM inv}, \eqref{eq:u_def}, and the definition of $a'$ in \eqref{eq:abcd tag} were used.
%%  in the last step.
Using the definition of $q$ in \eqref{eq:def q} and simplifying, the factor in brackets can be written as
% Note: details of this boring calculation are in notes_091012.pdf, page 3
\[
\frac{c}{d \!-\! c} + \frac{b^2}{q}
\eq \frac{cq + (d \!-\! c) \ist b^2}{(d \!-\! c) \ist q} \notag\\
%% \eq  \frac{(d-c)b^2 - cad - (r-1)ac^2 + rb^2c}{(d-c)q} \notag\\
\eq \frac{[d+(r \!-\! 1)c \ist ] \ist (b^2 \!-\rmv ac)}{(d \!-\! c) \ist q} \,.
\vspace{.5mm}
\]
Substituting back into \eqref{eq:bM inv dev}, we 
\vspace{.5mm}
obtain
\[
\bM^{-1} \rmv+ \frac{b^2}{e} \ist \bM^{-1} \one\one^T \bM^{-1} \rmv
\eq\ist \frac{1}{d \!-\! c} \ist \I_r - \frac{b^2\!-\rmv ac}{(d \!-\! c) \ist q} \one\one^T
\!\eq\ist \frac{1}{d \!-\! c} \ist \I_r + c' \one\one^T.
\]
Thus, within the $r \!\times\! r$ lower-right block of $\P^{-1}$, the off-diagonal entries all equal $c'$, as required. Furthermore, the diagonal entries in this block are given by
% Note: details of this boring calculation are in notes_091012.pdf, page 4
\[
\frac{1}{d \!-\! c} - \frac{b^2\!-\rmv ac}{(d \!-\! c) \ist q}
\eq \frac{(r \!-\! 1) \ist b^2 - ad - (r \!-\rmv 2) \ist ac}{(d \!-\! c) \ist q}
\eq d'
\]
which completes the proof of the lemma. \hfill  $\Box$
%% \vspace{3mm}
%% \end{proof}

%%%%%%%%%%%%%%%%%%%%%%%%%%%%%%%%%%%%%%%%
\section{Proof of Lemma \ref{equ_opt_k_in_supp}} 
\label{app_proof_opt_est_k_in_supp}
%%%%%%%%%%%%%%%%%%%%%%%%%%%%%%%%%%%%%%%%

\vspace{1mm}

%% \pagebreak %%%%%%%%%%

%% \begin{proof}[Proof of Theorem \ref{equ_opt_k_in_supp}]
%\footnote{ALEX: I would consider moving this to an appendix too.}
Let $\xtrue \!\in\! \mathcal{X}_{S}$ with ${\| \xtrue \|}_{0} \! =Ê\! S$ and consider a fixed $k \in \supp(\xtrue)$. We have to show that 
%% the unique LMVU of $x_k$, i.e., 
a solution of \eqref{equ_scalar_opt}, i.e.,
\vspace{-1.5mm}
\begin{equation}
\label{equ_scalar_opt_recall}
\argmin_{\hat{x}(\cdot) \ist\in\, \mathcal{U}^k} \ist \mathsf{E}_{\xtrue} \big\{ ( \hat{x} ( \mathbf{y} ) )^{2} \big\} \,, \qquad \text{with} \;\; 
\mathcal{U}^{k} \rmv= \big\{\hat{x}(\cdot) \, \big| \, \mathsf{E}_{\tilde{\mathbf{x}}} \{ \hat{x}( \mathbf{y} ) \} = \tilde x_k   \;\, \text{for all}\; \tilde{\mathbf{x}} \!\in\! \mathcal{X}_{S} \big\}
\vspace{1mm}
\end{equation}
is given by $\hat{x}^{(\xtrue)}_{k}(\mathbf{y}) = y_{k}$. Let $\varepsilon_0 \triangleq \min_{\hat{x}(\cdot) \ist\in\, \mathcal{U}^k} \ist \mathsf{E}_{\xtrue} \big\{ ( \hat{x}( \mathbf{y} ) )^{2} \big\}$ denote the mean power of the LMVU estimator defined by \eqref{equ_scalar_opt_recall}. We will show that $\varepsilon_0 \ge \sigma^{2} \rmv+ x_k^{2}$ and, furthermore, that $\sigma^{2} \rmv+ x_k^{2}$ is achieved by the estimator $\hat{x}^{(\xtrue)}_{k}(\mathbf{y}) = y_{k}$.

Let $\mathcal{C}_{\xtrue}^k\rmv$ denote the set of all $S$-sparse vectors $\tilde{\mathbf{x}}$ which equal $\xtrue$ except possibly for the $k$th component, i.e., 
$\mathcal{C}_{\xtrue}^k \triangleq \big\{ \tilde{\mathbf{x}} \!\in\! \mathcal{X}_{S} \,\big|\, \tilde{x}_l = x_l \;\, \text{for all}\;\,  l 
%% \in \supp(\xtrue) \!\setminus\! \{k\} 
\rmv\not=\rmv k \big\}$. Consider the modified optimization 
\vspace{.5mm}
problem
\be
\label{equ_opt_proof_k_in_supp_1}
\argmin_{\hat{x}(\cdot) \ist\in\, \mathcal{U}^{k}_{\xtrue}} \ist \mathsf{E}_{\xtrue} \big\{ ( \hat{x}( \mathbf{y} ) )^{2} \big\} \,, \qquad \text{with} \;\; 
\mathcal{U}^{k}_{\xtrue} \,\triangleq\, \big\{\hat{x}(\cdot) \, \big| \, \mathsf{E}_{\tilde{\mathbf{x}}} \{ \hat{x}( \mathbf{y} ) \} = \tilde x_k \;\,\text{for all}\;\, \tilde{\mathbf{x}} \!\in\rmv \mathcal{C}^k_{\xtrue} \big\}
\vspace{1.5mm}
\ee
and let $\varepsilon_0' \triangleq \min_{\hat{x}(\cdot) \ist\in\, \mathcal{U}^k_{\xtrue}} \ist \mathsf{E}_{\xtrue} \big\{ ( \hat{x}( \mathbf{y} ) )^{2} \big\}$ denote the mean power of the estimator defined by \eqref{equ_opt_proof_k_in_supp_1}. Note the distinction between $\mathcal{U}^{k}$ and $\mathcal{U}_{\xtrue}^{k}$: %%% AJ 
$\mathcal{U}^{k}$ is the set of estimators of $x_k$ which %% AJ
are unbiased for all $\tilde{\mathbf{x}} \!\in\! \mathcal{X}_{S}$ whereas $\mathcal{U}_{\xtrue}^{k}$ is the set of estimators of $x_k$ which %% AJ
are unbiased for all $\tilde{\mathbf{x}} \!\in\! \mathcal{X}_{S}$ which equal a given, fixed $\xtrue$ except possibly for the $k$th component. %%% AJ 
%% such that all components $\tilde{x}_l$ except for the $k$th are fixed (namely, equal to the respective components $x_l$ of a given $\xtrue$). 
Therefore, the unbiasedness requirement expressed by $\mathcal{U}^{k}$ is more restrictive than that expressed by $\mathcal{U}^{k}_{\xtrue}$, i.e., $\mathcal{U}^{k}\! \subseteq \mathcal{U}^{k}_{\xtrue}$, which implies that 
\be
\label{equ_epsilon0-bound}
\varepsilon_0' \rmv \leq \varepsilon_{0} \,.
\ee

We will use the following result, which is proved at the end of this appendix.

\begin{lemma}
\label{lemma_cond_expect}
Given an arbitrary estimator $\hat{x}( \mathbf{y}) \rmv\rmv\in \mathcal{U}^{k}_{\xtrue}$, the estimator
\be
\label{equ_new-est}
\hat{x}_{c} ( y_k ) \,\triangleq\, \mathsf{E}_{\xtrue} \{ \hat{x}(\mathbf{y}) |\ist  y_k \} %\frac{1}{(2 \pi \sigma^{2})^{(N-1)/2}} \int_{\tilde{\mathbf{y}}, \left( \tilde{\mathbf{y}}Ê\right)_{k} =\left(\mathbf{y}Ê\right)_{k}  } \hat{x}(\tilde{\mathbf{y}}) \prod_{l \in \{1,...,N\} \setminus \{k\}} e^{-\frac{1}{2 \sigma^{2}} | \left( \tilde{\mathbf{y}} - \xtrue \right)_{l} |^{2} } d \left( \tilde{\mathbf{y}} \right)_{l}
\ee
also satisfies the constraint $\hat{x}_{c}( y_k ) \rmv\rmv\in \mathcal{U}^{k}_{\xtrue}$, and its mean power does not exceed that obtained by $\hat{x}$, i.e., 
$\mathsf{E}_{\xtrue} \{ ( \hat{x}_{c}( y_k ) )^{2} \} \le \mathsf{E}_{\xtrue} \{ ( \hat{x}( \mathbf{y} ) )^{2} \}$.
\end{lemma}

%% \pagebreak %%%%%%%%%%

Thus, to each estimator $\hat{x}( \mathbf{y}) \rmv\rmv\in \mathcal{U}^{k}_{\xtrue}$ which %%% AJ 
depends on the entire observation $\mathbf{y}$, we can always find at least one estimator $\hat{x}_{c}( y_k) \rmv\rmv\in \mathcal{U}^{k}_{\xtrue}$ which %%% AJ 
depends only on the observation component $y_k$ and is at least as good. Therefore, with no loss in optimality, we can restrict the optimization problem \eqref{equ_opt_proof_k_in_supp_1} to estimators $\hat{x}( y_k) \rmv\rmv\in \mathcal{U}^{k}_{\xtrue}$ which depend on $\mathbf{y}$ only via its $k$th component $y_k$. %%% AJ 
This means that \eqref{equ_opt_proof_k_in_supp_1} can be replaced 
\vspace{.5mm}
by
\be
\label{equ_opt_problem_equiv_3}
\argmin_{\hat{x}(\cdot) \ist\in\, \widetilde{\mathcal{U}}^{k}} \ist \mathsf{E}_{\xtrue} \big\{ ( \hat{x}( y_k ) )^{2} \big\} \,, \qquad \text{with} \;\; 
\widetilde{\mathcal{U}}^{k} \ist\triangleq\, \big\{\hat{x}(\cdot) \, \big| \, \mathsf{E}_{\tilde{\mathbf{x}}} \{ \hat{x}( y_k ) \} = \tilde x_k \;\,\text{for all}\;\, \tilde{\mathbf{x}} \!\in\! \mathbb{R}^{N} \big\} \,.
\vspace{.7mm}
\ee
Note that in the definition of $\,\widetilde{\mathcal{U}}^{k}\rmv$, we can use the requirement $\tilde{\mathbf{x}} \!\in\! \mathbb{R}^{N}$ instead of $\tilde{\mathbf{x}} \!\in\! \mathcal{C}^k_{\xtrue}$ since the expectation $\mathsf{E}_{\tilde{\mathbf{x}}} \{ \hat{x}( y_k ) \}$ does not depend on the components $\tilde{x}_{l}$ with $l \rmv\neq\rmv k$. 
The corresponding minimum mean power $\min_{\hat{x}(\cdot) \ist\in\, \widetilde{\mathcal{U}}^{k}} \mathsf{E}_{\xtrue} \big\{ ( \hat{x}( y_k ) )^{2} \big\}$ is still equal to $\varepsilon_0'$. However, the new problem %%% AJ
\eqref{equ_opt_problem_equiv_3} is equivalent to the classical problem of finding the LMVU estimator 
%% at $x_k$ 
of a scalar $x_{k}$ based on the observation $y_k = x_k + n_k$, with $n_{k} \sim \mathcal{N} ( 0 , \sigma^{2} )$. A solution of this latter problem is the estimator $\hat{x}(y_k) = y_k$, whose variance and mean power are $\sigma^{2}$ and $\sigma^{2} \rmv+ x_k^{2}$, respectively \cite{scharf91}. Thus, a solution of \eqref{equ_opt_problem_equiv_3} or, equivalently, of \eqref{equ_opt_proof_k_in_supp_1} is the trivial estimator $\hat{x}(y_k) = y_k$, 
\vspace{-1mm}
and 
\be
\label{equ_epsilon0-result}
\varepsilon_0' \rmv \eq \sigma^{2} \rmv+ x_k^{2} \,.
\ee

Combining \eqref{equ_epsilon0-bound} and \eqref{equ_epsilon0-result}, we see that the minimum mean power for our original optimization problem \eqref{equ_scalar_opt_recall}
\vspace*{-3mm}
satisfies  %%% AJ %is lower-bounded as
\[
%% \label{equ_epsilon-bound}
\varepsilon_0 \,\geq\, \sigma^{2} \rmv+ x_k^{2} \,. 
\]
As we have shown, this lower bound is achieved by the estimator $\hat{x}(y_k) = y_k$. In addition, $\hat{x} (y_k) = y_k$ is an element of $\mathcal{U}^{k}$, the constraint set of \eqref{equ_scalar_opt_recall}. 
Therefore, it is a solution of 
\vspace{3mm}
\eqref{equ_scalar_opt_recall}. 
%% \end{proof}

%% \pagebreak %%%%%%%%%%
\emph{Proof of Lemma \ref{lemma_cond_expect}}:
Consider a fixed $\mathbf{x} \!\in\! \mathcal{X}_{S}$ and an estimator $\hat{x}( \mathbf{y}) \rmv\rmv\in \mathcal{U}^{k}_{\xtrue}$. In order to show the first statement of the lemma, $\hat{x}_{c}( y_k ) \rmv\rmv\in \mathcal{U}^{k}_{\xtrue}$, we first note that
%% \vspace{-3mm}
\be
\label{equ_tilde-or-not}
\mathsf{E}_{\xtrue} \{ \hat{x}(\mathbf{y}) |\ist y_k \} \eq \mathsf{E}_{\tilde{\mathbf{x}}} \{ \hat{x}(\mathbf{y}) |\ist y_k \} \,, \qquad \text{for any} \;\, \tilde{\mathbf{x}} \!\in\rmv \mathcal{C}_{\xtrue}^k \,.
\vspace{-.5mm}
\ee
We now have for $\tilde{\mathbf{x}} \!\in\rmv \mathcal{C}_{\xtrue}^k$
\[
\mathsf{E}_{\tilde{\mathbf{x}}} \{ \hat{x}_{c} ( y_k ) \} 
\ist\stackrel{(a)}{\eq}\ist \mathsf{E}_{\tilde{\mathbf{x}}} \{  \mathsf{E}_{\xtrue} \{ \hat{x}(\mathbf{y}) |\ist y_k \} \} 
\ist\stackrel{(b)}{\eq}\ist \mathsf{E}_{\tilde{\mathbf{x}}} \{  \mathsf{E}_{\tilde{\mathbf{x}}} \{ \hat{x}(\mathbf{y}) |\ist y_k \} \} 
\ist\stackrel{(c)}{\eq}\ist \mathsf{E}_{\tilde{\mathbf{x}}} \{  \hat{x}(\mathbf{y}) \} 
\ist\stackrel{(d)}{\eq}\ist \tilde{x}_{k}
 % \mathsf{E}_{\mathbf{x}} \{ \hat{x}' ( y_k ) \} & = \frac{1}{(2 \pi \sigma^{2})^{N/2}} \int_{ y_k } d y_ke^{-\frac{1}{2 \sigma^{2}} (  \left( \mathbf{y} - \mathbf{x} \right)_{k})^{2}}   \int_{\tilde{\mathbf{y}},   \left( \tilde{\mathbf{y}}Ê\right)_{k} =\left(\mathbf{y}Ê\right)_{k} } \hat{x}(\mathbf{y}) \prod_{l \in \{1,...,N\} \setminus \{k\}}  e^{-\frac{1}{2 \sigma^{2}} | y_l - x_l |^{2} } d \left( \tilde{\mathbf{y}} \right)_{l} \\Ê
% & = \frac{1}{(2 \pi \sigma^{2})^{N/2}}  \int_{\mathbf{y}} \hat{x}(\mathbf{y}) e^{- \frac{1}{2 \sigma^{2}} \| \mathbf{y} - \mathbf{x} \|^{2}} d \mathbf{y} =  \mathsf{E}_{\mathbf{x}} \{ \hat{x}(\mathbf{y}) \} = \left( \mathbf{x} \right)_{k}
\vspace{.5mm}
\]
where we used the definition \eqref{equ_new-est} in $(a)$,
the identity \eqref{equ_tilde-or-not} in $(b)$, 
the law of total probability \cite{papoulis} in $(c)$, 
and our assumption $\hat{x}( \mathbf{y}) \rmv\rmv\in \mathcal{U}^{k}_{\xtrue}$ in $(d)$. Thus, $\hat{x}_{c}( y_k ) \rmv\rmv\in \mathcal{U}^{k}_{\xtrue}$.

Next, the inequality $\mathsf{E}_{\xtrue} \{ ( \hat{x}_{c} ( y_k ) )^{2} \} \le \mathsf{E}_{\xtrue} \{ ( \hat{x}( \mathbf{y} ) )^{2} \}$ is proved as 
\vspace{1mm}
follows:
\[
\mathsf{E}_{\xtrue} \{ ( \hat{x}( \mathbf{y} ) )^{2} \}  
\ist\stackrel{(a)}{\eq}\ist \mathsf{E}_{\xtrue} \{ \mathsf{E}_{\xtrue} \{ ( \hat{x}( \mathbf{y} ) )^{2} |\ist y_k \} \} 
\ist\stackrel{(b)}{\,\geq\,}\ist \mathsf{E}_{\xtrue} \{ ( \mathsf{E}_{\xtrue} \{ \hat{x} ( \mathbf{y} ) |\ist y_k \}  )^{2} \} 
\ist\stackrel{(c)}{\eq}\ist \mathsf{E}_{\xtrue} \{ ( \hat{x}_{c} ( y_k ) )^{2} \}
%\mathsf{E}_{\xtrue} \{ | \hat{x} ( \mathbf{y} ) |^{2} \}  & = %\normconstgauss \int_{\mathbf{y}} |\hat{x}(\mathbf{y})|^{2} e^{- \frac{1}{2 \sigma^{2}} \| \mathbf{y} - \xtrue \|} d \mathbf{y} \\Ê
% \normconstgauss \int_{\mathbf{y}_k} d y_k e^{-\frac{1}{2 \sigma^{2}} | \left( \mathbf{y} - \xtrue\right)_{k}|^{2}} \int_{\tilde{\mathbf{y}}, \left( \tilde{\mathbf{y}}Ê\right)_{k} =\left(\mathbf{y}Ê\right)_{k}  }  |\hat{x}(\tilde{\mathbf{y}})|^{2} \prod_{l \in \{1,...,N\} \setminus \{k\}}
% e^{-\frac{1}{2 \sigma^{2}} | \left( \tilde{\mathbf{y}} - \xtrue \right)_{l}|^{2} } d  y_l \\Ê
%& \geq \normconstgauss \int_{\mathbf{y}_k} d \mathbf{y}_{k} e^{-\frac{1}{2 \sigma^{2}} | \left( \mathbf{y} - \xtrue \right)_{k}|^{2}} \bigg| \int_{\mathbf{y}_{l} \quad l \neq k}   \hat{x}(\mathbf{y}) \prod_{l}  e^{-\frac{1}{2 \sigma^{2}} |\mathbf{y}_l - \mathbf{x}_{0,l}|^{2} } d \left(\mathbf{y}\right)_{l}  \bigg |^{2} \\Ê
%& = \normconstgauss \int_{\mathbf{y}_k} d \mathbf{y}_{k} e^{-\frac{1}{2 \sigma^{2}} | \left( \mathbf{y} - \xtrue\right)_{k}|^{2}} \bigg|  \hat{x}'(y_k)  \bigg |^{2} \\Ê
%& = \mathsf{E}_{\xtrue} \{ | \hat{x}' (  y_k ) |^{2} \}
\]
where we used the law of total probability in $(a)$, 
Jensen's inequality for convex functions \cite{BoydConvexBook} in $(b)$,
and the definition \eqref{equ_new-est} in $(c)$.
\hfill  $\Box$
%% \end{proof}

%%%%%%%%%%%%%%%%%%%%%%%%%%%%%%%%%%%%%%%%
\section{Proof of Lemma \ref{equ_upper_bound_comp}} 
\label{app_proof_equ_upper_bound_comp}
%%%%%%%%%%%%%%%%%%%%%%%%%%%%%%%%%%%%%%%%

\vspace{1mm}

%% \begin{proof}[Proof of Theorem \ref{equ_upper_bound_comp}]
We wish %%% AJ 
to solve the componentwise %%% AJ 
optimization problem \eqref{equ_scalar_opt_modified}, i.e., 
$\argmin_{\hat{x}(\cdot) \ist\in\, \mathcal{U}^{k} \ist\cap\ist \mathcal{A}_{\xtrue}^{k}} \mathsf{E}_{\xtrue} \big\{ ( \hat{x} ( \mathbf{y} ) )^{2} \big\}$, for $k \notin \supp(\xtrue)$.
%% Observe first 
Note that $x_{k} \!=\! 0$ and, thus, the variance equals the mean power $\mathsf{E}_{\xtrue} \big\{ ( \hat{x} ( \mathbf{y} ) )^{2} \big\}$.

We first observe that the constraint $\hat{x} \!  \in \! \mathcal{A}_{\xtrue}^{k}$ implies that the estimator $\hat{x}$ is unbiased, 
and thus
%% i.e., 
$\mathcal{U}^{k} \rmv\cap \mathcal{A}_{\xtrue}^{k} = \mathcal{A}_{\xtrue}^{k}$. 
Indeed, using \eqref{equ_scalar_est_additive} and $x_k \! = \! 0$, we have
%% To see this, note that
\begin{align}
\mathsf{E}_{\xtrue} \{ \hat{x} ( \mathbf{y} ) \}  
& \eq \underbrace{ \mathsf{E}_{\xtrue} \{ y_k \} }_{x_k \,(=0)} \,+\; \mathsf{E}_{\xtrue} \{ \hat{x}' ( \mathbf{y}) \} \nonumber\\[0mm]
%& = x_k  + \normconstgauss \int_{\mathbf{y}} \hat{x}'(\mathbf{y})  e^{-\frac{1}{2 \sigma^{2}} \|\mathbf{y} - \mathbf{x_{0}} \|^{2}_{2}} d \mathbf{y} \\
& \eq x_k  \ist+\, \normconstgauss  \int_{\mathbb{R}^{N}} \!\hat{x}' ( \mathbf{y} ) \,
  e^{- \|\mathbf{y} - \mathbf{x} \|^{2}_{2}/(2 \sigma^{2})}  \ist d \mathbf{y} \nonumber\\[.5mm]
& \eq x_k  \ist+\, \normconstgauss  \int_{\mathbb{R}^{N-1}} \! e^{- \|\mathbf{y}_{\sim k} - \mathbf{x}_{\sim k} \|^{2}_{2}/(2 \sigma^{2})} \ist 
  \Bigg[ \underbrace{ \int_{-\infty}^\infty \!\hat{x}' ( \mathbf{y} ) \, e^{- ( y_k-0)^{2}/(2 \sigma^{2})} \ist d y_k }_0 \Bigg] d\mathbf{y}_{\sim k} \nonumber\\[-4.5mm]
& \eq  x_k
\label{equ_scalar_unbiased}
\end{align}
where $\mathbf{x}_{\sim k}$ and $\mathbf{y}_{\sim k}$ denote the $(N\!-\!1)$-dimensional vectors obtained from $\xtrue$ and $\mathbf{y}$ by removing the $k$th component $x_k$ and $y_k$, respectively, 
%% coefficient, 
and the result in \eqref{equ_scalar_unbiased}
%% last line 
follows because
$\int_{-\infty}^\infty \hat{x}' ( \mathbf{y} ) \, e^{- y_k^{2}/(2 \sigma^{2})} \ist d y_k = 0$ due to the odd symmetry assumption \eqref{equ_constr_odd_symm}.
Thus, we can replace the constraint $\hat{x}(\cdot) \rmv\in\ist \mathcal{U}^{k} \rmv \cap  \mathcal{A}_{\xtrue}^{k}$ in \eqref{equ_scalar_opt_modified} by $\hat{x}(\cdot) \! \in \! \mathcal{A}_{\xtrue}^{k}$.

A solution of \eqref{equ_scalar_opt_modified} can now be found by noting that for any  $\hat{x}(\cdot) \! \in \! \mathcal{A}_{\xtrue}^{k}$, we 
\vspace{1mm}
have
\begin{align*}
\mathsf{E}_{\xtrue} \big\{ ( \hat{x} ( \mathbf{y} ) )^{2} \big\}  %& =  \normconstgauss \int_{\mathbf{y}} | \hat{x} ( \mathbf{y} ) | e^{- \frac{1}{2 \sigma^{2}} \| \mathbf{y} - \xtrue \|^{2}} d \mathbf{y} \\Ê
& \eq  \normconstgauss \int_{\mathbb{R}^N} \!\! \big( y_k + \hat{x}'(\mathbf{y}) \big)^{2} \, e^{- \| \mathbf{y} - \xtrue \|_2^{2}/(2 \sigma^{2})} \ist d \mathbf{y}\\[1mm]Ê
& \eq \normconstgauss \int_{\mathbb{R}^N} \rmv  y_k^{2} \, e^{- \| \mathbf{y} - \xtrue \|_2^{2}/(2 \sigma^{2})} \ist d \mathbf{y} \\[.5mm]
& \rule{16mm}{0mm} + \normconstgauss \int_{\mathbb{R}^N} \!  \big[ 2 \ist y_k \ist \hat{x}'(\mathbf{y}) + ( \hat{x}'(\mathbf{y}) )^{2} \big] \, 
  e^{- \| \mathbf{y} - \xtrue \|_2^{2}/(2 \sigma^{2})} \ist d \mathbf{y} .   \nonumber  \\[-9mm]
& \nonumber
\end{align*}
The first term is equal to $\sigma^{2}+ x_{k}^{2} = \sigma^{2}$. Regarding the second term, let $\mathbf{y}_k$ be the length-$(S \rmv+\! 1)$ subvector of $\mathbf{y}$ that comprises all $y_{l}$ with $l \in \{k\} \cup \supp(\xtrue)$. Due to \eqref{equ_constr_invariance}, $\hat{x}'(\mathbf{y})$ depends only on $\mathbf{y}_k$ and can thus be written (with some abuse of notation) as $\hat{x}'(\mathbf{y}_k)$. Let $\bar{\mathbf{y}}_k$ denote the complementary subvector of $\mathbf{y}$, i.e., the length-$(N \!-\! S \!-\! 1)$ subvector comprising all $y_{l}$ with $l \not\in \{k\} \cup \supp(\xtrue)$. Furthermore, let $\mathbf{x}_k$ and $\bar{\mathbf{x}}_k$ denote the analogous subvectors of $\mathbf{x}$. The second integral can then be written as the product
\begin{align*}
&\frac{1}{(2 \pi \sigma^{2})^{(S+1)/2}} \int_{\mathbb{R}^{S+1}} \! \big[ 2 \ist y_k \ist \hat{x}'(\mathbf{y}_k) + ( \hat{x}'(\mathbf{y}_k) )^{2} \big] \, 
  e^{- \| \mathbf{y}_k - \xtrue_k \|_2^{2}/(2 \sigma^{2})} \ist d \mathbf{y}_k\\
& \rule{36mm}{0mm} \times \frac{1}{(2 \pi \sigma^{2})^{(N-S-1)/2}} \int_{\mathbb{R}^{N-S-1}} \!  e^{- \| \bar{\mathbf{y}}_k - \bar{\mathbf{x}}_k \|_2^{2}/(2 \sigma^{2})} \ist d \bar{\mathbf{y}}_k \,. \\[-9mm]
& \nonumber
\end{align*}
The second factor is 1, and thus we have
\be
 \label{equ_integral_upper_bound} 
\mathsf{E}_{\xtrue} \big\{ ( \hat{x} ( \mathbf{y} ) )^{2} \big\} \eq \sigma^{2} +\, \frac{1}{(2 \pi \sigma^{2})^{(S+1)/2}} 
  \int_{\mathbb{R}^{S+1}} \! \big[ 2 \ist y_k \ist \hat{x}'(\mathbf{y}_k) + ( \hat{x}'(\mathbf{y}_k) )^{2} \big] \, 
    e^{- \| \mathbf{y}_k - \xtrue_k \|_2^{2}/(2 \sigma^{2})} \ist d \mathbf{y}_k \,.
%% \label{equ_quad_form_scalar_power}
\ee
Using the symmetry property \eqref{equ_constr_odd_symm}, this can be written as
%\begin{align}
%\label{equ_expr_quad_form_comp}
%\mathsf{E} _{\xtrue} \{  |y_k  +\hat{x}' ( \mathbf{y} ) |^{2} \} & = \sigma^{2}  \nonumber \\Ê
% + \frac{4}{(2 \pi \sigma^{2})^{(S+1)/2}} & \int_{\mathbf{y} \in \mathcal{I}} \hat{x}'(\mathbf{y}) y_k e^{- \frac{1}{2 \sigma^{2}} | y_k|^{2}}  d y_k \prod_{l \in \supp(\xtrue)} \left( e^{- \frac{1}{2 \sigma^{2}} | y_l - x_l|^{2}} - e^{- \frac{1}{2 \sigma^{2}} | y_l + x_l|^{2}} \right) d y_l \\
%+ \frac{2}{(2 \pi \sigma^{2})^{(S+1)/2}} & \int_{\mathbf{y} \in \mathcal{I} } | \hat{x}'(\mathbf{y})|^{2} e^{- \frac{1}{2 \sigma^{2}} | y_k|^{2}} d y_k \prod_{l \in \supp(\xtrue)} \left( e^{- \frac{1}{2 \sigma^{2}} | y_l - x_l|^{2}} + e^{- \frac{1}{2 \sigma^{2}} | y_l + x_l|^{2}} \right) d y_l
%\end{align}.
%Lets rewrite \eqref{equ_expr_quad_form_comp} more compactly as:
\be
\label{equ_constr_bb_objective_1}
\mathsf{E}_{\xtrue} \big\{ ( \hat{x} ( \mathbf{y} ) )^{2} \big\} \eq \sigma^{2} +\, \frac{2}{(2 \pi \sigma^{2})^{(S+1)/2}} 
  \int_{\mathbb{R}^{S+1}_+} \! \big[ 2 \ist \hat{x}' (\mathbf{y}_k) \ist b(\mathbf{y}_k) + (\hat{x}'(\mathbf{y}_k))^{2} \ist c(\mathbf{y}_k) \big] \ist d \mathbf{y}_k \,,
\ee
with
\vspace{-2mm}
\begin{align}
b(\mathbf{y}_{k}) &\,\triangleq\, y_k \, e^{- y_k^{2}/(2 \sigma^{2})} \! \prod_{l \ist\in\ist \supp(\xtrue)} \! \big[ e^{- ( y_l - x_l)^{2}/(2 \sigma^{2})} - e^{- ( y_l + x_l)^{2}/(2 \sigma^{2})} \big] \label{equ_diff_short_b}\\
c(\mathbf{y}_{k}) &\,\triangleq\, e^{- y_k^{2}/(2 \sigma^{2})} \! \prod_{l \ist\in\ist \supp(\xtrue)} \! \big[ e^{- ( y_l - x_l)^{2}/(2 \sigma^{2})} + e^{- ( y_l + x_l)^{2}/(2 \sigma^{2})} \big] \,. \label{equ_diff_short_c} \\[-10mm]
& \nonumber
\end{align}
We sketch the derivation of expressions  \eqref{equ_diff_short_b} and \eqref{equ_diff_short_c} by showing the first of $S+1$ similar sequential calculations. For simplicity of notation and without loss of generality, we assume for this derivation that $k \rmv=\! 1$ and 
$\supp(\mathbf{x}) = \{2,\ldots,S \rmv\rmv+\! 1\}$. The integral in \eqref{equ_integral_upper_bound} then becomes 
\begin{align} 
  & \int_{\mathbb{R}^{S+1}} \! \big[ 2 \ist y_k \ist \hat{x}'(\mathbf{y}_k) + ( \hat{x}'(\mathbf{y}_k) )^{2} \big] \, 
    e^{- \| \mathbf{y}_k - \xtrue_k \|_2^{2}/(2 \sigma^{2})} \ist d \mathbf{y}_k  \nonumber \\[-.7mm]Ê
    &\hspace*{30mm} = \int_{\mathbb{R}^{S+1}} \! \big[ 2 \ist y_1 \ist \hat{x}'(\mathbf{y}_1) + ( \hat{x}'(\mathbf{y}_1) )^{2} \big] \, 
     \Bigg[ \prod_{l=1}^{S+1} e^{- ( y_l - x_l )^{2}/(2 \sigma^{2})} \Bigg] \ist d \mathbf{y}_1  \,. \label{equ_integral_upper_bound_1}Ê
\end{align} 
The $\int_{\mathbb{R}^{S+1}}$ integration can now be represented as $\int_{\mathbb{R}^{S} \times (\mathbb{R}_{+} \cup \mathbb{R}_{-})}$,
where the component $\int_{\mathbb{R}^{S}}$ refers to $y_1,\ldots,y_S$ and the component $\int_{\mathbb{R}_{+} \cup \mathbb{R}_{-}}$ refers to $y_{S+1}$.
Then \eqref{equ_integral_upper_bound_1} can be further processed as
\begin{align} 
    &\int_{\mathbb{R}^{S} \times \mathbb{R}_{+}} \!\! \big[ 2 \ist y_1 \ist \hat{x}'(\mathbf{y}_1) + ( \hat{x}'(\mathbf{y}_1) )^{2} \big] \, 
     \Bigg[ \prod_{l=1}^{S+1} e^{- ( y_l - x_l )^{2}/(2 \sigma^{2})} \Bigg]  \ist d \mathbf{y}_1  \nonumber \\Ê
    & \hspace*{30mm} + \int_{\mathbb{R}^{S} \times \mathbb{R}_{-}} \!\! \big[ 2 \ist y_1 \ist \hat{x}'(\mathbf{y}_1)+ ( \hat{x}'(\mathbf{y}_1) )^{2} \big] \, 
     \Bigg[ \prod_{l=1}^{S+1} e^{- ( y_l - x_l )^{2}/(2 \sigma^{2})} \Bigg]  \ist d \mathbf{y}_1  \nonumber \\[1mm] 
    &\stackrel{(*)}{=}   \int_{\mathbb{R}^{S} \times \mathbb{R}_{+}} \!\! \big[ 2 \ist y_1 \ist \hat{x}'(\mathbf{y}_1) \, \big(e^{- ( y_{S+1} - x_{S+1} )^{2}/(2 \sigma^{2})}
      - e^{- ( y_{S+1} + x_{S+1} )^{2}/(2 \sigma^{2})} \big) \nonumber \\[-1mm]Ê
     &\hspace*{30mm} +\ist\ist ( \hat{x}'(\mathbf{y}_1) )^{2} \, \big(e^{- ( y_{S+1} - x_{S+1} )^{2}/(2 \sigma^{2})}+ e^{- ( y_{S+1} + x_{S+1} )^{2}/(2 \sigma^{2})} \big) \big] \, 
     \Bigg[ \prod_{l=1}^{S} e^{- ( y_l - x_l )^{2}/(2 \sigma^{2})} \Bigg] \ist d \mathbf{y}_1  \nonumber
\end{align} 
where the odd symmetry property \eqref{equ_constr_odd_symm} was used in $(*)$. After performing this type of manipulation $S$ times, the integral is obtained in the 
\vspace{1mm}
form 
\begin{align} 
%%   & \int_{\mathbb{R}^{S+1}} \! \big[ 2 \ist y_1 \ist \hat{x}'(\mathbf{y}_1) + ( \hat{x}'(\mathbf{y}_1) )^{2} \big] \, 
%%     e^{- \| \mathbf{y}_1 - \xtrue_1 \|_2^{2}/(2 \sigma^{2})} \, d \mathbf{y}_1  \nonumber \\Ê
    &\int_{\mathbb{R} \times \mathbb{R}^{S}_{+}} \! \Bigg[ 2 \ist y_1 \ist \hat{x}'(\mathbf{y}_1) \prod_{l=2}^{S+1} \big(e^{- ( y_l - x_l )^{2}/(2 \sigma^{2})}- e^{- ( y_l + x_l )^{2}/(2 \sigma^{2})}\big) \nonumber \\Ê
     &\hspace*{30mm}+\ist\ist ( \hat{x}'(\mathbf{y}_1) )^{2} \prod_{l=2}^{S+1} \big(e^{- ( y_l - x_l )^{2}/(2 \sigma^{2})}+ e^{- ( y_l + x_l )^{2}/(2 \sigma^{2})}\big)  \Bigg] \ist 
     e^{- y_1^{2}/(2 \sigma^{2})} \, d \mathbf{y}_1  \nonumber \\[-9.5mm]
& \nonumber
\end{align} 
where $x_1 \!=\rmv 0$ was used. With $y_1 \ist\ist \hat{x}'(y_1,\ldots) = (-y_1) \ist\ist \hat{x}'(-y_1,\ldots)$, this becomes 
\vspace{1mm}
further 
\begin{align} 
      &\int_{\mathbb{R}_{+} \times \mathbb{R}^{S}_{+}} \! \Bigg[ 2 \ist y_1 \ist \hat{x}'(\mathbf{y}_1) \ist\ist  2 \ist \ist e^{- y_1^{2}/(2 \sigma^{2})} \prod_{l=2}^{S+1} \big(e^{- ( y_l - x_l )^{2}/(2 \sigma^{2})}- e^{- ( y_l + x_l )^{2}/(2 \sigma^{2})}\big) \nonumber \\Ê
     &\hspace*{30mm}+\ist\ist ( \hat{x}'(\mathbf{y}_1) )^{2} \ist\ist  2 \ist\ist  e^{- y_1^{2}/(2 \sigma^{2})} 
       \prod_{l=2}^{S+1} \big(e^{- ( y_l - x_l )^{2}/(2 \sigma^{2})}+ e^{- ( y_l + x_l )^{2}/(2 \sigma^{2})}\big)  \Bigg] \ist d \mathbf{y}_1 \,. \nonumberÊ \\[-10mm]
& \nonumber
\end{align} 
Finally, removing our ``notational simplicity'' assumptions $k \rmv=\! 1$ and $\supp(\mathbf{x}) = \{2,\ldots,S \rmv\rmv+\! 1\}$, this can be written for 
general $k$ and $\supp(\mathbf{x})$ 
\vspace{1mm}
as
\begin{align} 
            &2\ist e^{- y_k^{2}/(2 \sigma^{2})} \int_{\mathbb{R}^{S+1}_{+}} \! \Bigg[ 2 \ist y_k \ist \hat{x}'(\mathbf{y}_k) 
              \!\rmv\prod_{l \in \supp(\xtrue)} \!\!\! \big(e^{- ( y_l - x_l )^{2}/(2 \sigma^{2})}- e^{- ( y_l + x_l )^{2}/(2 \sigma^{2})} \big) \nonumber \\Ê
     &\hspace*{30mm} +\ist\ist ( \hat{x}'(\mathbf{y}_k) )^{2} \!\!\prod_{l \in \supp(\xtrue)} \!\!\! \big(e^{- ( y_l - x_l )^{2}/(2 \sigma^{2})}+ e^{- ( y_l + x_l )^{2}/(2 \sigma^{2})}\big) \Bigg]\ist d \mathbf{y}_k \,. \label{equ_expr_integral_upepr_bound_der}Ê
\end{align}
Inserting \eqref{equ_expr_integral_upepr_bound_der} into \eqref{equ_integral_upper_bound} yields \eqref{equ_constr_bb_objective_1}. 

The integral $\int_{\mathbb{R}^{S+1}_+} \! \big[ 2 \ist \hat{x}' (\mathbf{y}_k) \ist b(\mathbf{y}_k) + (\hat{x}'(\mathbf{y}_k))^{2} \ist c(\mathbf{y}_k) \big] \ist d \mathbf{y}_k$ is minimized with respect to $\hat{x}' (\mathbf{y}_k)$ by minimizing the integrand $2 \ist \hat{x}' (\mathbf{y}_k) \ist b(\mathbf{y}_k) + (\hat{x}'(\mathbf{y}_k))^{2} \ist c(\mathbf{y}_k)$ pointwise for each value of $\mathbf{y}_k \! \in \! \mathbb{R}_{+}^{S+1}$. This is easily done by completing the square
in
%% s with respect to the quantity 
$\hat{x}'(\mathbf{y}_k)$, yielding the optimization problem 
$\min_{\hat{x}'(\mathbf{y}_k)} \!\big[ \hat{x}'(\mathbf{y}_k)+b (\mathbf{y}_{k}) / c(\mathbf{y}_{k}) \big]^2\rmv$.
Thus, the optimal $\hat{x}'(\mathbf{y}_k)$ is obtained 
\vspace{.7mm}
as
%we obtain the following pointwise characterization of the function $\hat{x}'_{k, \xtrue}(\mathbf{y})$ that minimizes the expected power
\[
\hat{x}'_{k,\xtrue}(\mathbf{y}_k) \,\triangleq\,  -  \frac{b (\mathbf{y}_{k})}{c(\mathbf{y}_{k})} \eq - \, y_k \!\prod_{l \ist\in\ist \supp(\xtrue)} \!\!\mbox{tanh} \bigg( \frac{x_l y_l}{\sigma^2} \bigg)  
  \qquad \text{for all }\, \mathbf{y}_k \! \in \! \mathbb{R}^{S+1}_+ 
\vspace{.7mm}
\]
%%% BEGIN ALEX 08042010
and the corresponding pointwise minimum of the integrand is given by $-(b(\mathbf{y}_{k}))^{2}/ c(\mathbf{y}_{k})$. 
The extension $\hat{x}'_{k,\xtrue}(\mathbf{y})$ to all $\mathbf{y} \! \in \! \mathbb{R}^N$ is then obtained using the properties \eqref{equ_constr_odd_symm} and \eqref{equ_constr_invariance}, and the optimal component estimator solving \eqref{equ_scalar_opt_modified} follows as $\hat{x}_{k,\xtrue}(\mathbf{y}) = y_k  + \hat{x}'_{k, \xtrue}(\mathbf{y})$.
%$\mathsf{E}_{\xtrue} \{  |y_k  +\hat{x}' ( \mathbf{y} ) |^{2} \} $:
The corresponding minimum variance, denoted by $\mbox{BB}_{\text{c}}^{k}(\xtrue)$, is obtained by substituting the minimum value of the integrand, 
$-(b(\mathbf{y}_{k}))^{2}/ c(\mathbf{y}_{k})$, in \eqref{equ_constr_bb_objective_1}. This 
\pagebreak %%%%%%%%%
yields 
\be
\label{equ_expr_bbc_comp1_app}
\mbox{BB}^{k}_{\text{c}}(\xtrue) \,\triangleq\, \mathsf{E}_{\xtrue} \big\{ ( \hat{x}_{k,\xtrue}(\mathbf{y}) )^{2} \big\}  
\eq \sigma^{2} - \, \frac{2}{(2 \pi \sigma^{2})^{(S+1)/2}} 
  \int_{\mathbb{R}^{S+1}_+} \! \frac{ (b(\mathbf{y}_{k}))^{2}}{ c(\mathbf{y}_{k})}d \mathbf{y}_k \,. 
\vspace{1mm}
\ee
Inserting \eqref{equ_diff_short_b} and \eqref{equ_diff_short_c} into \eqref{equ_expr_bbc_comp1_app} and simplifying gives \eqref{equ_expr_bbc_comp1}.
%%% END ALEX 08042010
%% \end{proof}

%% \pagebreak %%%%%%%%%%

%%%%%%%%%%%%%%%%%%%%%%%%%%%%%%%%%%%%%%%%
\section{Proof of Equation \eqref{equ_exp_conv_bb}} 
\label{app_proof_exponential}
%%%%%%%%%%%%%%%%%%%%%%%%%%%%%%%%%%%%%%%%

\vspace{1mm}

%% \begin{proof}[Proof of Equation \eqref{equ_exp_conv_bb}]
To show 
%% Equation 
\eqref{equ_exp_conv_bb}, we consider
$g(x;\sigma^{2})$ for $x \geq 0$ (this is sufficient
%% valid 
since $g(-x;\sigma^{2}) = g(x;\sigma^{2})$), and we use
%% start from 
the simple bound $\mbox{tanh}(x) \ist\geq\ist 1 - e^{-x}$, which can be verified using elementary calculus.
We then obtain from \eqref{equ_expr_bbc_g}, for 
\vspace{1.5mm}
$x \geq 0$,
%%% BEGIN ALEX 19042010
\begin{align*}
g(x; \sigma^{2}) %& =  \int_{y \in \mathbb{R}_{+}}   \frac{1}{\sqrt{2 \pi \sigma^{2}} } e^{-\frac{1}{2 \sigma^{2}} (x^2 + y^2 ) } \mbox{sinh} \bigg( \frac{x y}{\sigma^{2}} \bigg) \mbox{tanh} \bigg( \frac{x y}{\sigma^{2}} \bigg) dy \\Ê
& \,\geq\,\frac{1}{\sqrt{2 \pi \sigma^{2}} } \int_{0}^\infty \! e^{-(x^2 + y^2)/(2 \sigma^{2}) } \,\mbox{sinh} \bigg( \frac{x y}{\sigma^{2}} \bigg) \big(1 - e^{- x y/\sigma^{2}} \big) \ist\ist dy \\[2.5mm]
& \eq \frac{1}{\sqrt{2 \pi \sigma^{2}} } \int_{0}^\infty \! \big[ e^{-(x-y)^{2}/(2 \sigma^{2})} -\ist e^{-(x+y)^{2}/(2 \sigma^{2})} \big] \big(1 -\ist e^{- x y/\sigma^{2}} \big) \ist\ist dy\\[2.5mm]
& \eq \frac{1}{\sqrt{2 \pi \sigma^{2}} } \int_{0}^\infty \! \big[ 
    e^{-(x-y)^{2}/(2 \sigma^{2})} 
    \ist-\ist\ist  e^{-(x^{2}+y^{2})/(2 \sigma^{2})} 
    \ist-\ist\ist e^{-(x+y)^{2}/(2 \sigma^{2})} 
    \ist+\ist\ist  e^{-(x+y)^{2}/(2 \sigma^{2})} e^{- x y/\sigma^{2}}
    \big] \ist\ist dy \\[2.5mm]
%% & \eq \frac{1}{\sqrt{2 \pi \sigma^{2}} } \int_{0}^\infty \! \big[ e^{-(x-y)^{2}/(2 \sigma^{2})} \ist+\ist\ist  e^{-(x+y)^{2}/(2 \sigma^{2})} e^{- x y/\sigma^{2}}
%%     \!-\ist\ist e^{-(x+y)^{2}/(2 \sigma^{2})} \ist-\ist\ist  e^{-(x^{2}+y^{2})/(2 \sigma^{2})} \big] \ist\ist dy \\[2.5mm]
& \,\geq\, \frac{1}{\sqrt{2 \pi \sigma^{2}} } \int_{0}^\infty \! \big[ 
    e^{-(x-y)^{2}/(2 \sigma^{2})}
    \ist-\ist\ist  e^{-(x^{2}+y^{2})/(2 \sigma^{2})} 
    \ist-\ist\ist e^{-(x+y)^{2}/(2 \sigma^{2})} 
    \big] \ist\ist dy \\[2.5mm]
%& \,\stackrel{(1)}{=}\, \hspace*{-43mm} \underbrace{\frac{1}{\sqrt{2 \pi \sigma^{2}} } \int_{0}^\infty \! e^{-(x-y)^{2}/(2 \sigma^{2})} \ist\ist dy}_{ %% \displaystyle
%              \hspace*{43mm} \rule{0mm}{3mm} 1 \;-\; \frac{1}{\sqrt{2 \pi \sigma^{2}} } \int_{-\infty}^0 \rmv e^{-(x-y)^{2}/(2 \sigma^{2})} \ist\ist dy 
%               \;\eq\; 1 \;-\; \frac{1}{\sqrt{2 \pi \sigma^{2}} } \int_{0}^\infty \! e^{-(x+y)^{2}/(2 \sigma^{2})} \ist\ist dy } 
%     \hspace*{-42mm} -\;\, \frac{1}{\sqrt{2 \pi \sigma^{2}} } \int_{0}^\infty \! \big[\ist\ist  e^{-(x^{2}+y^{2})/(2 \sigma^{2})} \ist+\ist\ist e^{-(x+y)^{2}/(2 \sigma^{2})} 
%    \big] \ist\ist dy \\[1.5mm]    
& \,=\, {\frac{1}{\sqrt{2 \pi \sigma^{2}} } \int_{0}^\infty \! e^{-(x-y)^{2}/(2 \sigma^{2})} \ist\ist dy}%_{ %% \displaystyle
          %    \hspace*{43mm} \rule{0mm}{3mm} 1 \;-\; \frac{1}{\sqrt{2 \pi \sigma^{2}} } \int_{-\infty}^0 \rmv e^{-(x-y)^{2}/(2 \sigma^{2})} \ist\ist dy 
             %  \;\eq\; 1 \;-\; \frac{1}{\sqrt{2 \pi \sigma^{2}} } \int_{0}^\infty \! e^{-(x+y)^{2}/(2 \sigma^{2})} \ist\ist dy } 
      \,-\, \frac{1}{\sqrt{2 \pi \sigma^{2}} } \int_{0}^\infty \! \big[\ist\ist  e^{-(x^{2}+y^{2})/(2 \sigma^{2})} \ist+\ist\ist e^{-(x+y)^{2}/(2 \sigma^{2})} 
    \big] \ist\ist dy \,. \\[-8mm]
& \nonumber  
\end{align*}
The first integral can be written as ${\frac{1}{\sqrt{2 \pi \sigma^{2}} } \int_{0}^\infty \! e^{-(x-y)^{2}/(2 \sigma^{2})} \ist\ist dy} = 1 - \frac{1}{\sqrt{2 \pi \sigma^{2}} } \int_{-\infty}^0 \rmv e^{-(x-y)^{2}/(2 \sigma^{2})} \ist\ist dy  = 1 - \frac{1}{\sqrt{2 \pi \sigma^{2}} } \int_{0}^\infty \! e^{-(x+y)^{2}/(2 \sigma^{2})} \ist\ist dy$. The bound thus 
\vspace{1mm}
becomes
\begin{align*}
g(x; \sigma^{2}) & \,\geq\, 1-\ist\ist \frac{1}{\sqrt{2 \pi \sigma^{2}} } \int_{0}^\infty \! \big[   2 \ist e^{-(x+y)^{2}/(2 \sigma^{2})} +  e^{-(x^{2}+y^{2})/(2 \sigma^{2})}\big] \ist\ist dy    \, \\[2.5mm]
%% & \eq 1-\ist\ist \frac{1}{\sqrt{2 \pi \sigma^{2}} } \int_{0}^\infty \! \big[   2 \ist  e^{-(x^{2}+2xy+y^{2})/(2 \sigma^{2})} +  e^{-(x^{2}+y^{2})/(2 \sigma^{2})} \big] \ist\ist dy  \, \\[2.5mm]
& \eq 1-\ist\ist \frac{1}{\sqrt{2 \pi \sigma^{2}} } \int_{0}^\infty \! \big[ 2 \ist  e^{-2xy/(2 \sigma^{2})} \,\ist+\, 1 \big] \, e^{-(x^{2}+y^{2})/(2 \sigma^{2})} 
 \ist\ist dy  \, \\[1.5mm]
%& \eq \frac{1}{\sqrt{2 \pi \sigma^{2}} } \int_{0}^\infty \! ??? \, dy  \quad \text{***intermediate step!!}\\[2.5mm]
& \,\stackrel{(*)}{\geq}\, 1-\ist\ist \frac{1}{\sqrt{2 \pi \sigma^{2}} } \int_{0}^\infty \!   3 \ist\ist e^{-(x^{2}+y^{2})/(2 \sigma^{2})}  \ist\ist dy  \, \\[2.5mm]
& \eq 1 -\ist\ist \frac{3}{\sqrt{2 \pi \sigma^{2}}} \, e^{- x^{2}/(2 \sigma^{2})} \!\int_{0}^\infty \! e^{- y^2/(2 \sigma^{2}) } \, dy \\[3mm]
& \eq  1 -\ist\ist \frac{3}{2}  \, e^{- x^{2}/(2 \sigma^{2})} 
%% \\
%%                             & \eq       1 - \left( \frac{2}{\sqrt{2 \pi \sigma^{2}}} + \frac{1}{2} \right) e^{- \frac{x^{2}}{2 \sigma^{2}}}
\end{align*}
%% that is, 
%% \[
%% g(x; \sigma^{2}) \,\geq\, 1 - \frac{3}{2}  \, e^{- x^{2}/(2 \sigma^{2})}.
%% \]
%% \pagebreak %%%%%%%%%
where $e^{-2xy/(2 \sigma^{2})} \leq 1$ was used in $(*)$. This bound on $g(x;\sigma^{2})$ is actually valid for all $x \rmv\in\rmv \mathbb{R}$ because 
$g(-x;\sigma^{2}) = g(x;\sigma^{2})$. Inserting it in \eqref{equ_expr_bbc_comp1}, we obtain
\be
\label{equ_bound_bb_c_k}
\mbox{BB}_{\text{c}}^{k}(\xtrue) \,\leq\, \bigg[ 1 -\! \prod_{l \ist\in\ist \supp(\xtrue) } \!\! \bigg( 1 - \frac{3}{2} \, e^{- x_l^{2}/(2 \sigma^{2})} \bigg) \bigg] \ist \sigma^{2}  \,.
\ee
%% where $c \triangleq  \left( \frac{2}{\sqrt{2 \pi \sigma^{2}}} + \frac{1}{2} \right)$.
%% From \eqref{equ_bound_bb_c_k} and 
The statement in \eqref{equ_exp_conv_bb} follows since we have (note that $\sum_{\mathcal{I} \ist\subseteq\ist \supp(\xtrue) }$ denotes the sum over all possible subsets $\mathcal{I}$ of $\supp(\xtrue)$, including $\supp(\xtrue)$ and the empty set 
\vspace{1mm}
$\emptyset$)
\begin{align*}
1 -\! \prod_{l \ist\in\ist \supp(\xtrue) } \!\! \bigg( 1 - \frac{3}{2} \, e^{- x_l^{2}/(2 \sigma^{2})} \bigg) 
&\eq 1 -\! \sum_{\mathcal{I} \ist\subseteq\ist \supp(\xtrue) } \, \prod_{l \ist\in\ist \mathcal{I} } \bigg( \!\! -\rmv \frac{3}{2} \, e^{- x_l^{2}/(2 \sigma^{2})} \bigg) \\[.5mm]
&\eq  -\! \sum_{ \mathcal{I} \ist\subseteq\ist \supp(\xtrue) ,\, \mathcal{I} \ist\neq\ist \emptyset } \,\ist \prod_{l \ist\in\ist \mathcal{I} } \bigg( \!\! -\rmv \frac{3}{2} \, e^{- x_l^{2}/(2 \sigma^{2})} \bigg) \\[.5mm] 
&\,\leq\, \sum_{ \mathcal{I} \ist\subseteq\ist \supp(\xtrue) ,\, \mathcal{I} \ist\neq\ist \emptyset } \,\ist \prod_{l \ist\in\ist \mathcal{I} } \bigg( \frac{3}{2} \, e^{- x_l^{2}/(2 \sigma^{2})} \bigg) \\[.5mm] 
&\,\leq\, \sum_{ \mathcal{I} \ist\subseteq\ist \supp(\xtrue) ,\, \mathcal{I} \ist\neq\ist \emptyset } \,\ist \prod_{l \ist\in\ist \mathcal{I} } \bigg( \frac{3}{2} \, e^{- \xi^{2}/(2 \sigma^{2})} \bigg) \\[.5mm] 
&\eq \sum_{ \mathcal{I} \ist\subseteq\ist \supp(\xtrue) ,\, \mathcal{I} \ist\neq\ist \emptyset } \! \bigg( \frac{3}{2} \, e^{- \xi^{2}/(2 \sigma^{2})} \bigg)^{\!|\mathcal{I}|} \\[.5mm] 
&\,\leq\, \sum_{ \mathcal{I} \ist\subseteq\ist \supp(\xtrue) ,\, \mathcal{I} \ist\neq\ist \emptyset } \! \bigg( \frac{3}{2} \bigg)^{\!\!S} e^{- \xi^{2}/(2 \sigma^{2})} \\[2mm]
& \,\leq\, 2^{S} \bigg( \frac{3}{2} \bigg)^{\!\!S} e^{- \xi^{2}/(2 \sigma^{2})}\\[1mm]
&\eq 3^{S} \ist e^{- \xi^{2}/(2 \sigma^{2})}
\end{align*}
where 
%% in the last step 
we have used the fact that the number of different subsets $\mathcal{I} \subseteq \supp(\mathbf{x})$ is $2^{| \supp(\mathbf{x})|} \! = \!2^{S}$.
%%% END ALEX 19042010
Inserting the last bound in \eqref{equ_bound_bb_c_k} and, in turn, the resulting bound on $\mbox{BB}_{\text{c}}^{k}(\xtrue)$ in \eqref{equ_upper_bound_bbc} yields \eqref{equ_exp_conv_bb}.

%%%%%%%%%%%%%%%%%%%%%%%%%%%%%%%%%%%%%%%%
\section{MSE of the ML Estimator}
\label{app_MSE_ML}
%%%%%%%%%%%%%%%%%%%%%%%%%%%%%%%%%%%%%%%%

\vspace{1mm}

%ALEX: Here too I would replace integrals with expectations whenever possible.

We calculate the MSE $\varepsilon(\xtrue; \ML_est)$ of the ML estimator $\ML_est$ 
%% as defined 
in \eqref{equ_LS_ML_coincide}. Let $\hat{x}_{\text{ML},k}$ denote the $k$th component of $\ML_est$.
We
\vspace{-1mm}
have
\begin{align}
\varepsilon(\xtrue; \ML_est) %& = \mathsf{E}_{\xtrue}Ê\{  \| \xtrue - \ML_est \|^{2}_{2} \} \\
&\eq \sum_{k=1}^{N} \mathsf{E}_{\xtrue} \big\{ ( \hat{x}_{\text{ML},k} \rmv-\rmv x_k )^{2} \big\} \nonumber \\[.5mm]
&\eq \sum_{k=1}^{N} \big[ \mathsf{E}_{\xtrue}  \big\{  \hat{x}^{2}_{\text{ML},k}  \big\} - 2 \, \mathsf{E}_{\xtrue} \big\{  \hat{x}_{\text{ML},k} \big \} \ist\ist x_k   +x^{2}_k \ist \big]   \nonumber \\[.5mm]
&\eq \sum_{k=1}^{N} \big[ \mathsf{E}_{\xtrue} \{ \hat{x}_{\text{ML},k}^{2} \} 
  \ist +\ist \big( \mathsf{E}_{\xtrue} \{ \hat{x}_{\text{ML},k}  \} \rmv- x_k \big)^{2} 
  -\ist \big( \mathsf{E}_{\xtrue} \{  \hat{x}_{\text{ML},k}  \}  \big)^{2} \big] \,.
\label{equ_expr_ml_risk}
\end{align}
Thus, we have to calculate the quantities $\mathsf{E}_{\xtrue} \{ \hat{x}_{\text{ML},k} \}$ and $\mathsf{E}_{\xtrue} \{ \hat{x}_{\text{ML},k}^{2} \}$.
%We sketch only the calculation of $\mathsf{E}_{\xtrue} \{  \left( \ML_est \right)_{k}  \}$ since the calculation of $\mathsf{E}_{\xtrue} \{ \big| \left( \ML_est \right)_{k} \big|^{2} \}$ is very similar.

%% \pagebreak %%%%%%%%%%

We recall that $\hat{x}_{\text{ML},k}(\mathbf{y}) = \big( {\mathsf P}_{\!S} (\mathbf{y}) \big)_k$, where $\mathsf{P}_{\! S}$ is an operator that retains the $S$ largest (in magnitude) components and zeros out all others. Let $\mathcal{L}_{k}$ denotes the set of vectors $\mathbf{y}$ for which $y_k$ is not among the $S$ largest (in magnitude) components. We then 
\vspace{-.5mm}
have
\[
\hat{x}_{\text{ML},k}(\mathbf{y}) \eq\! \begin{cases}
  y_k \,, & \mathbf{y} \!\not\in\! \mathcal{L}_{k}\\[-2mm]
  0 \,, & \mathbf{y} \!\in\! \mathcal{L}_{k} \,.
  \end{cases}
\]
Equivalently, 
$\hat{x}_{\text{ML},k}(\mathbf{y}) = y_k \ist [1 \rmv-\text{I}( \mathbf{y} \!\in\! \mathcal{L}_{k} )]$, where $\text{I}( \mathbf{y} \!\in\! \mathcal{L}_{k} )$ is the indicator function of the event $\{ \mathbf{y} \!\in\! \mathcal{L}_{k} \}$ (i.e., $\text{I}( \mathbf{y} \!\in\! \mathcal{L}_{k} )$ is $1$ if $\mathbf{y} \!\in\! \mathcal{L}_{k}$ and $0$ else).
Thus, we obtain $\mathsf{E}_{\xtrue} \{ \hat{x}_{\text{ML},k} \}$ as
%%% BEGIN ALEX 15042010
\begin{align}
\mathsf{E}_{\xtrue} \{ \hat{x}_{\text{ML},k} \} 
&\eq \mathsf{E}_{\xtrue} \big\{  y_k \ist [1 \rmv-\text{I}( \mathbf{y} \!\in\! \mathcal{L}_{k} )] \big\} \nonumber\\Ê
&\eq x_k - \mathsf{E}_{\xtrue} \big\{ y_{k} \ist\ist \text{I}( \mathbf{y} \!\in\! \mathcal{L}_{k} ) \big\} \nonumber\\[.5mm]Ê
&\stackrel{(a)}{\eq} x_k - \mathsf{E}_{\xtrue}^{(y_{k})} \big\{ \mathsf{E}_{\xtrue}^{(\mathbf{y}_{\sim k})} \big\{  y_{k} \ist\ist \text{I}( \mathbf{y} \!\in\! \mathcal{L}_{k} ) \big|\ist y_{k}  \big\} \rmv\rmv\big\} \nonumber\\[.5mm]Ê
&\stackrel{(b)}{\eq} x_k - \mathsf{E}_{\xtrue}^{(y_{k})} \big\{     y_{k} \, \mathsf{E}_{\xtrue}^{(\mathbf{y}_{\sim k})}  \big\{  \ist\ist \text{I}( \mathbf{y} \!\in\! \mathcal{L}_{k} ) \big|\ist y_{k}  \big\} \rmv\rmv\big\} \nonumber\\[.5mm]Ê
&\eq x_k -\mathsf{E}_{\xtrue}^{(y_{k})} \big\{  y_{k}  \, \text{P}_{\xtrue} ( \mathbf{y} \!\in\! \mathcal{L}_{k} |\ist y_k ) \big\}                                                                                               
\label{equ_expr_ml_risk_mean_ml}
\end{align}
where the notations $\mathsf{E}_{\xtrue}^{(y_{k})}$ and $\mathsf{E}_{\xtrue}^{(\mathbf{y}_{\sim k})}$ indicate that the expectation is taken with respect to the random quantities $y_{k}$ and $\mathbf{y}_{\sim k}$, respectively (here, $\mathbf{y}_{\sim k}$ denotes $\mathbf{y}$ without the component $y_{k}$) and $\text{P}_{\xtrue} ( \mathbf{y} \!\in\! \mathcal{L}_{k} |\ist y_k )$ is the conditional probability that $\mathbf{y} \!\in\! \mathcal{L}_{k}$, given $y_k$.
%\begin{align}
%\mathsf{E}_{\xtrue} \big\{  y_{k} \ist\ist \text{I}( \mathbf{y} \!\in\! \mathcal{L}_{k} ) \big|\ist y_{k}  \big\}  & = \int_{\mathbf{y} \in \mathbb{R}^{N}} f(\mathbf{y}|y_{k};\xtrue) y_{k} \text{I}( \mathbf{y} \!\in\! \mathcal{L}_{k} ) d \mathbf{y}  =  \int_{\mathbf{y} \in \mathbb{R}^{N}} f(\mathbf{y}|y_{k};\xtrue) y_{k} \text{I}( \mathbf{y} \!\in\! \mathcal{L}_{k} ) d \mathbf{y} \\
%& = y_{k} \int_{\mathbf{y} \in \mathbb{R}^{N}} f(\mathbf{y}|y_{k};\xtrue) \text{I}( \mathbf{y} \!\in\! \mathcal{L}_{k} ) d \mathbf{y} = y_{k}\text{P}_{\xtrue} ( \mathbf{y} \!\in\! \mathcal{L}_{k} |\ist y_k ). 
%\end{align}
Furthermore, we used the law of total probability in $(a)$ and the fact that 
%% the value of 
$y_{k}$ is held constant in the conditional expectation $ \mathsf{E}_{\xtrue} \big\{  y_{k} \ist\ist \text{I}( \mathbf{y} \!\in\! \mathcal{L}_{k} ) \big|\ist y_{k}  \big\}$ in $(b)$. 
%%% END ALEX 15042010
Similarily,
\vspace{-4mm}
\begin{align}
\mathsf{E}_{\xtrue} \{ \hat{x}_{\text{ML},k}^{2} \} 
&\eq \mathsf{E}_{\xtrue} \big\{  y_k^2 \, [1 \rmv- \text{I}( \mathbf{y} \!\in\! \mathcal{L}_{k} ) ]^2 \big\} \nonumber\\Ê
&\eq \mathsf{E}_{\xtrue} \big\{  y_k^2 \, [1 \rmv- \text{I}( \mathbf{y} \!\in\! \mathcal{L}_{k} ) ] \big\} \nonumber\\Ê
&\eq \sigma^{2}  + x_k^2  - \mathsf{E}_{\xtrue} \big\{  y_k^2 \, \text{I}( \mathbf{y} \!\in\! \mathcal{L}_{k} ) \big\} \nonumber\\Ê
&\eq  \sigma^{2}+ x_k^2  - \mathsf{E}_{\xtrue}^{(y_{k})} \big\{  y_k^2 \, \text{P}_{\xtrue} ( \mathbf{y} \!\in\! \mathcal{L}_{k} |\ist y_k ) \big\} \,.
 \label{equ_expr_ml_risk_power_ml}
\end{align}
Calculating $\mathsf{E}_{\xtrue} \{ \hat{x}_{\text{ML},k} \}$ and $\mathsf{E}_{\xtrue} \{ \hat{x}_{\text{ML},k}^{2} \}$ is thus reduced to calculating the conditional probability $\text{P}_{\xtrue} ( \mathbf{y} \!\in\! \mathcal{L}_{k} |\ist y_k )$. 

%% \pagebreak %%%%%%%%%%

Let $\mathcal{M}_k \triangleq \{1,\ldots,N\} \rmv\setminus\rmv \{k\}$, and let $\mathcal{P}$ denote the set of all binary partitions $( \mathcal{A} , \mathcal{B} )$ of the set $\mathcal{M}_k$, where $\mathcal{A}$ is at least of cardinality $S$: 
%% the following set of ordered pairs $(\mathcal{A} , \mathcal{B})$ with $\mathcal{A} , \mathcal{B} \subseteq \mathcal{M}_k$: 
%% %% the set $\mathcal{P}$ which is defined as a certain collection of ordered pairs where the first and second coefficient are subsets of $\mathcal{M}_k \triangleq \{1,2...,N\} \setminus \{k\}$:
\[
\mathcal{P} \,\triangleq\, \big\{ ( \mathcal{A} , \mathcal{B} ) \ist \big|\, \mathcal{A} \!\subseteq\! \mathcal{M}_k, \ist \mathcal{B} \!\subseteq\! \mathcal{M}_k, \ist \mathcal{A} \rmv\cap\rmv \mathcal{B} \rmv=\rmv \emptyset, \ist \mathcal{A} \rmv\cup\rmv \mathcal{B} \rmv=\rmv \mathcal{M}_k, \ist |\mathcal{A}| \rmv\geq\rmv S \big\} \,.
\]
%% That is, we consider
%%% BEGIN ALEX 19042010
In order to evaluate the conditional probability $\text{P}_{\xtrue} ( \mathbf{y} \!\in\! \mathcal{L}_{k} |\ist y_k )$ 
of the event $\{ \mathbf{y} \!\in\! \mathcal{L}_{k} \}$, i.e., of the event that a given 
$y_{k}$ is not among the $S$ largest (in magnitude) components of $\mathbf{y}$, we split the event $\{ \mathbf{y} \!\in\! \mathcal{L}_{k} \}$ into 
several elementary events.
More specifically, let $\mathcal{E}_{\rmv\mathcal{A},\mathcal{B}}$ denote the event that every component $y_{l}$ with $l \!\in\! \mathcal{A}$ satisfies $|y_l| > |y_k|$ and every component  $y_{l}$ with $l \!\in\! \mathcal{B}$ satisfies $|y_l| \leq |y_k|$.
The events $\mathcal{E}_{\rmv\mathcal{A},\mathcal{B}}$ for all $(\mathcal{A},\mathcal{B}) \!\in\! \mathcal{P}$ 
%% for varying $\mathcal{A}$ and $\mathcal{B}$ 
are mutually exclusive, i.e., 
$(\mathcal{A},\mathcal{B}) \neq (\mathcal{A}',\mathcal{B}') \,\Rightarrow\, \mathcal{E}_{\rmv\mathcal{A},\mathcal{B}} \cap \mathcal{E}_{\rmv\mathcal{A}',\mathcal{B}'}  =  \emptyset$,
and their union corresponds to
%% they are exhaustive conditioned on 
the event $\{ \mathbf{y} \! \in \! \mathcal{L}_{k} \}$, i.e., $\bigcup_{(\mathcal{A}, \mathcal{B}) \ist\in\ist \mathcal{P}} \mathcal{E}_{\rmv\mathcal{A},\mathcal{B}} \ist=\ist \{ \mathbf{y} \! \in \! \mathcal{L}_{k} \}$.
%% for all $(\mathcal{A},\mathcal{B}), (\mathcal{A}',\mathcal{B}') \in \mathcal{P}$.
Consequently, 
%% we have according to the law of total probability that
\begin{align}
%\mathbf{P}_{\xtrue} ( \mathbf{y} \in \mathcal{L}_{k} | y_k = y) & = \sum\limits_{(\mathcal{A}, \mathcal{B}) \in \mathcal{P}} \prod_{l \in \mathcal{A} \cap \supp( \xtrue) } \mathbf{P}_{\xtrue} ( | y_l | > |y|) \prod_{l \in \mathcal{B} \cap \supp( \xtrue) } \mathbf{P}_{\xtrue} ( | y_l | < |y|) \nonumber \\Ê
%&  \prod_{l \in \mathcal{A} \setminus \supp( \xtrue) } 2 \mathcal{Q}(|y|/\sigma^{2})  \prod_{l \in \mathcal{B} \setminus \supp( \xtrue) }  (1 - 2 \mathcal{Q}(|y|/\sigma^{2}))  \label{equ_prob_expansion_ml_risk}
\text{P}_{\xtrue} ( \mathbf{y} \!\in\! \mathcal{L}_{k} |\ist y_k\rmv\rmv=\rmv y) 
&\eq\rmv \sum\limits_{(\mathcal{A}, \mathcal{B}) \ist\in\ist \mathcal{P}} \text{P}_{\xtrue}  (  \mathcal{E}_{\rmv\mathcal{A},\mathcal{B}} \ist|\ist y_k\rmv\rmv=\rmv y ) \nonumber\\[1mm] 
&\eq\rmv \sum\limits_{(\mathcal{A}, \mathcal{B}) \ist\in\ist \mathcal{P}} \, \prod_{l \ist\in\ist \mathcal{A}} \text{P}_{\xtrue}  \big( |y_{l}| \rmv>\rmv |y_{k}| \ist\big|\ist\ist y_k\rmv\rmv=\rmv y \big) 
  \prod_{m \ist\in\ist \mathcal{B}} \! \text{P}_{\xtrue}  \big( |y_{m}| \rmv\leq\rmv |y_{k}| \ist\big| \ist\ist y_k\rmv\rmv=\rmv y \big) \nonumber\\[1mm] 
&\eq\rmv \sum\limits_{(\mathcal{A}, \mathcal{B}) \ist\in\ist \mathcal{P}}  \, \prod_{l \ist\in\ist \mathcal{A}} \text{P}_{\xtrue}  ( |y_{l}| \rmv>\rmv |y| ) 
  \prod_{m \ist\in\ist \mathcal{B}} \!\text{P}_{\xtrue}  ( |y_{m}| \rmv\leq\rmv |y|  ) \nonumber\\[1mm] 
&\eq\rmv \sum\limits_{(\mathcal{A}, \mathcal{B}) \ist\in\ist \mathcal{P}} \, \prod_{l \ist\in\ist \mathcal{A}\,\cap\, \supp( \xtrue) } \!\!\text{P}_{\xtrue}  ( |y_{l}| \rmv>\rmv |y| ) 
  \!\!\prod_{m \ist\in\ist \mathcal{B}\,\cap\, \supp( \xtrue) } \!\!\text{P}_{\xtrue}  ( |y_{m}| \rmv\leq\rmv |y|  ) \nonumber \\[.5mm] 
& \hspace*{25mm} \times \!\! \prod_{n \ist\in\ist \mathcal{A}\ist\setminus \ist \supp( \xtrue) } \!\!\text{P}_{\xtrue}  ( |y_{n}| \rmv>\rmv |y| ) 
  \!\!\prod_{p \ist\in\ist \mathcal{B}\ist\setminus \ist \supp( \xtrue) } \!\!\text{P}_{\xtrue}  ( |y_{p}| \rmv\leq\rmv |y|  )  \nonumber\\[2mm]
&\eq\rmv \sum\limits_{(\mathcal{A}, \mathcal{B}) \ist\in\ist \mathcal{P}} \, \prod_{l \ist\in\ist \mathcal{A} \,\cap\, \supp( \xtrue) } 
  \bigg[ Q\bigg( \frac{|y| \!-\rmv x_l}{\sigma} \bigg) + 1 - Q\bigg( \frac{-|y| \!-\rmv x_l}{\sigma} \bigg) \bigg] \nonumber \\ 
&\hspace*{25mm} \times\!\! \prod_{m \ist\in\ist \mathcal{B} \,\cap\, \supp( \xtrue) } \bigg[ \!\rmv- Q\bigg( \frac{|y| \!-\rmv x_m}{\sigma} \bigg) + Q\bigg( \frac{-|y| \!-\rmv x_m}{\sigma} \bigg) \bigg] \nonumber \\Ê
&\hspace*{25mm} \times\!\! \prod_{n \ist\in\ist \mathcal{A} \ist\setminus\ist \supp( \xtrue) } \!\! 2 \ist\ist Q\bigg( \frac{|y|}{\sigma} \bigg) 
  \!\prod_{p \ist\in\ist \mathcal{B} \ist\setminus\ist \supp( \xtrue) } \bigg[ 1 - 2 \ist\ist Q\bigg( \frac{|y|}{\sigma} \bigg) \bigg] 
%%% END ALEX 19042010
\label{equ_prob_expansion_ml_risk} \\[-8.5mm]Ê
&\nonumber
\end{align}
where we have used the fact that the $y_l$ 
%% of the observed vector $\mathbf{y}$ 
are independent and $k \!\notin\! \mathcal{M}_{k}$; furthermore, $Q(y) \triangleq \frac{1}{\sqrt{2 \pi}} \int_{y}^{\infty} e^{-x^{2}/2} \ist dx$ is the right tail probability of a standard Gaussian random variable.
Plugging \eqref{equ_prob_expansion_ml_risk} into \eqref{equ_expr_ml_risk_mean_ml} and \eqref{equ_expr_ml_risk_power_ml}  and, in turn, the resulting expressions into
\eqref{equ_expr_ml_risk} yields a (very complicated) 
%% closed-form 
expression of $\varepsilon(\xtrue; \ML_est)$. This expression is evaluated numerically
in Section \ref{sec_sim}.

%%%%%%%%%%%%%%%%%%%%%%%%%%%%%%%%%%%%%%%%
\bibliographystyle{IEEEtran}
\bibliography{literatur_journal_SLM,literatur_Asilomar} %%IEEEabrv,
%%%%%%%%%%%%%%%%%%%%%%%%%%%%%%%%%%%%%%%%

\end{document}